\documentclass[graybox]{svmult}

\usepackage{subcaption}
\usepackage{type1cm}        
\usepackage{makeidx}         
\usepackage{graphicx}        
\usepackage{multicol}        
\usepackage[bottom]{footmisc}

\usepackage{newtxtext}       %
\usepackage[varvw]{newtxmath}       

\makeindex             

\usepackage{tikz}
\usepackage{pgfplots}
\usetikzlibrary{shapes.arrows}
\usetikzlibrary{arrows.meta}
\usetikzlibrary{calc}
\usetikzlibrary{hobby,decorations.markings}
\usetikzlibrary {perspective}
\usepackage{pgfplots}
\pgfplotsset{compat=1.18}
\usepackage{mathtools}
\usepackage{enumitem}
\usepackage[onelanguage,linesnumbered,boxed]{algorithm2e}
\SetNlSkip{0em} 
\SetNlSty{}{\hspace{1em}}{:}

\SetAlgoNlRelativeSize{-3} 
\renewcommand{\theAlgoLine}{\hspace{1em}\ifnum\value{AlgoLine}<10\noindent{\bf 0}\fi\textbf{\arabic{AlgoLine}}}
\usepackage{alltt}
\usepackage{multicol}
\definecolor{darkcyan}{rgb}{0.0, 0.55, 0.55}

\usepackage[markup=underlined]{changes} 
\definechangesauthor[name= ads, color=red]{ads}
\definechangesauthor[name= vD, color=blue]{vD}

\usepackage[mathlines,pagewise]{lineno} 
\begin{document}

\title*{Practical Introduction to FEM with GMSH: A MATLAB/Octave Perspective}
\author{Víctor Domínguez\orcidID{0000-0002-6095-619X} and\\
Alejandro Duque-Salazar\orcidID{0009-0000-1339-5281} }

\institute{V. Domínguez \at  Dep. Estadística, Informática y Matemáticas, Universidad Pública de Navarra, Avda Tarazona s/n, 31500 Spain. \email{victor.dominguez@unavarra.es} \and
A. Duque-Salazar \at Basque Center of Applied Mathematics, Mazarredo Zumarkalea, 14, Abando, 48009 Bilbo, Bizkaia, Spain, \email{aduque@bcamath.org}}

\maketitle

\abstract{The Finite Element Method (FEM) is a powerful computational tool for solving partial differential equations (PDEs). Although commercial and open-source FEM software packages are widely available,  an independent implementation of FEM provides significant educational value, provides a deeper understanding of the method, and enables the development of custom solutions tailored to specialized applications or integration with other solvers.
This work introduces a 3D $\mathbb{P}_m$-element FEM implementation in MATLAB/Octave that is designed to balance educational clarity with computational efficiency. A key feature is its integration with GMSH, an open-source 3D mesh generator with CAD capabilities that streamlines mesh generation for complex geometries. By leveraging GMSH data structures, we provide a seamless connection between geometric modeling and numerical simulation. The implementation focuses on solving the general convection-diffusion-advection equation and serves as a flexible foundation for addressing advanced problems, including elasticity, mixed formulations, and integration with other numerical methods.}

\section{Introduction}

The Finite Element Method (FEM) is a powerful computational technique widely used for the numerical solution of partial differential equations (PDEs), particularly for elliptic problems. Originally developed to address complex elasticity and structural modeling, FEM has become popular in many scientific and engineering disciplines. Starting from a variational formulation—where the solution minimizes an energy functional and/or satisfies specific integral test with sufficiently smooth functions—FEM is defined by restricting the test to a finite-dimensional subspace made from piecewise polynomial functions defined on a partition of the original domain. This partition is commonly referred to as the mesh in the finite element method. The error of the FEM solution is then proven, by Cea's Lemma-like results, to be of the same order as the best approximation in this discrete space. These fundamentals of FEM are extensively covered in advanced graduate textbooks, including classical works as \cite{Ciarlet1978} as well as more recent contributions detailed in \cite{brenner_scott,MR2050138,MR4242224, MR4269305}.

We emphasize that although commercial (e.g., ANSYS, Abaqus) and free (e.g., FreeFEM \cite{MR3043640}, FEniCS \cite{AlnaesEtal2015}) FEM software options are available,  implementing the Finite Element Method (FEM) offers significant advantages. These include educational benefits, such as gaining a comprehensive understanding of the operation of the method, its strengths and limitations, and the parameters that influence its behavior. Additionally, it allows for the creation of customized solutions for specialized applications, such as FEM-BEM coupling (cf. \cite{doi:https://doi.org/10.1002/9781119176817.ecm2008} and references therein).

As noted in introductory FEM implementation materials, the primary challenge is obtaining a mesh for the problem domain. While theoretically achievable for complex geometries, generating a mesh can be difficult, even for 2D cases. However, once a mesh is available—regardless of its simplicity or the limited geometric information it provides—the FEM can be implemented with just a few lines of code.

The choice of MATLAB as the programming language is deliberate for several compelling reasons. It is straightforward yet powerful, primarily when the code is written with the proper characteristics of this environment in mind. MATLAB offers a wide range of matrix-like operations, including structures for sparse matrix manipulation and efficiently implemented built-in functions, such as direct or iterative solvers. In addition, MATLAB's robust graphical capabilities facilitate the visualization of meshes and solutions, making it an ideal tool for displaying meshes, computing FE solutions, and many other applications. MATLAB's support for advanced programming features, including object-oriented programming (OOP), further enhances its flexibility, allowing programmers to incrementally increase the complexity of their code, starting with the essentials of the algorithms.

For example, we cite the pioneering work of \cite{MR1709562}, which demonstrates that a $\mathbb{P}_1$ (linear) 2D FEM implementation for mixed boundary conditions of the Laplacian can be achieved in just 50 lines of MATLAB code. Works such as \cite{ MR2194203, MR1935965,funken2011fem,  MR3674245} extend these goals, but with two additional features: writing more optimal code by using vectorization, which results in faster execution, and extending the implementation to more sophisticated FEM algorithms. 

On the other hand, GNU Octave is another viable alternative to MATLAB. Most MATLAB code can be executed in Octave with no or minimal modification. Being open-source, Octave is developed, improved, and maintained by a community of developers and users contributing to its source code at no cost. This makes Octave an especially attractive option for those seeking a cost-effective solution with broad community support.

This work has two primary objectives: first, to develop a 3D $\mathbb{P}_m$-element FEM implementation in MATLAB/Octave that highlights algorithmic principles for educational purposes; and second, to ensure that the implementation is efficient enough to solve efficiently the general reaction-diffusion-advection equation as a first step towards more general problems (e.g., elasticity, evolution problems and mixed formulations), more sophisticated elements or getting its integrations with other numerical schemes. In essence, our efforts align with the goals of the works referenced in \cite{MR1709562,  MR2194203, MR1935965, funken2011fem, MR3674245}, as well as the aspirations of Pancho’s team, led by our dearly missed colleague and friend Francisco-Javier Sayas \cite{sayas2015fem, sayas_fem_tools_3d} (see \cite{teampancho} for a comprehensive overview of his legacy).
A significant novelty of this work is addressing the challenge of mesh construction by integrating our solver with GMSH cf. \cite{geuzaine2009GMSH}, an open-source 3D FEM mesh generator equipped with a built-in CAD engine. GMSH streamlines the handling of complex geometries, facilitating the solution of realistic and challenging problems.

GMSH is a powerful tool renowned for its advanced mesh generation and geometric modeling capabilities. Although it is well-documented with extensive resources, its sophisticated features and broad functionality can present a steep learning curve for neophytes. However, the effort invested in mastering GMSH is often rewarded, as it provides precise and efficient tools for handling complex simulations. As a secondary goal, this work also aims to introduce GMSH, with a particular focus on the data structures used for mesh management and how this data structure can be enhanced to facilitate our FEM implementation.

This paper is structured as follows: In \S2, we present the model problem considered for numerical solution, as well as its weak formulation—the first step toward applying FEM. In \S3, we provide a sufficiently detailed introduction to 3D classical continuous Lagrange finite elements on tetrahedra. Then, \S4 describes how GMSH defines and processes domains and meshes, and how this information can be read and handled within our own code. In \S5, we revisit FEM and outline how the method reduces to the computation of a linear system of equations. We discuss the corresponding algorithms in Section 7 and briefly highlight the key features of our software package. Specifically, we demonstrate how mesh data obtained from GMSH can be efficiently preprocessed and how the method is implemented in a computationally efficient way—both tasks accomplished in a loop-free manner, which is essential for achieving this efficiency objective. Finally, in \S8, we provide numerical experiments that demonstrate the capabilities of our code.

\section{The (not-so) toy 3D elliptic boundary value problem}

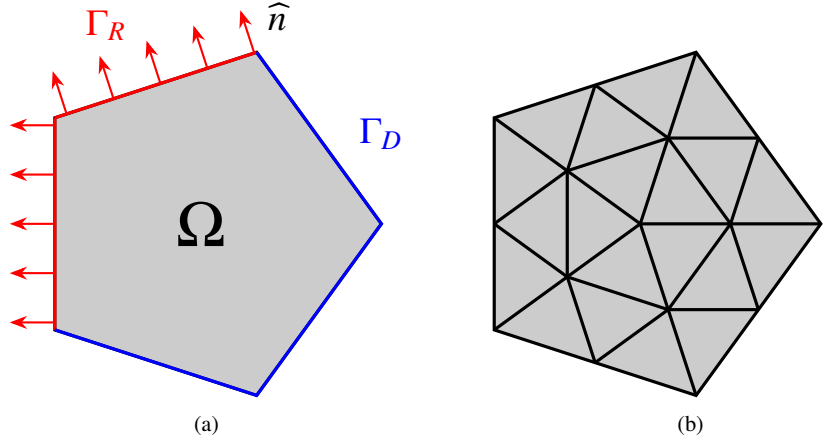
\begin{figure}[!htb]
    \centering
    \begin{subfigure}{0.45\textwidth}
       \begin{tikzpicture}[very thick, use Hobby shortcut,scale=0.6]
  \def\n{5} 
  \def\r{4} 
  \def\m{5} 
  \def\l{1} 

  \draw[fill=gray!40] (360/\n:\r)
    \foreach \x in {2,...,\n} {
      -- (360/\n*\x:\r)
    } -- cycle;
\draw[red] (360/\n*2:\r) -- (360/\n*3:\r);
\draw[red] (360/\n*1:\r) -- (360/\n*2:\r);
\draw[blue] (360/\n*3:\r) -- (360/\n*4:\r)-- (360/\n*5:\r)-- (360/\n*6:\r);

\node  at (-2.1,4.4) {\Large \color{red}$\Gamma_R$};
\node[blue]  at (4,2) {\Large $\Gamma_D$};
\node  at (1.7,4.6) {\Large $\widehat{n}$};

  \path[decoration={
    markings,
    mark=between positions 0.18 and 0.4 step 0.05 with { \draw[thick, -Stealth] (0,0) -- (0,-0.6);},
  }, postaction={decorate},red] (0:2) foreach \x in {1,...,4} { -- (360/5*\x:\r) } -- cycle;

  \path[decoration={
    markings,
    mark=between positions 0.4 and 0.6 step 0.05 with { \draw[thick, -Stealth] (0,0) -- (0,-0.6);},
  }, postaction={decorate},red] (0:2) foreach \x in {1,...,4} { -- (360/5*\x:\r) } -- cycle;

  \node  at (0,0) {\Huge $\Omega$};
\end{tikzpicture}
    \caption{\label{fig: Domain3d}}
    \end{subfigure}\hfill
    \begin{subfigure}{0.45\textwidth}
        \begin{tikzpicture}[very thick, use Hobby shortcut,  scale=0.6]
  \def\n{5} 
  \def\r{4} 
  \def\m{5} 
  \def\l{1} 

  \draw[fill=gray!40] (360/\n:\r)
    \foreach \x in {2,...,\n} {
      -- (360/\n*\x:\r)
    } -- cycle;

\foreach \x in {1,...,\n}{
    \draw (360/\n*\x:\r) -- (0,0);
}
\foreach \x in {1,...,\n}{
    \pgfmathsetmacro{\next}{mod(\x,\n)+1}
    \draw ($(360/\n*\x:\r)!0.5!(360/\n*\next:\r)$) -- ($(360/\n*\x:\r)!0.5!(0,0)$);
    \draw ($(360/\n*\next:\r)!0.5!(0,0)$) -- ($(360/\n*\x:\r)!0.5!(360/\n*\next:\r)$);
}
\draw (360/\n:2)
    \foreach \x in {2,...,\n} {
      -- (360/\n*\x:2)
    } -- cycle;

\end{tikzpicture}
    \caption{\label{fig: DiscDomain}}
    \end{subfigure}\\
\caption{Left: (2D) Sketch of the domain $\Omega$. Right: example of conformal simplicial mesh of $\Omega$. }
\end{figure}

Let $\Omega \subset \mathbb{R}^3$ be a compact polygonal domain with boundary  $\Gamma = \Gamma_D  \cup \Gamma_R$,   $\Gamma_D,\ \Gamma_R$  disjoint, and  outward unit normal vector $\widehat{\boldsymbol{n}}$ (see Figure \ref{fig: Domain3d}).
We consider the following  boundary value problem: for given functions $f:\Omega\to\mathbb{R}$,
$u_D:\Gamma_D\to \mathbb{R}$ (Dirichlet data) and
$g_R:\Gamma_R\to\mathbb{R}$ (Robin data), find $u:\Omega\to\mathbb{R}$ such that
\begin{equation}\label{eq:FemProblem}
\left\{\begin{array}{rcll}
-\nabla\cdot \left(\underline{\boldsymbol{\kappa}} (\nabla u)^\top \right) +   \boldsymbol{\beta} \cdot \nabla u + c u &=& f,  &\text{in }\Omega,\\
u &=& u_{D}, &\text{on }\Gamma_{D},\\
\left(\underline{\boldsymbol{\kappa}} (\nabla u)^\top \right)\cdot\widehat{\boldsymbol{n}}+ \alpha u &=& g_R, & \text{on } \Gamma_R.
\end{array}\right.
\end{equation}
Here, $\underline{\boldsymbol{\kappa}}$ is a $3\times 3$   matrix function (the diffusion matrix), $\boldsymbol{\beta}$ is a vector function (the velocity field) and  $ c$ and $\alpha$ are scalar functions (the reaction term and the Robin coefficient respectively). Notice that the gradient $\nabla $ is taken as a {\em row} operator and that the scalar product ``$\cdot$'' is defined in a standard manner for vectors, regardless of whether they are represented as column or row vectors. 
 
Problem  \eqref{eq:FemProblem} is sufficiently general to include diffusion and advection equations, for which the matrix $\kappa$ is symmetric and positive definite, Helmholtz equations, and, with minor modifications, linear elasticity problem.

The variational formulation for  \eqref{eq:FemProblem} is obtained by multiplying by sufficiently smooth functions and applying Stokes Theorem. Indeed, let
\[
\begin{aligned}
H^1(\Omega) & = \{u\in L^2(\Omega) \ : \ \nabla u \in L^2(\Omega) \}
\end{aligned}
\]
be the Sobolev space of order one. Then any solution of  \eqref{eq:FemProblem} solves
\begin{equation}\label{eq:weak}
\begin{cases}
\text{Find }u\in H^1(\Omega) , \quad \text{with }u= u_D,\quad \text{on } {\Gamma_{D}}\\
\begin{aligned}
a(u,w) &= \ell(w), \quad \forall w\in H^1(\Omega) \ \text{with }w|_{\Gamma_D}=0
\end{aligned}
\end{cases}
\end{equation}
where
\begin{equation}
\begin{aligned}
a(u,w) &:= \iiint_{\Omega} \left[\nabla w \cdot (\underline{\boldsymbol{\kappa}}(\nabla u)^\top)   + w(\boldsymbol{\beta}\cdot \nabla u)+   cw u\right]+\iint_{\Gamma_{R}} \alpha  w u,\\
\ell(w) & := \iiint_{\Omega} fw+\iint_{\Gamma_{R}} g_Rw.
\end{aligned}
\end{equation}

The Finite Element Method (FEM), in a nutshell, begins with   \eqref{eq:weak} by replacing the continuous infinite-dimensional space $H^1(\Omega)$
 with a discrete, finite-dimensional function space. This approach, where the trial space for the solution essentially coincides with the test space, is known as a Galerkin scheme, unlike Petrov-Galerkin schemes, where trial and test spaces are distinct.  

\section{Finite element space}

Finite element spaces typically comprise simple piecewise functions, often polynomials, defined over a suitable domain partition, referred to as the elements. This partition, known as the mesh, will be described using an approach similar to how GMSH defines and manages it.

This exposition focuses on tetrahedral meshes, which partition the domain into tetrahedra. While GMSH also supports other types of meshes—such as rectangular hexahedra, prisms, or pyramids—that can be useful in specific contexts (e.g., layered, structured geometries or wave equations), tetrahedral meshes generally provide greater flexibility for arbitrary domains. GMSH handles these alternative mesh types in a similar way, allowing readers to extend the principles discussed here to other mesh structures with relative ease. 

 \subsection{Tetrahedral meshes}

We denote
  \[
\mathcal{T}_h=\{K_j\}_{j=1}^{\tt nTtrh}
,\ \text{with } K_\ell\cap K_m = \emptyset \ \text{ for } m\ne \ell \quad \mathrm{s.t.}\quad
\overline{\Omega} = \bigcup_{\overline{K}_j\in \mathcal{T}_h }\overline{K}_j,
\]
a partition of the original domain $\Omega$. Here $K_j$ is an (open) tetrahedron (3-simplex) for the discretization process, which forms the basis of the $\mathbb{P}_m$-finite element method. The number of tetrahedra in the resulting mesh is denoted henceforth by {\tt nTtrh}.

Additional assumptions will be made for ${\cal T}_h$. Specifically, each pair of different tetrahedra $\overline{K}_{\ell}$ and $\overline{K}_m$ with $\ell \neq m$ must either share a vertex, one complete edge, one complete face or have an empty intersection. In other words, we only consider {\em conformal} meshes (see Figure \ref{fig: DiscDomain}). We also assume that the mesh is non-degenerate, which means that for a given collection of meshes $\{{\cal T}_h\}_{h}$,
\[
\max_{h}\max_{K\in{\cal T}_h} \frac{h_K}{\rho_K}<\gamma
\]
where $h_K$ and $\rho_K$ are the diameter and inradius of $K$. Constant $\gamma$ is the so-called chunkiness parameter  of the mesh \cite{brenner_scott}, so that non-large values of this parameter ensure that {\em regular} (non-elongated or flat) tetrahedra are being used.

The parameter $h$ in the notation ${\cal T}_h$ is weakly related to the size of the mesh, so we write $h \to 0$ for a collection of meshes $\{{\cal T}_h\}_{h}$ to mean that the maximum of the diameters of $K_j \in \{{\cal T}_h\}$ tends to zero.

Finally, if an element has one face lying on the boundary $\Gamma$, that face is entirely contained in either $\Gamma_D$ or $\Gamma_R$. Although not necessary, it is often assumed that an element cannot have two or more faces on the boundary.  The case of non-polygonal and smooth domains $\Omega$ can also be considered using this method. In this case, an additional error occurs due to the approximation of the curved boundary.

\subsection{The $\mathbb{P}_m-$finite element space: General setting}\label{sec: FEMspace}

The finite element space on $\mathcal{P}_h^m$, continuous in $\Omega$, is defined as:
\begin{equation*}
\mathcal{P}_h^m :=\{v_h\in \mathcal{C}(\overline{\Omega}) : v_h|_{K} \in \mathbb{P}_m, \, \forall K\in \mathcal{T}_h \},
\end{equation*}
where $\mathbb{P}_m$ is the space of (trivariate) polynomials of degree less than or equal to $m$. Hence, in this space, we consider continuous functions that are polynomials on each element in the mesh. This restriction is enough to ensure that $\mathcal{P}_h^m\subset H^1(\Omega)$.

 \begin{example}
For simplicity, let us consider linear (first-degree) finite elements, i.e., $m=1$, and denote  
\[
\{\boldsymbol{x}_r^K\}_{r=0}^3
\]
the set of vertices of an element $K \subset \mathcal{T}_h$. Next, consider the set of all different {\tt nNodes} vertices of the tetrahedra in the mesh ${\cal T}_h$ appropriately indexed: 
\[
\{\boldsymbol{x}_i\}_{i=1}^{\tt nNodes} = \bigcup_{K\in{\cal T}_h} \{\boldsymbol{x}_r^K\}_{r=0}^3.
\]
 Then (i) a linear function, polynomial of degree $1$, in ${x, y, z}$ is uniquely determined by its value at four non co-planar points, and (ii) if it vanishes at three non co-linear points, respectively two different points, it also vanishes in the plane (and then the triangle) they determine resp. in the straight line (and then the edge) they define. Property (i) ensures that a polynomial function can be constructed on any $K$ provided the values at the vertices are known. Property (ii) guarantees that this piecewise polynomial construction, performed element by element, results in a continuous function across the mesh (since vertices are obviously not co-linear) and, therefore, an element of $\mathcal{P}_h^1$.

These conditions essentially constitute the classical definition of finite elements (see \cite{Ciarlet1978}). Indeed, condition (i) defines a finite element: a domain $K$, a finite-dimensional space $\mathbb{P}_1$, and a local interpolant—a way to determine $\mathbb{P}_1$; in this case, interpolation at the vertices. Condition (ii) ensures that the global interpolant, constructed on a mesh by interpolating locally on each element is continuous, and its image, $\mathcal{P}_h^1$, is a continuous finite element space.

On the other hand, notice that we can construct
\[
\{\varphi_j\}_{j=1}^{\tt nNodes}\subset \mathcal{P}_h^1, \quad {\text{with}\ }\quad
\varphi_j(\boldsymbol{x}_k) = \delta_{jk}=\begin{cases}
1, & j=k, \\
0, & \text{otherwise}.
\end{cases}
\]
The function $\varphi_j$ has its support on the union of elements that share the vertex $\boldsymbol{x}_j$. The set $\{\varphi_j\}_{j=1}^{\tt nNodes}$ constitutes the Lagrange basis for the global interpolant, ensuring that for any $v_h \in \mathcal{P}_h^1$, we have
\[
v_h = \sum_{j=1}^{\tt nNodes} v_h(\boldsymbol{x}_j)\varphi_j.
\]
\phantom{a}\hfill$\Box$
  \end{example}

 As we will now show, higher-degree elements are defined by extending this construction. For any $K$ and $m=1,2,\ldots,$ we consider a set of nodes
\[
\{\boldsymbol{x}_r^{K}\}_{r = 0}^{{\rm dofK}-1}, \quad  {\tt dofK}:= \binom{m+3}{3}.
\]
satisfying the following assumptions:
\begin{enumerate}[label = {(3.2.\Alph*)}, ref = {(3.2.\Alph*)}, leftmargin=*, align=left]
    \item Any polynomial of degree $m$ is uniquely determined by its values at $\{\boldsymbol{x}_r^{K}\}_{r = 0}^{{\rm dofK}-1}$. \label{CA}
    \item For any $\boldsymbol{x}_j, \boldsymbol{x}_k \in \{\boldsymbol{x}_i\}_{i = 1}^{{\rm nNodes}}$, there exists a unique element $\varphi_j \in \mathcal{P}_h^m$ such that \label{CB}
    \[
    \varphi_j(\boldsymbol{x}_k) = \delta_{jk}.
    \]
\end{enumerate}

 It is evident  under such hypothesis that
\[
\mathop{\rm supp} \varphi_j = \bigcup_{\boldsymbol{x}_j\in \overline{K}} \overline{K}
\]

We then distinguish two different indices for  nodes: for $\boldsymbol{x}_j =\boldsymbol{x}_r^K$, we have $j$,  the global index of the mesh $j\in\{1,\ldots,{\tt nNodes}\}$, and $r\in \{0,\ldots,{\tt dofK}-1\}$, the local index as node in $K$.  Similarly, we have for any element $K$ the local basis of $\mathbb{P}_m$
\[
N^K_r:=\varphi_j|_{K} \quad \text{s.t.}\quad \boldsymbol{x}_j =\boldsymbol{x}_r^K,
\]
which satisfies
\[
    N_r^K(\boldsymbol{x}_s^K)=\delta_{rs} \qquad \forall r,s\in \left\{0, \ldots, {\tt dofK}-1  \right\}.
\]
 Conversely,
 \[
 \varphi_j (\boldsymbol{x}) = \begin{cases}
N^K_r (\boldsymbol{x}), & \text{if }\boldsymbol{x}\in K \text{ and } \boldsymbol{x}_r^K =\boldsymbol{x}_j, \\
0,& \text{otherwise}.
 \end{cases}
 \]

\subsection{Local polynomial interpolation} \label{sub:section:3.3}
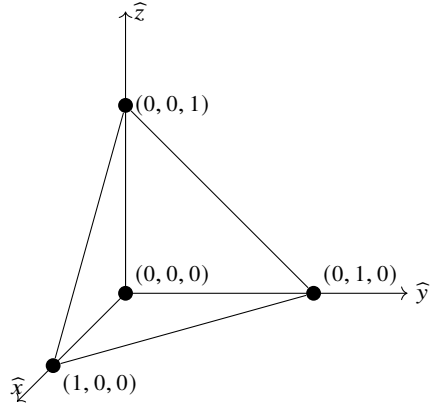
\begin{figure}[!htb]
    \centering
    \begin{tikzpicture}[scale=2.5]
  \node[circle,fill=black,scale=0.7] (A) at (0,0,0) {} ;
  \node[circle,fill=black,scale=0.7] (B) at (0,1,0) {};
  \node[circle,fill=black,scale=0.7] (C) at (1,0,0){};
  \node[circle,fill=black,scale=0.7] (D) at (0,0,1){};
    \node[above right] at (A) {$(0,0,0)$};
  \node[right] at (B) {$(0,0,1)$};
  \node[above right] at (C) {$(0,1,0)$};
  \node[below right] at (D) {$(1,0,0)$};

  \draw (A) -- (B) -- (C) -- (D) -- (B)--cycle;

  \draw[->] (0,0,0) -- (1.5,0,0) node[right] {$\widehat{y}$};
  \draw[->] (0,0,0) -- (0,1.5,0) node[right] {$\widehat{z}$};
  \draw[->] (0,0,0) -- (0,0,1.5) node[above] {$\widehat{x}$};

\end{tikzpicture}
    \caption{Reference tetrahedral element.}
    \label{fig:refElm}
\end{figure}

  $\mathbb{P}_m-$Lagrange finite elements (cf. \cite{Ciarlet1978}) are typically defined using barycentric coordinates or, equivalently, by defining a reference element and applying an affine mapping.
For the sake of simplicity, let us introduce the  {\bf reference element}
\[
\widehat{K} = \{\boldsymbol{\widehat{x}} = (\widehat{x},\widehat{y}, \widehat{z})\ \lvert\ \,0< 1-\widehat{x}-\widehat{y}-\widehat{z},\  0< \widehat{x} ,\ 0< \widehat{y},\ 0< \widehat{z} \}.
\]
(see Figure \ref{fig:refElm} for a sketch)  and the associated barycentric coordinates $\widehat{\boldsymbol{x}} = (\widehat{x},\,\widehat{y},\, \widehat{z})$ in the reference element by $\widehat{\boldsymbol\lambda} (\widehat{\boldsymbol{x}})=(\widehat{\lambda}_0 (\widehat{\boldsymbol{x}}),\,\widehat{\lambda}_1(\widehat{\boldsymbol{x}}),\,\widehat{\lambda}_2(\widehat{\boldsymbol{x}}), \,\widehat{\lambda}_3(\widehat{\boldsymbol{x}}))$ given by
\begin{equation}\label{eq:lambdahat}
\begin{aligned}
    \widehat{\lambda}_1(\boldsymbol{x})&= 1-\widehat{x}-\widehat{y}-\widehat{z},\ \quad &
    \widehat{\lambda}_2(\boldsymbol{x}) &= \widehat{x},\\
    \widehat{\lambda}_3(\boldsymbol{x}) &= \widehat{y},\ \quad &
    \widehat{\lambda}_4(\boldsymbol{x}) &= \widehat{z}.
\end{aligned}
\end{equation}

The vertices  of this element are given and indexed as follows:
\begin{equation}\label{eq: verticesTetrahedron}
\begin{aligned}
\widehat{\boldsymbol{x}}_0 & = (0,0,0),\ \quad &
\widehat{\boldsymbol{x}}_1 & = (1,0,0),\\
\widehat{\boldsymbol{x}}_2 & = (0,1,0),\ \quad &
\widehat{\boldsymbol{x}}_3 & = (0,0,1).
\end{aligned}
\end{equation}

We define the set of nodes for the interpolation problem on $\mathbb{P}_m$ in $\widehat{K}$ as
\begin{equation}\label{eq:barycentricNodes}
\{\widehat{\boldsymbol{x}}_r\}_{r=0}^{{\tt dofK}-1}=
\left\{
\left(\frac{i}m,\frac{j}{m},\frac{k}m\right) \ \Big|\
\text{ $i,j,k\in\mathbb{N}\cup\{0\}$ },\ i+j+k\leq m\right\},
\end{equation}
for $r= 0, \ldots, {\tt dofK}-1$. 

The sorting of nodes, specifically the local indexing of nodes (index $r$) within the reference tetrahedron, can undoubtedly be a subject for discussion. We adhere to the convention adopted in GMSH  displayed in Figure~\ref{fig: LagrangeNodes} for
$m=1,2,3,$ and $4$:

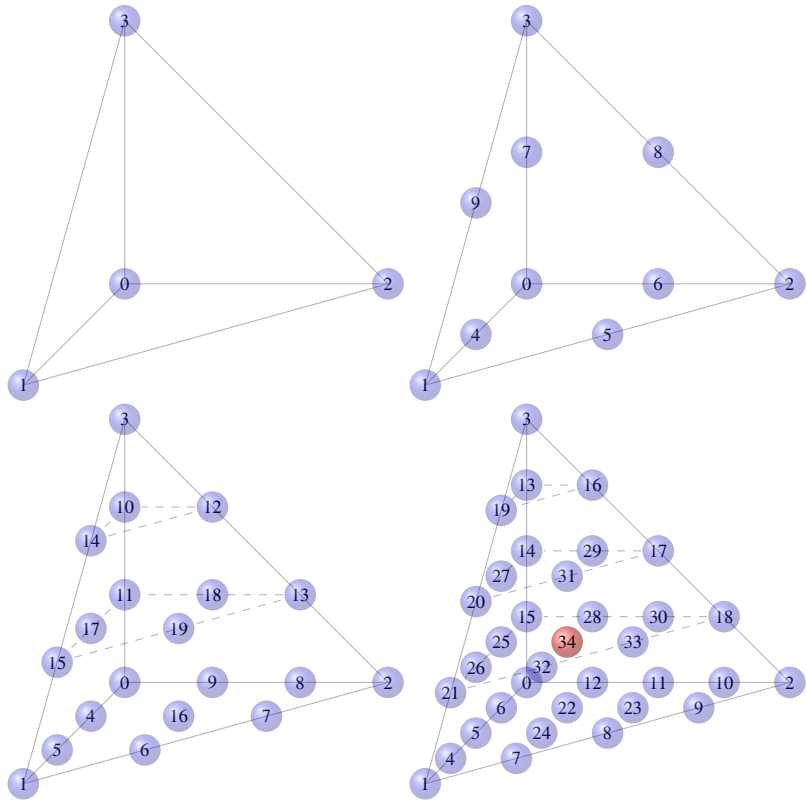
\begin{figure}
    \centering
    \begin{tikzpicture}[scale=3.5]

\coordinate (x0) at (0,0,0);
\coordinate (x2) at (1,0,0);
\coordinate (x3) at (0,1,0);
\coordinate (x1) at (0,0,1);

\draw[opacity = 0.3] (x0) -- (x1) -- (x2) -- cycle;
\draw[ opacity=0.3] (x0) -- (x3) -- (x1);
\draw[opacity = 0.3] (x3) -- (x2);

\foreach \i/\a/\b/\c/\d in {
    0/1/0/0/0,
    1/0/1/0/0,
    2/0/0/1/0,
    3/0/0/0/1}{
        \node[fill=none] at
     ($\a*(x0) + \b*(x1) + \c*(x2) + \d*(x3)$) {\scriptsize \i};
    \shade[ball color=blue, opacity=0.3]
        ($\a*(x0) + \b*(x1) + \c*(x2) + \d*(x3)$)
        circle[radius=0.06];
}
\end{tikzpicture}
    \begin{tikzpicture}[scale=3.5]

\coordinate (x0) at (0,0,0);
\coordinate (x2) at (1,0,0);
\coordinate (x3) at (0,1,0);
\coordinate (x1) at (0,0,1);

\draw[opacity = 0.3] (x0) -- (x1) -- (x2) -- cycle;
\draw[ opacity=0.3] (x0) -- (x3) -- (x1);
\draw[opacity = 0.3] (x3) -- (x2);

\foreach \i/\a/\b/\c/\d in {
    0/1/0/0/0,
    1/0/1/0/0,
    2/0/0/1/0,
    3/0/0/0/1,
    4/0.5      /0.5           /0            /0,
    5/0      /0.5      /0.5            /0,
    6/0.5            /0      /0.5            /0,
    7/0.5            /0            /0      /0.5,
    8/0            /0      /0.5      /0.5,
    9/0      /0.5            /0      /0.5}{
        \node[fill=none] at
     ($\a*(x0) + \b*(x1) + \c*(x2) + \d*(x3)$) {\scriptsize \i};
    \shade[ball color=blue, opacity=0.3]
        ($\a*(x0) + \b*(x1) + \c*(x2) + \d*(x3)$)
        circle[radius=0.06];
}

\end{tikzpicture}

    \begin{tikzpicture}[scale=3.5]

\coordinate (x0) at (0,0,0);
\coordinate (x2) at (1,0,0);
\coordinate (x3) at (0,1,0);
\coordinate (x1) at (0,0,1);

\draw[opacity = 0.3] (x0) -- (x1) -- (x2) -- cycle;
\draw[ opacity=0.3] (x0) -- (x3) -- (x1);
\draw[opacity = 0.3] (x3) -- (x2);

\foreach \i/\a/\b/\c/\d in {
    0/1/0/0/0,
    1/0/1/0/0,
    2/0/0/1/0,
    3/0/0/0/1,
    4/0.66667      /0.33333            /0            /0,
    5/0.33333      /0.66667            /0            /0,
    6/0      /0.66667      /0.33333            /0,
    7/0      /0.33333      /0.66667            /0,
    8/0.33333            /0      /0.66667            /0,
    9/0.66667            /0      /0.33333            /0,
    10/0.33333            /0            /0      /0.66667,
    11/0.66667            /0            /0      /0.33333,
    12/0            /0      /0.33333      /0.66667,
    13/0            /0      /0.66667      /0.33333,
    14/0      /0.33333            /0      /0.66667,
    15/0      /0.66667            /0      /0.33333,
    16/0.33333      /0.33333      /0.33333            /0,
    17/0.33333      /0.33333            /0      /0.33333,
    18/0.33333            /0      /0.33333      /0.33333,
    19/0      /0.33333      /0.33333      /0.33333}{
        \node[fill=none] (node\i) at
     ($\a*(x0) + \b*(x1) + \c*(x2) + \d*(x3)$) {\scriptsize \i};
    \shade[ball color=blue, opacity=0.3]
        ($\a*(x0) + \b*(x1) + \c*(x2) + \d*(x3)$)
        circle[radius=0.06];
}

  \draw[opacity = 0.3,dashed] (node10) -- (node14) -- (node12) -- (node10);
  \draw[opacity = 0.3,dashed] (node11) -- (node15) -- (node13) -- (node11);
\end{tikzpicture}
    \begin{tikzpicture}[scale=3.5]

\coordinate (x0) at (0,0,0);
\coordinate (x2) at (1,0,0);
\coordinate (x3) at (0,1,0);
\coordinate (x1) at (0,0,1);

\draw[opacity = 0.3] (x0) -- (x1) -- (x2) -- cycle;
\draw[ opacity=0.3] (x0) -- (x3) -- (x1);
\draw[opacity = 0.3] (x3) -- (x2);

\foreach \i/\a/\b/\c/\d in {
    0/1/0/0/0,
    1/0/1/0/0,
    2/0/0/1/0,
    3/0/0/0/1,
    4/0.25      /0.75            /0            /0,
    5/0.5      /0.5            /0            /0,
    6/0.75      /0.25            /0            /0,
    7/0      /0.75      /0.25            /0,
    8/0      /0.5      /0.5            /0,
    9/0      /0.25      /0.75            /0,
    12/0.25            /0      /0.25            /0,
    11/0.5            /0      /0.5            /0,
    10/0.75            /0      /0.75            /0,
    13/0.25            /0            /0      /0.75,
    14/0.5            /0            /0      /0.5,
    15/0.75            /0            /0      /0.25,
    16/0            /0      /0.25      /0.75,
    17/0            /0      /0.5      /0.5,
    18/0            /0      /0.75      /0.25,
    19/0      /0.25          /0.  /0.75     ,
    20/0      /0.5         /0.    /0.5      ,
    21/0      /0.75     /0.      /0.25,
    22/0.5      /0.25   /0.25               /0,
    24/0.25      /0.5   /0.25               /0,
    23/0.25      /0.25   /0.5               /0,
    25/0.5      /0.25    /0        /0.25        ,
    26/ 0.25      /0.5    /0        /0.25        ,
    27/ 0.25      /0.25   /0        /0.5,
    28/ 0.5  /0          /0.25       /0.25         ,
    29/ 0.25 /0          /0.25       /0.5         ,
    30/ 0.25 /0          /0.5        /0.25,
    31/   0/0.25           /0.25      /0.5      ,
    32/ 0/0.5            /0.25      /0.25        ,
    33/ 0/0.25           /0.5       /0.25
    }{
        \node[fill=none]  (node\i) at
     ($\a*(x0) + \b*(x1) + \c*(x2) + \d*(x3)$) {\scriptsize \i};
 \shade[ball color=blue, opacity=0.3] ($\a*(x0) + \b*(x1) + \c*(x2) + \d*(x3)$) circle[radius=0.06]
;
}

  \draw[opacity = 0.3,dashed] (node13) -- (node16) -- (node19) -- (node13);
  \draw[opacity = 0.3,dashed] (node14) -- (node20) -- (node17) -- (node14);
  \draw[opacity = 0.3,dashed] (node15) -- (node21) -- (node18) -- (node15);
  \shade[ball color=red, opacity=0.5] ($0.25*(x0) + 0.25*(x1) + 0.25*(x2) + 0.25*(x3)$)  circle[radius=0.06];
\node[fill=none] at ($0.25*(x0) + 0.25*(x1) + 0.25*(x2) + 0.25*(x3)$)  {\scriptsize 34};
\end{tikzpicture}

    \caption{Node indexing on the $\mathbb{P}_1$ (top-left), $\mathbb{P}_2$ (top-right), $\mathbb{P}_3$ (bottom-left) and $\mathbb{P}_4$ (bottom-right) reference elements. The first inner node appears with $\mathbb{P}_4$ and it is marked in red}
    \label{fig: LagrangeNodes}

\end{figure}

\begin{itemize}
 \item The nodes $\widehat{\boldsymbol{x}}_i$ for $i=0,1,2, 3$, are indexed according to the numeration of the vertices in   \eqref{eq: verticesTetrahedron}.
\item The nodes at the interior of the edges are indexed in the following order: first, the nodes along the edge from vertex $\widehat{\boldsymbol{x}}_0$ to vertex $\widehat{\boldsymbol{x}}_1$; then, the nodes along the edge from vertex $\widehat{\boldsymbol{x}}_1$ to vertex $\widehat{\boldsymbol{x}}_2$; next, the nodes along the edge from vertex $\widehat{\boldsymbol{x}}_2$ to vertex $\widehat{\boldsymbol{x}}_0$. The remaining nodes are indexed starting from the vertex $\widehat{\boldsymbol{x}}_3$ to all other vertices, starting from the edge shared with vertex $\widehat{\boldsymbol{x}}_0$, followed by the edge shared with vertex $\widehat{\boldsymbol{x}}_1$, and ending with the edge  shared with vertex $\widehat{\boldsymbol{x}}_2$.
\item The nodes inside the faces of the tetrahedron are indexed such that the normal vector points outwards, following the right-hand rule.
\item The indexing of the internal nodes follows the convention of lower-order tetrahedral elements. For the finite element space $\mathbb{P}_4$, there is only one node in the interior of $\widehat{K}$ and its index is trivial. For $m\geq 5$, the internal nodes are distributed and indexed following the idea of the $\mathbb{P}_{m-4}$ reference element.
\end{itemize}

To show the well-posedness of the interpolation problem, that is, that a polynomial of degree $m$ is uniquely determined by its value at this set of points, it suffices to construct the associated Lagrange basis $\{\widehat{N}_r\}_{r=0}^{\tt dofK-1}\subset\mathbb{P}_m$. That is,  $m$th degree polynomials that vanish at all nodes except for one.  For $\mathbb{P}_1$ it is almost trivial: the associated Lagrange basis is\index{$\widehat{N}_\ell$}
\[
\widehat{N}_{r}(\widehat{\boldsymbol{x}}) = \widehat{\lambda}_{r}(\boldsymbol{x}),\quad r = 0,1,2,3,
\]
i.e., the barycentric coordinates themselves, since in this case, the vertices $\{\widehat{\boldsymbol{x}}_r\}_{r=0}^3$  of $\widehat{K}$  are just the interpolating nodes   {(see   \eqref{eq:lambdahat})} and therefore
\begin{equation*}
    \widehat{N}_r(\widehat{\boldsymbol{x}}_s)=\delta_{rs}, \qquad  r,s \in \left\{0,1, 2,3 \right\}.
\end{equation*}

Consider now the case of higher degree polynomials, i.e. $m\ge 2$.  Let
$\widehat{\boldsymbol{x}}_r$ be the $r$th node with barycentric coordinates given by
\[
\widehat{\boldsymbol{\lambda}}(\widehat{\boldsymbol{x}}_r)= \left(\widehat{\lambda}_0(\widehat{\boldsymbol{x}}_r),\widehat{\lambda}_1(\widehat{\boldsymbol{x}}_r),\widehat{\lambda}_2(\widehat{\boldsymbol{x}}_r),\widehat{\lambda}_3(\widehat{\boldsymbol{x}}_r)\right)=\left(\frac{i_0}m,\frac{i_1}m,\frac{i_2}m,\frac{i_3}m\right).
\]
From   \eqref{eq:barycentricNodes}, the barycentric coordinates satisfy $i_0+i_1+i_2+i_3=m$. The Lagrange function associated to this node is given by
\begin{equation}\label{eq:Nj3D}
\widehat{N}_r(\widehat{\boldsymbol{x}}) = \prod_{n=0}^3\,\prod_{\ell=0}^{i_n-1}\, \frac{\widehat{\lambda}_n(\widehat{\boldsymbol{x}})-\frac{\ell}{m}}{\widehat{\lambda}_n(\widehat{\boldsymbol{x}}_r)-\tfrac{\ell}{m}}.
 \end{equation}
These functions are commonly known as the shape functions associated to the set of nodes ${\boldsymbol x}_r$ (cf. \cite{zienkiewicz_fem}).  It is a simple exercise to check that $\widehat{N}_{r}$ is a polynomial of degree $m$ which vanishes at each node $\widehat{\boldsymbol{x}}_s$ except at $\widehat{\boldsymbol{x}}_r$ where it takes the value 1. In other words, $\{\widehat{N}_{r}\}_{r=0}^{{\tt dofK}-1}$ is the Lagrange basis of the interpolation problem at this set of nodes. Therefore
\begin{equation}\label{eq:polRefElem}
    p_m^{\widehat{K}}(\widehat{\boldsymbol{x}})=\sum_{r=0}^{{\tt dofK}-1} p_m^{\widehat{K}}(\widehat{\boldsymbol{x}}_r)\widehat{N}_r(\widehat{\boldsymbol{x}} ), \quad 
    \forall p_m^{\widehat{K}}(\widehat{\boldsymbol{x}})\in\mathbb{P}_m.
 \end{equation}
Furthermore, for $\widehat{\pi}_s = \{{\boldsymbol x}\ : \ \widehat{\lambda}_s(\boldsymbol{x})=0\}$,  $s=0,1,2,3$, one of the four planes containing the faces of $\widehat{K}$,

\begin{equation}\label{eq:polRefElem:02}
   p_m^{\widehat{K}}\big|_{\widehat{\pi}_s}=\sum_ { \boldsymbol{x}_r\in \widehat{\pi}_s}  p_m^{\widehat{K}}(\widehat{\boldsymbol{x}}_r)\widehat{N}_r|_{\widehat{\pi}_s}.
\end{equation}
In other words,  $p_m^{\widehat{K}}\big|_{\widehat{\pi}_s}$ is fully determined by the values the polynomial 
takes at the subset of nodes that lie in $\pi_s$. A similar condition applies to the straight lines that contain the edges of the tetrahedron.

\subsection{The Lagrange $\mathbb{P}_m-$finite element space}\label{sec: Affine Mapping}

\begin{figure}[!htb]
    \centering
\begin{tikzpicture}[scale  =1]
  \draw[-latex,blue] (0,0,0) -- (2,0,0) node[above] {$\widehat{y}$};
  \draw[-latex, above,blue] (0,0,0) -- (0,2,0) node[below left] {$\widehat{z}$};
  \draw[-latex, blue] (0,0,0) -- (0,0,2) node[left] {$\widehat{x}$};
  \coordinate (A) at (0,0,0);
  \coordinate (B) at (1.5,0,0);
  \coordinate (C) at (0,1.5,0);
  \coordinate (D) at (0,0,1.5);

  \draw[thick] (A) -- (B) -- (C) -- cycle;
  \draw[thick] (A) -- (D) -- (B);
  \draw[thick] (D) -- (C);

  \node[left] at (A) {\small$\widehat{\boldsymbol{x}}_0$};
  \node[above] at (B) {\small$\widehat{\boldsymbol{x}}_2$};
  \node[left] at (C) {\small$\widehat{\boldsymbol{x}}_2$};
  \node[above left] at (D) {\small$\widehat{\boldsymbol{x}}_1$};
  \node[circle,fill=blue, scale=0.3] at (A) {};
  \node[circle,fill=blue, scale=0.3] at (B) {};
  \node[circle,fill=blue, scale=0.3] at (C) {};
  \node[circle,fill=blue, scale=0.3] at (D) {};

  \coordinate (E) at (6.5,3,3);
  \coordinate (F) at (8,2,4);
  \coordinate (G) at (7,5,4);
  \coordinate (H) at (6,4,6);

  \draw[thick, dashed,dash pattern=on 1pt off 2pt] (E) -- (F) -- (G) -- cycle;
  \draw[thick,dashed,dash pattern=on 1pt off 2pt] (E) -- (H) -- (F);
  \draw[thick] (H) -- (G)--(F)--cycle;

  \node[ right] at (E) {\small$\boldsymbol{x}_0$};
  \node[below] at (F) {\small$\boldsymbol{x}_2$};
  \node[right] at (G) {\small$\boldsymbol{x}_3$};
  \node[left] at (H) {\small$\boldsymbol{x}_1$};
  \node[circle,fill=red, scale=0.3] at (E) {};
  \node[circle,fill=red, scale=0.3] at (F) {};
  \node[circle,fill=red, scale=0.3] at (G) {};
  \node[circle,fill=red, scale=0.3] at (H) {};

  \draw[-latex, bend right=70,dashed] (A) to (E);
  \draw[-latex, bend right=40,dashed] (B) to (F);
  \draw[-latex, bend left=30,dashed] (C) to (G);
  \draw[-latex, bend right=50,dashed] (D) to (H);

\coordinate (A1) at ($ (A)!1/4!(B)!1/4!(C)!1/4!(D) $);
\node[circle,fill=purple, scale=0.3] at (A1) {};

\coordinate (A2) at ($ (E)!1/4!(F)!1/4!(G)!1/4!(H) $);
\node[circle,fill=purple, scale=0.3] at (A2) {};

\coordinate (Barycenter) at ($ (A1)!1/2!(A2) $);

\draw[-latex, purple, bend left=65, thick] (A1) node [right] {\small$\widehat{\boldsymbol{x}}$} to[bend left=45] node[midway, above, sloped] {$F_K$} (A2) node [above] {\small$\boldsymbol{x}$};

\end{tikzpicture}

    \caption{Mapping from $\widehat{K}$ to $K$.}
    \label{fig:FK}
\end{figure}
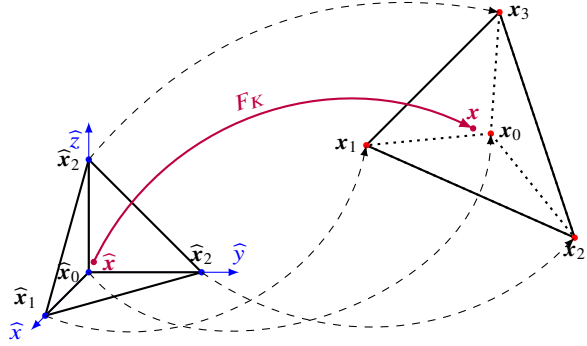

Let $K$ be an arbitrary tetrahedron of vertices $\{\boldsymbol{x}_0, \boldsymbol{x}_1, \boldsymbol{x}_2, \boldsymbol{x}_3\}$.
The affine transformation \index{$F_K$}$F_K : \widehat{K} \longrightarrow K$    defined as
\begin{align}\label{eq: FT}
F_K(\widehat{\boldsymbol{x}})
&=\begin{bmatrix}
    x_1-x_0& x_2-x_0& x_3-x_0\\
    y_1-y_0 &y_2-y_0 &y_3-y_0\\
    z_1-z_0 &z_2-z_0 &z_3-z_0
\end{bmatrix}\begin{bmatrix}
    \widehat{x}\\\widehat{y}\\\widehat{z}
\end{bmatrix}+ \begin{bmatrix}
    x_0\\y_0\\z_0
\end{bmatrix} = \begin{bmatrix}
    \boldsymbol{x}_1-\boldsymbol{x}_0\ & \boldsymbol{x}_2-\boldsymbol{x}_0\ & \boldsymbol{x}_3-\boldsymbol{x}_0
\end{bmatrix}\widehat{\boldsymbol{x}}+\boldsymbol{x}_0\nonumber\\
   &= B_K \widehat{\boldsymbol{x}} + \boldsymbol{x}_0,
\end{align}
is a bijection (see  Figure \ref{fig:FK}) that maps a point in the reference element to a point $\boldsymbol{x}$ inside $K$, i.e., $\boldsymbol{x}= F_K(\widehat{\boldsymbol{x}})$. Alternatively, we can write  \eqref{eq: FT} as 
\begin{equation}\label{eq: baryc_mapping}
    \boldsymbol{x}=F_K(\widehat{\boldsymbol{x}}) = \sum_{n=0}^3 \widehat{\lambda}_n(\widehat{\boldsymbol{x}}) \boldsymbol{x}_n
\end{equation}
which, in turn, shows that $F_K$ preserves the barycentric coordinates.

For any $K\in {\cal T}_h$ we construct $F_K$ as in   \eqref{eq: FT}. We then define
\[
\boldsymbol{x}_r^{K} = F_K(\widehat{\boldsymbol{x}}_r), \quad  r = 0,\ldots, {\tt dofK}-1,  \text{with } {\tt dofK} = \binom{m+3}{3},
\]
the (local) set of nodes in $K$  and, consequently, the nodes in the mesh. 

 Any  polynomial $p_m^K(\boldsymbol{x})\in\mathbb{P}_m$ with support in $K$ satisfies  
\[
  p_m^{\widehat{K}}=p_m^K  \circ F_K  \in \mathbb{P}_m.
\]
Therefore, from   \eqref{eq:polRefElem}:
\[
     p^K_m (\boldsymbol{x}) = \sum_{r=0}^{{\tt dofK}-1} p_m^K(F_K(\widehat{\boldsymbol{x}}_r))  (\widehat{N}_r \circ F_K^{-1}(\boldsymbol{x})).
\]
That is, Condition  \ref{CA} is  automatically satisfied.

Besides, if we let  ${\boldsymbol x}_j \in e= \overline{K}_1\cap \overline{K}_2$ be a node in the mesh, there exists a pair of local indices $r_1,\ r_2$ such that
\[
\boldsymbol{x}_j= \boldsymbol{x}_{r_1}^{K_1}= \boldsymbol{x}_{r_2}^{K_2}.
\]
That is, the set of local nodes on each element matches common faces and edges.  Therefore, we have that for any pair of polynomials $p_m^{K_1}, p_m^{K_2}\in \mathbb{P}_m$:
\begin{equation}\label{eq: polLocalBasis3}
p_m^{K_1}\big|_{e} = p_m^{K_2}\big|_{e}\quad \text{if and only if}\quad  p_m^{K_1}  (\boldsymbol{x}_{r_1}^{K_1})=p_m^{K_2} (\boldsymbol{x}_{r_2}^{K_2}), \quad \text{for}\quad
\boldsymbol{x}_{r_1}^{K_1} = \boldsymbol{x}_{r_2}^{K_2}\in e, 
\end{equation}
that is, the polynomials  match at their interface $e$ if and only if they match at the nodes in the mesh in $e$. 

Now, let  us define a piecewise polynomial function $\psi_j$ satisfying:
\[
\psi_j|_{K} \in\mathbb{P}_m, \quad \psi_j(\boldsymbol{x}_k^K) = \begin{cases}
1, & \text{if }\boldsymbol{x}_{j}^K = \boldsymbol{x}_k^K\\
0, & \text{otherwise}
\end{cases}
\]
This function is uniquely  defined by Condition \ref{CA}, it is indeed continuous by  \eqref{eq: polLocalBasis3} and  therefore it belongs to $\mathcal{P}_h^m$, and so it is the basis function $\varphi_j$. In other words,  Condition \ref{CB} is also satisfied.

\subsection{Lagrange $\mathbb{P}_m$-elements at the boundary}\label{subsec: ElmBoundaries}

The treatment of the boundary conditions in the FEM, particularly Dirichlet but much more  Robin conditions, will require a suitable description of the inherited partition of the boundary into triangles, which are the faces of the elements that lie on the boundary. This corresponds to a continuous 2D $\mathbb{P}_m$-finite element space on a triangular mesh, which is entirely analogous to the 3D case discussed in previous subsections. We will provide a quick and simple description of this space, including the local set of nodes and their relationship with the nodes in one of the faces in $K$.

Let the reference triangle be defined as
\[
\widehat{A} = \{ \widehat{\boldsymbol{x}} =(\widehat{x},\widehat{y})\ \lvert\ \, 0< 1-\widehat{x}-\widehat{y},\ 0 <\widehat{x},\ 0 <\widehat{y}\}, 
\]
with vertices 
\[
\begin{aligned}
\widehat{\boldsymbol{x}}_0   = (0,0),\quad
\widehat{\boldsymbol{x}}_1   = (1,0),\quad 
\widehat{\boldsymbol{x}}_2  = (0,1).
\end{aligned}
\]
For any triangle $A\subset\mathbb{R}^3$ with vertices $\{\boldsymbol{x}_0,\boldsymbol{x}_1,\boldsymbol{x}_2\}$, we define the affine mapping $F_A:\widehat{A}\longrightarrow A$ as
\begin{align}\label{eq: FA}
F_A(\widehat{\boldsymbol{x}})
&=\begin{bmatrix}
    \boldsymbol{x}_1-\boldsymbol{x}_0& \boldsymbol{x}_2-\boldsymbol{x}_0
    \end{bmatrix} \begin{bmatrix}
    \widehat{x}\\ \widehat{y}\\ 
\end{bmatrix}+ \boldsymbol{x}_0 =  \begin{bmatrix}
    \boldsymbol{\ell}_1\ &\boldsymbol{\ell}_2 
    \end{bmatrix} \begin{bmatrix}
    \widehat{x}\\ \widehat{y}\\ 
\end{bmatrix}+ \boldsymbol{x}_0,
\end{align}
which can be recast as
\[
F_A(\widehat{\boldsymbol{x}}) = \sum_{j=0}^2 \widehat{\lambda}_j(\widehat{\boldsymbol{x}}) \boldsymbol{x}_j, \quad
\widehat{\lambda}_j(\widehat{\boldsymbol{x}}) = 
\begin{cases} 
1 - \widehat{x} - \widehat{y}, & j=0, \\ 
\widehat{x}, & j=1, \\ 
\widehat{y}, & j=2.
\end{cases}
\]
Then 
\[
\widehat{\boldsymbol{\lambda}}(\widehat{\boldsymbol{x}})= (\widehat{\lambda}_0(\widehat{\boldsymbol{x}}),\widehat{\lambda}_1(\widehat{\boldsymbol{x}}),\widehat{\lambda}_2(\widehat{\boldsymbol{x}}))
\]
are the (2D) barycentric coordinates of both  $\widehat{\boldsymbol{x}}$ with respect to $\{\widehat{\boldsymbol{x}}_0,\widehat{\boldsymbol{x}}_1,\widehat{\boldsymbol{x}}_2\}$ and $\boldsymbol{x}=F_A(\widehat{\boldsymbol{x}}) $ with respect to $\{\boldsymbol{x}_0, \boldsymbol{x}_1, \boldsymbol{x}_2\}$.

\begin{figure}[!htb]
\begin{subfigure}{0.5\textwidth}
\[
\begin{tikzpicture}[scale=3.25]

\coordinate (x0) at (0,-0.1,-.5);
\coordinate (x2) at (1,0.2,-0.3);
\coordinate (x3) at (.2,1,0.2);
\coordinate (x1) at (0.1,0.1,1);
\node[left=5pt] at (x0) {$x_0$};
\node[left=5pt] at (x1) {$x_1$};
\node[right=5pt] at (x2) {$x_2$};
\node[above=5pt] at (x3) {$x_3$};
\draw[opacity=0.3] (x0) -- (x1) -- (x2) -- cycle;
\draw[opacity=0.3] (x0) -- (x3) -- (x1);
\draw[opacity=0.3] (x3) -- (x2);

\foreach \i/\a/\b/\c/\d in {
    0/1/0/0/0,
    1/0/1/0/0,
    2/0/0/1/0,
    3/0/0/0/1,
    4/0.66667/0.33333/0/0,
    5/0.33333/0.66667/0/0,
    6/0/0.66667/0.33333/0,
    7/0/0.33333/0.66667/0,
    8/0.33333/0/0.66667/0,
    9/0.66667/0/0.33333/0,
    10/0.33333/0/0/0.66667,
    11/0.66667/0/0/0.33333,
    12/0/0/0.33333/0.66667,
    13/0/0/0.66667/0.33333,
    14/0/0.33333/0/0.66667,
    15/0/0.66667/0/0.33333,
    16/0.33333/0.33333/0.33333/0,
    17/0.33333/0.33333/0/0.33333,
    18/0.33333/0/0.33333/0.33333,
    19/0/0.33333/0.33333/0.33333}{
    \node[fill=none,opacity = 0.4] (node\i) at
     ($\a*(x0) + \b*(x1) + \c*(x2) + \d*(x3)$) {}; 
    \shade[ball color=blue, opacity=0.1]
        ($\a*(x0) + \b*(x1) + \c*(x2) + \d*(x3)$)
        circle[radius=0.06];
}

\foreach \i/\a/\b/\c/\d in {
    1/0/1/0/0,
    2/0/0/1/0,
    3/0/0/0/1,
    6/0/0.66667/0.33333/0,
    7/0/0.33333/0.66667/0,
    12/0/0/0.33333/0.66667,
    13/0/0/0.66667/0.33333,
    14/0/0.33333/0/0.66667,
    15/0/0.66667/0/0.33333,
    19/0/0.33333/0.33333/0.33333}{
    \node[fill=none] (node\i) at
     ($\a*(x0) + \b*(x1) + \c*(x2) + \d*(x3)$) {};
    \shade[ball color=blue, opacity=0.5]
        ($\a*(x0) + \b*(x1) + \c*(x2) + \d*(x3)$)
        circle[radius=0.06];
}
\end{tikzpicture}
\]
\end{subfigure}
\begin{subfigure}{0.45\textwidth}
\[
 \begin{tikzpicture}[scale=3.15]
 
  \node[circle,fill=black,scale=0.7] (A) at (0,0) {} ;
  \node[circle,fill=black,scale=0.7] (B) at (0,1) {};
  \node[circle,fill=black,scale=0.7] (C) at (1,0){};
  
    \node[below left] at (A) {$(0,0)$};
  \node[below] at (C) {$(1,0)$};
  \node[left] at (B) {$(0,1)$};

 \draw[-latex,blue] (-0.2,0) -- (1.3,0) node[above] {$\widehat{x}$};
  \draw[-latex, above,blue] (0,-0.2) -- (0,1.3) node[below left] {$\widehat{y}$};

\coordinate (x0) at (0,0,0);
\coordinate (x1) at (1,0,0);
\coordinate (x2) at (0,1,0);

\draw[opacity = 1, black] (x0) -- (x1) -- (x2) -- cycle;

\foreach \i/\a/\b/\c in {
    0/1/0/0,
    1/0/1/0,
    2/0/0/1,
    3/0.6666667      /0.3333333            /0         ,
    4/0.3333333      /0.6666667           /0         ,
    5/0      /0.6666667      /0.3333333           , 
    6/0      /0.3333333      /0.6666667           ,
    7/0.75            /0      /0.6666667     , 
    8/0.6666667            /0      /0.33333333333 ,
    9/0.33333333333      /0.33333333333   /0.33333333333    
    }{
 \draw[color=black, fill=blue!40, opacity=1] ($\a*(x0) + \b*(x1) + \c*(x2) $) circle[radius=0.04, color= blue];
        \node[fill=none]  (node\i) at
     ($\a*(x0) + \b*(x1) + \c*(x2) $) {};
}
\end{tikzpicture}
\]
\end{subfigure}
\caption{Examples of node heritage for $\mathbb{P}_2$ elements. \label{fig:A:K}Left: Nodes in a tetrahedron $K$ and the inherited nodes in a triangular face $A$. Right: Nodes in the reference triangle $\widehat{A}$}
\end{figure}
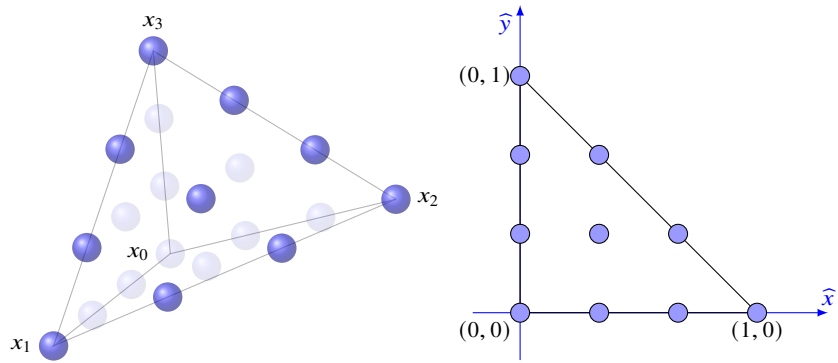

If $A$ is a face (triangle) of an element (tetrahedra) $K$ of the mesh, we have the set of local nodes $\{\widehat{\boldsymbol{x}}_j^A\}_{j=0}^{{\tt dofA}-1}$ inherited from that of $K$. Via $F_A$ it moves to the reference triangle $\widehat{A}$ and are given, similarly to   \eqref{eq:polRefElem:02}, by    (see Figure \ref{fig:A:K})
\[
\widehat{\boldsymbol{\lambda}}(\widehat{\boldsymbol{x}}_r)=\left(\frac{i_0}m,\frac{i_1}m,\frac{i_2}m\right),\quad i_1+i_2+i_3=m,\quad i_0,\ i_1,\ i_2\ge 0.
\] 

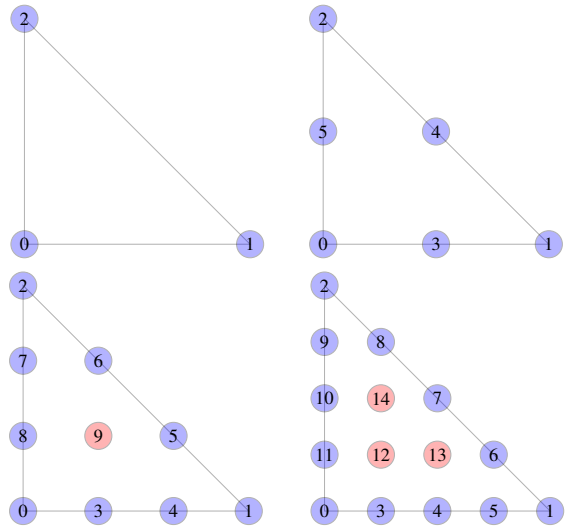
\begin{figure}
\[
    \begin{tikzpicture}[scale=3]
\coordinate (x0) at (0,0,0);
\coordinate (x1) at (1,0,0);
\coordinate (x2) at (0,1,0);

\draw[opacity = 0.3] (x0) -- (x1) -- (x2) -- cycle;

\foreach \i/\a/\b/\c in {
    0/1/0/0,
    1/0/1/0,
    2/0/0/1
    }{
 \draw[color=black, fill=blue, opacity=0.3] ($\a*(x0) + \b*(x1) + \c*(x2) $) circle[radius=0.06, color= blue];
        \node[fill=none]  (node\i) at
     ($\a*(x0) + \b*(x1) + \c*(x2) $) {\scriptsize \i};
}
\end{tikzpicture}\qquad
\begin{tikzpicture}[scale=3]

\coordinate (x0) at (0,0,0);
\coordinate (x1) at (1,0,0);
\coordinate (x2) at (0,1,0);

\draw[opacity = 0.3] (x0) -- (x1) -- (x2) -- cycle;

\foreach \i/\a/\b/\c in {
    0/1/0/0,
    1/0/1/0,
    2/0/0/1,
    3/0.5      /0.5            /0         , 
    4/0      /0.5      /0.5           , 
    5/0.5            /0      /0.5          
    }{
 \draw[color=black, fill=blue, opacity=0.3] ($\a*(x0) + \b*(x1) + \c*(x2) $) circle[radius=0.06, color= blue];
        \node[fill=none]  (node\i) at
     ($\a*(x0) + \b*(x1) + \c*(x2) $) {\scriptsize \i};
}
\end{tikzpicture}
\]
\[
 \begin{tikzpicture}[scale=3]

\coordinate (x0) at (0,0,0);
\coordinate (x1) at (1,0,0);
\coordinate (x2) at (0,1,0);

\draw[opacity = 0.3] (x0) -- (x1) -- (x2) -- cycle;

\foreach \i/\a/\b/\c in {
    0/1/0/0,
    1/0/1/0,
    2/0/0/1,
    3/0.6666667      /0.3333333            /0         ,
    4/0.3333333      /0.6666667           /0         ,
    5/0      /0.6666667      /0.3333333           , 
    6/0      /0.3333333      /0.6666667           ,
    7/0.75            /0      /0.6666667     , 
    8/0.6666667            /0      /0.33333333333    
    }{
 \draw[color=black, fill=blue, opacity=0.3] ($\a*(x0) + \b*(x1) + \c*(x2) $) circle[radius=0.06, color= blue];
        \node[fill=none]  (node\i) at
     ($\a*(x0) + \b*(x1) + \c*(x2) $) {\scriptsize \i};}

     \foreach \i/\a/\b/\c in {
    9/0.33333333333      /0.33333333333   /0.33333333333    
    }{
 \draw[color=black, fill=red, opacity=0.3] ($\a*(x0) + \b*(x1) + \c*(x2) $) circle[radius=0.06, color= blue];
        \node[fill=none]  (node\i) at
     ($\a*(x0) + \b*(x1) + \c*(x2) $) {\scriptsize \i};
}
\end{tikzpicture}\qquad
\begin{tikzpicture}[scale=3]

\coordinate (x0) at (0,0,0);
\coordinate (x1) at (1,0,0);
\coordinate (x2) at (0,1,0);

\draw[opacity = 0.3] (x0) -- (x1) -- (x2) -- cycle;

\foreach \i/\a/\b/\c in {
    0/1/0/0,
    1/0/1/0,
    2/0/0/1,
    3/0.75      /0.25            /0         ,
    4/0.5      /0.5            /0           ,
    5/0.25      /0.75            /0         ,
    6/0      /0.75      /0.25            ,
    7/0      /0.5      /0.5            ,
    8/0      /0.25      /0.75           ,
    9/0.75            /0      /0.75     ,
    10/0.5            /0      /0.5       ,
    11/0.25            /0      /0.25        
    }{
 \draw[color=black, fill=blue, opacity=0.3] ($\a*(x0) + \b*(x1) + \c*(x2) $) circle[radius=0.06, color= blue];
        \node[fill=none]  (node\i) at
     ($\a*(x0) + \b*(x1) + \c*(x2) $) {\scriptsize \i};}
     
     \foreach \i/\a/\b/\c in {
    12/0.5      /0.25   /0.25    ,
    13/0.25      /0.5   /0.25    ,
    14/0.25      /0.25   /0.5    
    }{
 \draw[color=black, fill=red, opacity=0.3] ($\a*(x0) + \b*(x1) + \c*(x2) $) circle[radius=0.06, color= blue];
        \node[fill=none]  (node\i) at
     ($\a*(x0) + \b*(x1) + \c*(x2) $) {\scriptsize \i};}
\end{tikzpicture} 
    \]
     \caption{\label{fig:LagrangeNodesBound} Nodes numbering on the $\mathbb{P}_1$ (top-left), $\mathbb{P}_2$ (top-right), $\mathbb{P}_3$ (bottom-left) and $\mathbb{P}_4$ (bottom-right) reference boundary (triangular) elements.}
\end{figure}

Once more, we must handle the ordering of the ${\tt dofA} = \binom{m+2}{2}$ nodes, which can be indexed in different ways. We follow the GMSH convention as illustrated in Figure \ref{fig:LagrangeNodesBound} for elements of degrees 1 through 4:
\begin{itemize}
\item Firstly, the vertices are indexed using the same order as above.
\item Next, the nodes inside each edge of the triangle are indexed, starting from the edge lying on the $\widehat{x}$ axis and following a counterclockwise indexing.
\item Finally, the interior nodes are defined recursively, taking them from the definition of lower-order elements. For $\mathbb{P}_m$ elements, the interior nodes are indexed following the rules for $\mathbb{P}_{m-3}$ elements.
\end{itemize}

Finally we have the Lagrange basis for the interpolation problem $\widehat{N}_r$, the  (2D) counterpart of the functions in the tetrahedra (see   \eqref{eq:Nj3D}): 
\begin{equation*}
\widehat{N}_r(\widehat{\boldsymbol{x}}) = \prod_{n=0}^2\,\prod_{\ell=0}^{i_r-1}\, \frac{\widehat{\lambda}_n(\widehat{\boldsymbol{x}})-\frac{\ell}{m}}{\widehat{\lambda}_n(\widehat{\boldsymbol{x}}_r)-\tfrac{\ell}{m}}.
 \end{equation*}

\section{Node and Element Data Extraction from GMSH mesh files}\label{sec: GMSH_msh}

We previously discussed the relationship between nodes in the reference element and their mapping to arbitrary elements. Note that the indexing is not unique, as the choice of the vertex $\boldsymbol{x}_0^K$ is arbitrary even after imposing a positive determinant in $B_K$. GMSH manages the meshing process and supplies the necessary data to define geometric and mesh properties. In version {\tt 4.1} (the latest version available at the time of this publication), the mesh file with extension {\tt *.msh} contains relevant geometric and mesh information.

\begin{figure}[!htb]
    \centering
    \includegraphics[scale=0.3]{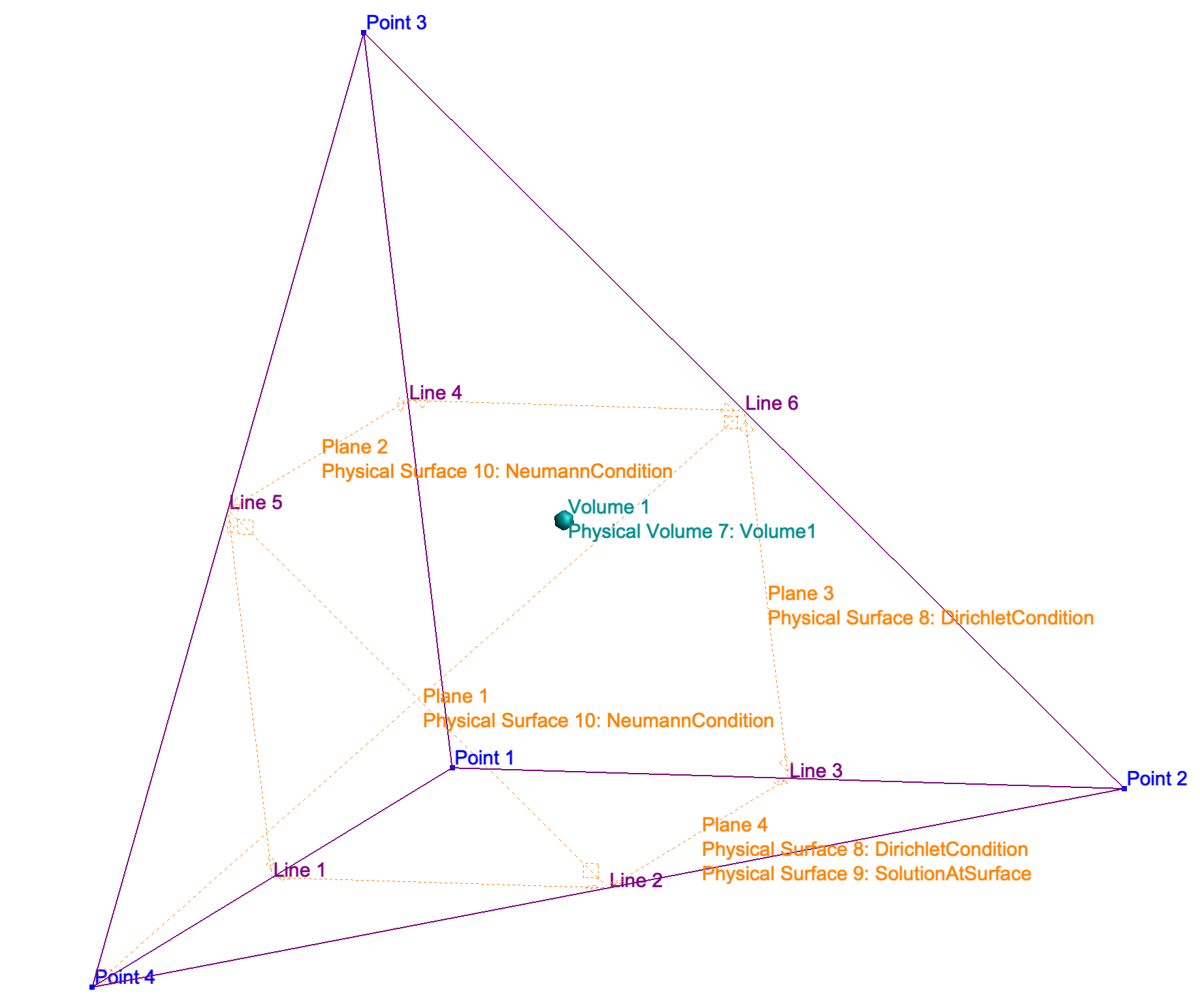}
    \caption{Sample geometry case in GMSH.}
    \label{fig:testGMSH}
\end{figure}

A general {\tt *.msh} file contains much more information than we will discuss here. For more details, see the GMSH documentation. Here, we focus on the standard file obtained after generating (not so) {\it simple} geometries and meshes. We use the geometry shown in Figure \ref{fig:testGMSH} to illustrate the data arrangement in the mesh file and explain how to store the mesh information. The geometry is a tetrahedral domain with $\mathbb{P}_1$ elements. During the {\tt *.msh} storage, we chose the {\tt Save All} options to store all the mesh information. For more detailed information, we refer the reader to \cite{duque2023integration} and the official GMSH tutorial.

\subsection{Physical Groups}\label{sec: PhysicalGroups}

During the geometry generation, it is possible to categorize the geometric entities used to construct $\Omega$ by organizing them into physical groups. These groups allow us to associate the generated entities with specific material properties, boundary conditions, or other physical attributes important for simulation and post-processing. Figure \ref{fig: PhysicalGroups} illustrates a typical example of this classification for the sample case. In the example considered, we observe that there are three 2D physical groups and only one 3D physical group.

\begin{figure}[!htb]
\begin{alltt}
$PhysicalNames
\textcolor{red}{4}
\textcolor{blue}{2} \textcolor{purple}{8} \textcolor{orange}{"DirichletCondition"}
\textcolor{blue}{2} \textcolor{purple}{9} \textcolor{orange}{"SolutionAtSurface"}
\textcolor{blue}{2} \textcolor{purple}{10} \textcolor{orange}{"NeumannCondition"}
\textcolor{blue}{3} \textcolor{purple}{7} \textcolor{orange}{"Volume1"}
$EndPhysicalNames
\end{alltt}
\caption{Physical Names section. The dimension of each physical group is indicated in blue. The number of Physical groups is shown in red. The numerical labels are shown in purple and the Physical Names are shown in orange.\label{fig: PhysicalGroups}  }
\end{figure}
First, GMSH provides the number of Physical groups defined for a given geometry. Then it shows a list of the Physical groups, ordered by the dimension of the entities associated to the group and followed by both the numerical and string label. Naturally, it is also possible to set Physical group to zero- and one-dimensional entities, but this is not considered in our implementation.

\subsection{Entities}\label{sec: Entities}
This section list all geometric entities used to define the domain $\Omega$ by using a boundary representation. Figure \ref{fig: entities} illustrates the structure of this section for the sample case.

\begin{figure}[!htb]
\begin{multicols}{2}
\begin{alltt}
\$Entities
\textcolor{blue}{4} \textcolor{purple}{6} \textcolor{orange}{4} \textcolor{darkcyan}{1}
\textcolor{blue}{1} 0 0 0 \textcolor{blue}{0} 
\textcolor{blue}{2} 1 0 0 \textcolor{blue}{0} 
\textcolor{blue}{3} 0 1 0 \textcolor{blue}{0}  
\textcolor{blue}{4} 0 0 1 \textcolor{blue}{0} 
\textcolor{purple}{1} 0 0 0 0 0 1 \textcolor{purple}{0} 2 1 -4 
\textcolor{purple}{2} 0 0 0 1 0 1 \textcolor{purple}{0} 2 4 -2 
\textcolor{purple}{3} 0 0 0 1 0 0 \textcolor{purple}{0} 2 2 -1 
\textcolor{purple}{4} 0 0 0 0 1 0 \textcolor{purple}{0} 2 1 -3 
\textcolor{purple}{5} 0 0 0 0 1 1 \textcolor{purple}{0} 2 3 -4 
\textcolor{purple}{6} 0 0 0 1 1 0 \textcolor{purple}{0} 2 3 -2 
\textcolor{orange}{1} 0 0 0 1 1 1 \textcolor{orange}{1 10} 3 5 2 -6 
\textcolor{orange}{2} 0 0 0 0 1 1 \textcolor{orange}{1 10} 3 1 -5 -4 
\textcolor{orange}{3} 0 0 0 1 1 0 \textcolor{orange}{1} \textcolor{orange}{8} 3 3 4 6 
\textcolor{orange}{4} 0 0 0 1 0 1 \textcolor{orange}{2} \textcolor{orange}{8} \textcolor{orange}{9} 3 1 2 3 
\textcolor{darkcyan}{1} 0 0 0 1 1 1 \textcolor{darkcyan}{1} \textcolor{darkcyan}{7} 4 1 2 4 3 
\$EndEntities
\end{alltt}
\end{multicols}
\caption{Entities section. First line contains the number of entities of every dimension. Each other row contains information of an entity used to define the geometry. Each color is used to represent entities with the same dimension, highlighting the quantities of interest. \label{fig: entities}}
\end{figure}

The first row shows the number of entities for each dimension, starting from zero-dimensional to three-dimensional entities. The geometry consists of four points, six lines, four surfaces and one volume.

The entities are listed by dimension and their numerical label (first highlighted number on each row) as assigned during the geometry generation. In this example, the numerical indexing of the entities starts at one and is continuous, but this is not always the case. After the numerical label, it shows the minimum and maximum bounds for each coordinate.

The following highlighted numbers indicate the number of Physical Groups associated with an entity, followed by the numerical labels of the associated Physical Groups.  In this example, all surfaces in the boundary and the volume have an associated Physical Group. Surfaces {\tt 1} and {\tt 2} are associated to \texttt{NeumannCondition}. Surface {\tt 3} and {\tt 4} are associated to {\tt DirichletCondition}. Surface{\tt 4} is also associated to {\tt SolutionAtSurface}. Volume {\tt 1} is associated to {\tt Volume1}.

The rest of the information in each row lists the entities of one dimension lower used to define the entity. The sign indicates whether the orientation of the embedded entity has been adjusted to match the orientation of the entity it is defining. If we refer back to Figure \ref{fig:testGMSH} and focus on {\tt Plane 2},  we observe that this surface is constructed from three curves: {\tt Line 1}, {\tt Line 4} and {\tt Line 5}  (as can be seen in {\tt 3 1 -5 -4}).  {\tt Line 1} is an oriented segment from {\tt Point 1} to {\tt Point 4}, {\tt Line 4} is oriented from {\tt Point 1} to {\tt Point 3} and {\tt Line 5} from {\tt Point 3} to {\tt Point 4} (see the respective rows). The orientation of the surface is defined following the direction of the first segment (in this case {\tt 1}). The negative signs in the other lines are to match the orientation of {\tt Line 1} and achieve a coherent orientation in Surface {\tt 2}. This definition for the orientation also applies for volumes using surfaces rather than curves.

 \subsection{Mesh}
The mesh consists of nodes organized in elements (points, faces, edges and tetrahedra). 
\subsubsection{Nodes}\label{NODES}
 The nodes within the mesh are given in Section {\tt Nodes}. Figure \ref{fig:nodes} illustrates the structure of this section for the sample case. 

\begin{figure}
    \centering
    \begin{multicols}{4}
\begin{alltt}
\$Nodes
15 7 1 7
\textcolor{blue}{0 1 0 1}
1
0 0 0
\textcolor{blue}{0 2 0 1}
2
1 0 0
\textcolor{blue}{0 3 0 1}
3
0 1 0
\textcolor{blue}{0 4 0 1}
4
0 0 1
\textcolor{purple}{1 1 0 0}
\textcolor{purple}{1 2 0 1}
5
0.5 0 0.5
\textcolor{purple}{1 3 0 0}
\textcolor{purple}{1 4 0 0}
\textcolor{purple}{1 5 0 1}
6
0 0.5 0.5
\textcolor{purple}{1 6 0 1}
7
0.5 0.5 0
\textcolor{orange}{2 1 0 0}
\textcolor{orange}{2 2 0 0}
\textcolor{orange}{2 3 0 0}
\textcolor{orange}{2 4 0 0}
\textcolor{darkcyan}{3 1 0 0}
\$EndNodes
\end{alltt}    
\end{multicols}
    \caption{Nodes section. Each color is used to represent entities with the same dimension, highlighting the quantities of interest.}
    \label{fig:nodes}
\end{figure}

The first line details the number of entities that have stored data in this section, the total number of nodes in the mesh, and the minimum and maximum indices of those nodes. In this mesh, there are seven nodes, indexed from 1 to 7, organized into fifteen blocks. Each node uniquely appears within the interior of each entity. For each entity, there is a row that includes its dimension, numerical label, a parameter indicating whether the entity is parametric (which is always set to 0 for the cases of interest), and the number of associated nodes. Following this, the indices of all the nodes belonging to the entity are listed, one per row, along with their coordinates. For example, the nodes inside {\tt Line 5} are gathered in the block  {\tt 1 5 0 1}. This representation shows the entity dimensionality, its numerical label, and that it contains one node. The node is labeled with {\tt 6} and with coordinates $(0,0.5,0.5)$.

 In general, the indexing of the nodes is sparse. This can occur, for instance, when {\tt Save All} is disabled, in which case nodes belonging to entities without an associated Physical Group are not stored, leading to missing nodes in the mesh file. To avoid this, we can enable {\tt Save All} or add Physical Groups to all two- and three-dimensional entities.  We redefine the indices so that the order is monotonic and there are no missing indices. The information of the coordinates of all nodes can then be stored in a  $ {\tt nNodes}\times 3$  matrix:
\begin{equation}\label{eq:coord}
    {\tt coord} = \begin{bmatrix}
        \boldsymbol{x}_1^\top\\
        \vdots\\
        \boldsymbol{x}_{\tt nNodes}^\top
    \end{bmatrix}.
\end{equation}

\subsubsection{Elements}\label{sec: GMSH_Elements}
This section contains all the elements defined in the geometry and the nodes belonging to them. Figure \ref{fig :GMSHElement} illustrates the structure of this section for the sample case.

\begin{figure}[!htb]
\begin{multicols}{4}
\begin{alltt}
\$Elements
15 27 1 27
\textcolor{blue}{0 1 15 1}
\textcolor{red}{1} \textcolor{gray}{1}
\textcolor{blue}{0 2 15 1}
\textcolor{red}{2} \textcolor{gray}{2}
\textcolor{blue}{0 3 15 1}
\textcolor{red}{3} \textcolor{gray}{3}
\textcolor{blue}{0 4 15 1}
\textcolor{red}{4} \textcolor{gray}{4}
\textcolor{purple}{1 1 1 1}
\textcolor{red}{5} \textcolor{gray}{1 4}
\textcolor{purple}{1 2 1 2}
\textcolor{red}{6} \textcolor{gray}{4 5}
\textcolor{red}{7} \textcolor{gray}{5 2}
\textcolor{purple}{1 3 1 1}
\textcolor{red}{8} \textcolor{gray}{2 1}
\textcolor{purple}{1 4 1 1}
\textcolor{red}{9} \textcolor{gray}{1 3}
\textcolor{purple}{1 5 1 2}
\textcolor{red}{10} \textcolor{gray}{3 6}
\textcolor{red}{11} \textcolor{gray}{6 4}
\textcolor{purple}{1 6 1 2}
\textcolor{red}{12} \textcolor{gray}{3 7}
\textcolor{red}{13} \textcolor{gray}{7 2}
\textcolor{orange}{2 1 2 4}
\textcolor{red}{14} \textcolor{gray}{7 5 2}
\textcolor{red}{15} \textcolor{gray}{6 7 3}
\textcolor{red}{16} \textcolor{gray}{5 6 4}
\textcolor{red}{17} \textcolor{gray}{5 7 6}
\textcolor{orange}{2 2 2 2}
\textcolor{red}{18} \textcolor{gray}{6 3 1}
\textcolor{red}{19} \textcolor{gray}{4 6 1}
\textcolor{orange}{2 3 2 2}
\textcolor{red}{20} \textcolor{gray}{7 2 1}
\textcolor{red}{21} \textcolor{gray}{3 7 1}
\textcolor{orange}{2 4 2 2}
\textcolor{red}{22} \textcolor{gray}{5 2 1}
\textcolor{red}{23} \textcolor{gray}{4 5 1}
\textcolor{darkcyan}{3 1 4 4}
\textcolor{red}{24} \textcolor{gray}{5 6 1 4}
\textcolor{red}{25} \textcolor{gray}{1 7 3 6}
\textcolor{red}{26} \textcolor{gray}{1 7 5 2}
\textcolor{red}{27} \textcolor{gray}{1 7 6 5}
\$EndElements
\end{alltt}
\end{multicols}
\caption{Element section.  The same color pattern used for nodes is used to represent elements belonging to entities of the same direction. In red, the indexing of the element. In gray, the indexing of the nodes as given in the previous section.\label{fig :GMSHElement}}
\end{figure} 
 
The first line specifies the number of element blocks (15), the total number of elements in the mesh (27), and the minimum (1) and maximum (27) indices of the elements in the mesh. 

Each block represents an entity. It begins with a row that specifies the dimensions of the elements, the entity that contains these elements, the type of element, and the number of elements in the block.

We see four blocks containing elements of dimension 0 (with element type 15, according to GMSH). Each block contains a single element with indices 1 to 4, corresponding to nodes 1, 2, 3, and 4, respectively. 

Following this are six consecutive blocks of dimension 1 (edges). These blocks contain 1, 2, 1, 1, 2, and 2 elements (with Element type 1), numbered from 5 to 13. Each block lists the nodes that define the edges. For example, the edges connect nodes 1 and 4, 4 and 5, 5 and 2, and so on.

Next, four blocks of dimension 2 (2D elements) are associated with the 2D entities with labels from 1 to 4. The first block contains four elements, while the remaining three contain two elements each. The elements are of type 2 ($\mathbb{P}_1-$triangles), numbered from 14 to 23. The nodes defining each element are also provided. For instance, the first triangle is defined by nodes 7, 5, and 2, while the second is defined by nodes 6, 7, and 3.

Finally, a single block of dimension 3 is associated with the only three-dimensional entity, which contains four elements of type 4 ($\mathbb{P}_1-$tetrahedra). The vertices defining these elements are also listed.

Triangular elements are oriented based on the orientation of the entities they belong to. To ensure that the elements are properly defined, the surfaces must be correctly oriented during the geometry generation. Our implementation assumes that this condition is met.

Even though the indexing of the elements can be sparse, we care about their order of appearance. We store the nodes belonging to each tetrahedral element in a  ${\tt nTtrh \times dofK}$ matrix (for the case of tetrahedral elements):
\begin{equation}\label{eq:def:ttrh}
    {\tt ttrh } = \begin{bmatrix}
       i_{1,1} &i_{1,2} &\cdots &i_{1,{\tt dofK}}\\
       i_{2,1} &i_{2,2} &\cdots&i_{2,{\tt dofK}}\\
        \vdots &\vdots &  \ddots & \vdots\\
        i_{{\tt nTtrh},1} &i_{{\tt nTtrh},2} &\cdots&i_{{\tt nTtrh},{\tt dofK}}\\
    \end{bmatrix},
\end{equation}
such that for any given tetrahedron $K_j$, the row ${\tt ttrh}(j,: )$ contains the global indices of the nodes belonigng to the $j-$th tetrahedron.
\[
\begin{bmatrix}
    \boldsymbol{x}_{{\tt ttrh}(j,1)} & \boldsymbol{x}_{{\tt ttrh}(j,2)}& \cdots &\boldsymbol{x}_{{\tt ttrh}(j,{\tt dofK})} \end{bmatrix}=\begin{bmatrix}
\boldsymbol{x}_{0}^{K_j}&\boldsymbol{x}_{1}^{K_j}&\cdots&\boldsymbol{x}_{{\tt dofK}-1}^{K_j}\end{bmatrix}.
\]

Similarly,  We can store the global indices for the triangles in the boundary in a ${\tt ntrB}\times {\tt dofA}$ matrix:
\begin{equation}\label{eq:def:trB}
    {\tt trB} = \begin{bmatrix}
       i_{1,1} &i_{1,2} &\cdots&i_{1,{\tt dofA}}\\
       i_{2,1} &i_{2,2} &\cdots&i_{2,{\tt dofA}}\\
        \vdots &\vdots & \ddots & \vdots\\
        i_{{\tt ntrB},1} &i_{{\tt ntrB},2} &\cdots&i_{{\tt ntrB},{\tt dofA}}\\
    \end{bmatrix}.
\end{equation}

\section{Finite Element Method}
The Finite Element Method can be stated as follows
\begin{equation}\label{eq: weak4}
\left|
\begin{array}{ll}
\text{Find }
u_h\in \mathcal{P}_h^m, \text{ with }u_h\approx u_D  \text{ on $\Gamma_{D}$ such that}\\[2ex]
\begin{aligned}
\underbrace{\iiint_{\Omega} \nabla w_h\cdot (\underline{\boldsymbol{\kappa}}\nabla u_h) +\iint_{\Gamma_{R}} \alpha  w_h u_h +\iiint_{\Omega} w_h(\boldsymbol{\beta}\cdot \nabla u_h)+ \iiint_{\Omega}  cw_h u_h}_{=:a(u_h,w_h)}\\
& \hspace{-3cm}=
\underbrace{\iiint_{\Omega} fw_h+\iint_{\Gamma_{R}} g_Rw_h}_{=:\ell(w_h)}, \\
\multicolumn{2}{r}{\forall w_h\in \{v_h\in \mathcal{P}_h^m\ |\  v_h=0 \, \text{ on } \Gamma_D\}.}
\end{aligned}
\end{array}
\right.
\end{equation}
This {\em discretizes} the problem in the following sense: to solve \eqref{eq: weak4}, i.e., to find $u_h$, we only need to determine the values of $u_h$ at the nodes outside the Dirichlet boundary. Thus, the number of unknowns is finite and is determined by requiring the solution to satisfy a test when multiplied by the elements of (a subspace of) $\mathcal{P}_h^m$, a finite-dimensional space.  Indeed, let us define\index{{\tt iD}} {\tt iD} and {\tt inD}\index{\tt niD} to be the vectors containing the global indices of Dirichlet (where the exact solution is known) and non-Dirichlet (where it is not) nodes respectively. That is, {\tt iD} and {\tt inD} are a partition of the global indices,
\[
{\tt iD} \cup {\tt inD}  = \{1,2,\ldots, {\tt nNodes}\},\quad {\tt iD} \cap {\tt inD} =\emptyset,
\]
and
\[
j\in {\tt iD} \text{ if and only if } {\boldsymbol{x}_j}\in \Gamma_D.
\]
Write $u_h\in \mathcal{P}_h^m$ as 
\begin{equation*}
u\approx u_h = u_{h,D}+u_{h,nD}=\sum_{k\in {\tt iD}} u_D(\boldsymbol{x}_j) \varphi_k + \sum_{j\in {\tt inD}} u_j \varphi_j.
\end{equation*}
The term $u_{h,D}$ is a {\em discrete lifting of $u_D$} since  $u_{h,D}(\boldsymbol{x}_j)=u_D(\boldsymbol{x}_j)$ for $j\in {\tt iD}$. Observe also that  $u_{h,nD}$ vanishes on $\Gamma_D$.  
Therefore,  with
\[
     \mathcal{P}_{h,D}^m:= \{v_h\in \mathcal{P}_h^m\ |\  v_h=0 \, \text{ on } \Gamma_D\},
\]
FEM scheme becomes: (i) solve
\begin{equation}\label{eq: weak3a}
\begin{cases}
\text{Find } u_{h,nD}\in  \mathcal{P}_{h,D}^m,\,\ \text{s.t.} \\
        a(u_{h,nD},w_h)= \ell(w_h)-a(u_{h,D},w_h),\quad \forall w_h\in  \mathcal{P}_{h,D}^m,
\end{cases}
\end{equation}
and (ii)  set next
\[
u_{h}:= u_{h,D}+u_{h,nD}.
\]
The true unknowns are now  
\[
\{ u_j\}_{j\in{\tt inD}}
\]
which determine $u_{h,nD}$ and, by construction, the pointwise approximation of $u$ at the set of nodes of $\mathcal{P}_h^m$. 

By linearity, we can restrict the test to the elements of any  basis of $ \mathcal{P}_{h,D}^m$, such as $\{\varphi_j\}_{j\in {\tt inD}}$, resulting in  the scheme
\begin{equation}\label{eq: weak3}
\begin{cases}
\text{Find } 
\{ u_j\}_{j\in{\tt inD}}\subset \mathbb{R},\ \text{s.t.}\\
\displaystyle \sum_{j\in {\tt inD}}  a(\varphi_j,\varphi_i)u_j= \ell(\varphi_i)-\sum_{k\in {\tt iD}}a(\varphi_k,\varphi_i)u_D(\boldsymbol{x}_k),\qquad \forall i\in {\tt inD}.
\end{cases}
\end{equation}

This means that the approximation of the solution $u$ via the Finite Element Method transforms the problem to the resolution of the linear system
\begin{align}\label{eq: systemFEM}
    \nonumber \sum_{j \in {\tt inD}} \big( \boldsymbol{S}_{\underline{\boldsymbol{\kappa}},ij} 
    + \boldsymbol{A}_{\boldsymbol{\beta},ij} 
    &+ \boldsymbol{M}_{c,ij}
    + \boldsymbol{R}_{\alpha,ij} \big) u_{j}  
    = \boldsymbol{b}_{f,i} + \boldsymbol{t}_{r,i} \\
    &\quad - \sum_{k \in {\tt iD}} \big( \boldsymbol{S}_{\underline{\boldsymbol{\kappa}},ik}  
    + \boldsymbol{A}_{\boldsymbol{\beta},ik} 
    + \boldsymbol{M}_{c,ik} 
    + \boldsymbol{R}_{\alpha,ik}\big) u_{k}, \quad \forall i \in {\tt inD},
\end{align}
where  (see \eqref{eq: weak4}) 
\begin{equation}\label{eq:systemFEM2}
\begin{aligned}
\boldsymbol{S}_{\underline{\boldsymbol{\kappa}},ij} &:= \iiint_{\Omega} \nabla \varphi_i \cdot \left( \underline{\boldsymbol{\kappa}} (\nabla \varphi_j)^\top \right), &
\boldsymbol{R}_{\alpha,ij} &:= \iint_{\Gamma_{R}} \alpha \varphi_i \varphi_j, \\
\boldsymbol{M}_{c,ij} &:= \iiint_{\Omega} c \varphi_i \varphi_j, &
\boldsymbol{b}_{f,i} &:= \iiint_{\Omega} f \varphi_i, \\
\boldsymbol{A}_{\boldsymbol{\beta},ij} &:= \iiint_{\Omega} \varphi_i \big( \boldsymbol{\beta} \cdot \nabla \varphi_j \big), &
\boldsymbol{t}_{R,i} &:= \iint_{\Gamma_{R}} g_R \varphi_i.
\end{aligned}
\end{equation}
These terms are referred to as follows: in the first column, the stiffness matrix, the mass matrix, and the convection matrix. In the first term of the second column, the boundary mass matrix. The vectors in the second column are known as the load vector and the Robin vector, respectively.

\section{Assembly}\label{sec: assembly}

We next focus on the different vectors and matrices involved in the FEM, such as they are presented in \eqref{eq:systemFEM2}, which can be computed.

\subsection{Load vector. Why element-by-element implementation}

 Let us consider the load vector: for $i=1,\ldots, {\tt nNodes}$. Due to the linearity of the integral, we can decompose the integral over $\Omega$ to the tetrahedra containing each node:
\[
\begin{aligned}
 \iiint_{\Omega} f \varphi_i \  = \ &
 \iiint_{\mathop{\rm supp} \varphi_i} f \varphi_i  = \sum_{\substack{m \\ \boldsymbol{x}_i\in K_{m}}} \iiint_{K_{m}} f \varphi_i\\
= \ & \sum_{\substack{m \\ \boldsymbol{x}_i\in K_{m}}}|\det B_K| \iiint_{\widehat{K}} (f\circ F_{K_m}) (\varphi_i\circ F_{K_m})\\
 =\ & |\det B_K| \sum_{\substack{m \\ \boldsymbol{x}_i\in K_{m}}}  \iiint_{\widehat{K}} (f\circ F_{K_m})\widehat{N}_{r_{i,m}},\quad \text{with }F_{K_m}(\widehat{\boldsymbol{x}}_{r_{i,m}}) = \boldsymbol{x}_i.
 \end{aligned}
 \]
In other words, $i = {\tt ttrh }(m,r_{i,m})$ (cf. \eqref{eq:def:ttrh}), i.e.,  $r_{i,m}$ represents the local index of the node $\boldsymbol{x}_i$ within the element $K_m$.
 
The last integral can be computed or approximated by a quadrature rule.   A possible implementation is given in Algorithm \ref{alg:01}.

\begin{algorithm}[H]
     \For{$i=1:{\tt nNodes}$}{
      Find $\{m_1,m_2,\ldots,m_r\}$ s.t. $\boldsymbol{x}_{i}\in K_{m_r}$\;\label{number2}
      \For{$m\in\{m_1,m_2,\ldots,m_r\}$}{
       Find $r$, the local index of $\boldsymbol{x}_i$ in $K_{m}$\;\label{number4}
       Compute $b_{i}^{m} := |\det B_{K_m}|\iiint_{\widehat{K}} (f\circ F_{K_{m}})\widehat{N}_{r} $\;
      }
      $\boldsymbol{b}_{f,i}=\displaystyle\sum_{\ell =1}^{r} b_{\ell}^{m_\ell}$\;
    }
    \caption{\label{alg:01}Naive assembly for the load vector.}
\end{algorithm}
This approach requires two steps: first, searching for the elements to which node $\boldsymbol{x}_i$ belongs, and then performing a loop to determine the local ordering of the node within $K$. The first step is the most time-consuming, as it requires exploring all elements in the mesh for each node. 

A more efficient method can be employed, assembling element by element. In essence, this approach can be understood as the result of interchanging the two loops in the initial version of the algorithm:

\begin{algorithm}[H]
     \For{$m=1:{\tt nTtrh}$}{\label{line:01:alg02}
      \For{$r =0:{\tt dofK}-1$}{
       Find $i$ s.t. $F_{K_m}(\widehat{\boldsymbol{x}}_r) =\boldsymbol{x}_i$\;\label{line:03:alg02}
       $\boldsymbol{b}_{f,i}  =\boldsymbol{b}_{f,i}+ |\det B_{K_m}|\iiint_{\widehat{K}} (f\circ F_{K_{m}})\widehat{N}_{r}$\;
      }}
    \caption{\label{alg:02}Element-by-element assembly for the load vector.}
\end{algorithm}
Line\hspace{-0.7em}\ref{line:01:alg02} in Algorithm \ref{alg:02} serves as the foundation for the element-by-element computation by iterating through the rows of the element data  {\tt ttrh}. Line\hspace{-0.7em}\ref{line:03:alg02} retrieves the global indices of the nodes associated with the $m-$th element. Due to the structure of the {\tt ttrh} matrix, the second loop simplifies into a loop through its columns, allowing to obtain the global index simply as $i={\tt ttrh}(m,r)$.

\subsection{Robin vector}

For the Robin vector we need to compute integrals
\[
\iint_{A_m} g_R \varphi_i
\]
where $A_m$ is a triangle with vertices $\{\boldsymbol{x}_0^{A_m}, \boldsymbol{x}_1^{A_m}, \boldsymbol{x}_2^{A_m}\}$ and edges $\boldsymbol{\ell}^{A_m}_i=\boldsymbol{x}_1^{A_m}-\boldsymbol{x}_0^{A_m}$ ($i=1,2$)  contained in the (Robin) boundary of the domain.

Clearly, $\varphi_i|_{A_m}=0$ unless $\boldsymbol{x}_j = \boldsymbol{x}_r^{A_m}$ for some $r$. Then
\begin{equation}\label{eq:RobinVector}
\begin{aligned}
\iint_{A_m} g_R \varphi_i \ &=\  \|\boldsymbol{\ell}_1^{A_m}\times  \boldsymbol{\ell}^{A_m}_2\|
\iint_{\widehat{A}} g_R \circ F_{A_m} \varphi_i\circ F_{A_m}
\\
&=\  \|\boldsymbol{\ell}^{A_m}_1\times  \boldsymbol{\ell}^{A_m}_2\|  \begin{cases}
\displaystyle \iint_{\widehat{A}} g_R \circ F_{A_m}  \widehat{N}_r& \text{if }F_{A_m}(\widehat{\boldsymbol{x}}_r) = \boldsymbol{x}_i\\
0,&\text{otherwise}.
\end{cases}
\end{aligned}
\end{equation}

 If {\tt ntrBR} denotes the triangles on the Robin boundary, we have:

\begin{algorithm}[H]
     \For{$m=1:{\tt nTrBd}$}{
      \For{$r =0:{\tt dofA}-1$}{
       Find $i$ s.t. $F_{A_m}(\widehat{\boldsymbol{x}}_r) =\boldsymbol{x}_i$\;\label{line03:alg:03}
       $\boldsymbol{t}_{R,i}  =\boldsymbol{t}_{R,i}+ \|\boldsymbol{\ell}^{m}_1\times  \boldsymbol{\ell}^{m}_2\|\iint_{\widehat{A}} (g_R \circ F_{A_{m}})\widehat{N}_{r}$\;
      }}
    \caption{\label{alg:03}Element-by-element assembly for the Robin vector.}
\end{algorithm} 
As in the case of the mass matrix, the information in Line\hspace{-0.7em}\ref{line03:alg:03} is readily available from \eqref{eq:def:ttrh}.

\subsection{Mass matrices}\label{sec: Mass}\index{$\boldsymbol{M}_c$}

Let us consider in this subsection the computation of the mass and Robin matrices 
\[
\boldsymbol{M}_{c,ij}=\iiint_{\Omega}c\varphi_i \varphi_j, \quad 
\boldsymbol{R}_{\alpha,ij}=\iint_{\Gamma_R}\alpha \varphi_i \varphi_j.
\]
Any entry of the first matrix can be expressed as
 \begin{subequations}\label{eq:Mcij}
\begin{equation}\label{eq:Mcij:a}
 \boldsymbol{M}_{c,ij}=\sum_{\substack{ m\\ \boldsymbol{x}_i, \boldsymbol{x}_j\in {K_m}}} 
 |\det B_{K}|\iiint_{\widehat{K}_m} \widehat{c}_{K_m} \widehat{N}_{r_{i,m}}\widehat{N}_{r_{j,m}}, \quad  i,j =1,\ldots,{\tt nNodes}
 \end{equation}
 where in the sum above
\begin{equation}\label{eq:Mcij:b}
\widehat{c}_{K_m} = c\circ F_{K_m}, \quad F_{K_m}(\widehat{\boldsymbol{x}}_{r_{i,m}})=\boldsymbol{x}_{i}.
\end{equation}
\end{subequations} 
(Once again  $r_{i,m}$  is the local index of the $i-$th global node  in the element $K$).
Similarly to the load and Robin vector, we identify  in  \eqref{eq:Mcij} three nested loops: the first two iterate  over $i$ and $j$, running on the nodes, and the third, the innermost one, is over $m$, to locate the tetrahedron $K_m$ to which $\boldsymbol{x}_i$ and $\boldsymbol{x}_j$ belong.
Using an element-by-element assembly, which effectively shifts the innermost loop to the outermost position, reorganizes these computations and simplifies the process, as demonstrated in the following algorithm

\begin{algorithm}[H]
     \For{$m=1:{\tt nTrrh}$}{
      \For{$r =0:{\tt dofK}-1$}{
      Find $i$ s.t. $F_{K_m}(\widehat{\boldsymbol{x}}_{r}) =\boldsymbol{x}_i$\;\label{alg:04:03}
      \For{$s =0:{\tt dofF}-1$}{
      Find $j$ s.t. $F_{K_m}(\widehat{\boldsymbol{x}}_{s}) =\boldsymbol{x}_j$\;\label{alg:04:05}
       $\displaystyle \boldsymbol{M}_{c,ij}  =\boldsymbol{M}_{c,ij}+ |\det B_{K_m}|\iiint_{\widehat{K}}  \widehat{c}_{K_m} \widehat{N}_{s}\widehat{N}_{r}$\;
      }}}
    \caption{\label{alg:04}Element-by-element assembly for the Mass matrix.}
\end{algorithm}
Notice that once again Lines\hspace{-0.7em}\ref{alg:04:03} and\hspace{-0.7em} \ref{alg:04:05} are just the $r-$th and $s-$th column of {\tt ttrh}$(m,:)$, respectively. 

The Robin matrix can be handled in a very similar manner, reducing, after the change of variables, the assembly to compute the integrals element by element to
\begin{equation} \label{eq:BoundaryMassMatrix}
\|\boldsymbol{\ell}^{A_m}_1\times  \boldsymbol{\ell}^{A_m}_2\| \iint_{\widehat{A}} \widehat{\alpha}_{A_m}\widehat{N}_r\widehat{N}_s,\quad \widehat{\alpha}_{A_m} := \alpha\circ F_{A_m}.
\end{equation}
As a result, we obtain the following algorithm:

\begin{algorithm}[H]
     \For{$m=1:{\tt nTrBd}$}{
      \For{$r =0:{\tt dofA}-1$}{
       Find $i$ s.t. $F_{A_m}(\widehat{\boldsymbol{x}}_{r}) =\boldsymbol{x}_i$\;
      \For{$s =0:{\tt dofA}-1$}{
       Find $j$ s.t. $F_{A_m}(\widehat{\boldsymbol{x}}_{j}) =\boldsymbol{x}_j$\;
       $\displaystyle \boldsymbol{R}_{\alpha,ij}  =\boldsymbol{R}_{\alpha,ij}+  \|\boldsymbol{\ell}^{m}_1\times  \boldsymbol{\ell}^{m}_2\||\iiint_{\widehat{K}} \widehat{\alpha}_{A_m}  \widehat{N}_{s}\widehat{N}_{r}$\;
      }}}
    \caption{\label{alg:05}Element-by-element assembly for the Robin matrix.}
\end{algorithm}

Notice that in Algorithms \ref{alg:04} and \ref{alg:05} and  we are implicitly computing on each element $K$ the matrices 
\[
\boldsymbol{M}_{c}^K=   \begin{bmatrix}
    \displaystyle\int_K \widehat{c}_K  \widehat{N}_{0}\widehat{N}_{0}& \cdots &\displaystyle \int_K \widehat{c}_K  \widehat{N}_{0}\widehat{N}_{{\tt dofK}-1}\\
    \vdots & \ddots & \vdots \\
    \displaystyle\int_K \widehat{c}_K  \widehat{N}_{{\tt dofK}-1}\widehat{N}_{0}& \cdots & \displaystyle\int_K \widehat{c}_K  \widehat{N}_{{\tt dofK}-1}\widehat{N}_{{\tt dofK}-1}\\
\end{bmatrix},\quad \widehat{c}_K = c\circ F_K,
\]
\[
\boldsymbol{R}_{\alpha}^A= \begin{bmatrix}
    \displaystyle\int_A \widehat{\alpha}_A  \widehat{N}_{0}\widehat{N}_{0}& \cdots &\displaystyle \int_A \widehat{\alpha}_A  \widehat{N}_{0}\widehat{N}_{{\tt dofA}-1}\\
    \vdots & \ddots & \vdots \\
    \displaystyle\int_A \widehat{\alpha}_A  \widehat{N}_{{\tt dofA}-1}\widehat{N}_{0}& \cdots & \displaystyle\int_A \widehat{\alpha}_A  \widehat{N}_{{\tt dofA}-1}\widehat{N}_{{\tt dofA}-1} 
\end{bmatrix}, \quad \widehat{\alpha}_A = \alpha\circ F_A,
\]
usually referred as the local mass and Robin matrices.  

\subsection{Stiffness Matrix}\label{sec: Stiffness}

The assembly of the stiffness matrix is more involved as a consequence of combining the change of variable, which moves the computation to the reference element and the chain rule. Indeed, for a fixed $K$ and two nodes with indices $i,j$, the stiffness matrix reduces to:
\[
\int_K \nabla \varphi_i \underline{\boldsymbol{\kappa}}(\nabla\varphi_j)^\top = 
\int_K \nabla (\widehat{N}_r\circ F_K^{-1})  \underline{\boldsymbol{\kappa}}  (\nabla(\widehat{N}_s\circ F_K^{-1}))^\top 
\]
where 
\[
\varphi_i\circ F_K\lvert_K= \widehat{N}_r, \quad
\varphi_j\circ F_K\lvert_K= \widehat{N}_s.
\]
(That is, nodes $i$ and $j$ are the $r-$th and $s-$th in the local index of element $K$). The chain rule and the change of variables theorem yield then 
\[
\int_K \nabla \varphi_i \underline{\boldsymbol{\kappa}}(\nabla\varphi_j)^\top= \iiint_{\widehat{K}} \widehat{\nabla} \widehat{N}_r \underline{\widehat{\boldsymbol{\kappa}}}^K (\widehat{\nabla} \widehat{N}_s)^\top,\quad 
\underline{\widehat{\boldsymbol{\kappa}}}^K  = |\det B_K| B_K^{-1} (\underline{\boldsymbol{\kappa}}\circ F_K) B_K^{-\top}.
\]
\index{$\boldsymbol{S}_{\underline{\boldsymbol{\kappa}}}^K$}
Hence,  denoting  as usual $\underline{\widehat{\boldsymbol{\kappa}}} _{ij}^K$ the entries of the matrix function
$\underline{\widehat{\boldsymbol{\kappa}}}^K$  it holds
\begin{align*}
    \int_K \nabla \varphi_i \underline{\boldsymbol{\kappa}} (\nabla\varphi_j)^\top 
    =& \sum_{k,\ell=1}^3\iiint_{\hat{K}} \underline{\widehat{\boldsymbol{\kappa}}} _{k\ell}^K 
     \frac{\partial \widehat{N}_i}{\partial \widehat{x}_k} \frac{\partial \widehat{N}_j}{\partial \widehat{x}_\ell} 
\end{align*}
Therefore, the local stiffness matrix can be written in terms of nine local stiffness matrices  
 \[
 \boldsymbol{S}_{\underline{\boldsymbol{\kappa}}}^K = \sum_{k,\ell=1}^3
 \boldsymbol{S}_{\underline{\boldsymbol{\kappa}},\widehat{x}_k\widehat{x}_\ell}^K 
 \]
where 
\[
\boldsymbol{S}_{\underline{\boldsymbol{\kappa}},\widehat{x}_k\widehat{x}_\ell}^K = 
   \begin{bmatrix}
    \displaystyle\iiint_{\widehat{K}} \underline{\widehat{\boldsymbol{\kappa}}} _{k\ell}^K   \partial_{\widehat{x}_k}\widehat{N}_{0}\partial_{\widehat{x}_\ell}\widehat{N}_{0}& \cdots &\displaystyle \displaystyle\iiint_{\widehat{K}} \underline{\widehat{\boldsymbol{\kappa}}}_{k\ell}^K  \partial_{\widehat{x}_k}\widehat{N}_{0}\partial_{\widehat{x}_\ell}\widehat{N}_{{\tt dofK}-1}\\
    \vdots & \ddots & \vdots \\
    \displaystyle\iiint_{\widehat{K}} \underline{\widehat{\boldsymbol{\kappa}}} _{k\ell}^K \partial_{\widehat{x}_k} \widehat{N}_{{\tt dofK}-1}\partial_{\widehat{x}_\ell}\widehat{N}_{0}& \cdots & \displaystyle\iiint_{\widehat{K}} \underline{\widehat{\boldsymbol{\kappa}}} _{k\ell}^K  \partial_{\widehat{x}_k}\widehat{N} _{{\tt dofK}-1}\partial_{\widehat{x}_\ell}\widehat{N}_{{\tt dofK}-1}\\
\end{bmatrix}
\] 

The assembly process is shown in Algorithm  \ref{alg:06}.
The steps are now essentially the same. The only key difference is that more terms have to be computed for the assembly of the local matrix. The integrals are again computed using cubature rules. 

\begin{algorithm}[H]
     \For{$m=1:{\tt nTtrh}$}{
      \For{$r =0:{\tt dofK}-1$}{
      Find $i$ s.t. $F_{K_m}(\widehat{\boldsymbol{x}}_{r}) =\boldsymbol{x}_i$\;
      \For{$s =0:{\tt dofK}-1$}{
      Find $j$ s.t. $F_{K_m}(\widehat{\boldsymbol{x}}_{s}) =\boldsymbol{x}_j$\;
       $\boldsymbol{S}_{\underline{\boldsymbol{\kappa}},ij}  = \boldsymbol{S}_{\underline{\boldsymbol{\kappa}},ij} $ $+
             \big(\boldsymbol{S}_{\underline{\boldsymbol{\kappa}},\widehat{x}_1\widehat{x}_1}^{K_m}\big)_{rs}+
       \big(\boldsymbol{S}_{\underline{\boldsymbol{\kappa}},\widehat{x}_1\widehat{x}_2}^{K_m}\big)_{rs}$ $+       \big(\boldsymbol{S}_{\underline{\boldsymbol{\kappa}},\widehat{x}_1\widehat{x}_3}^{K_m}\big)_{rs}+       \big(\boldsymbol{S}_{\underline{\boldsymbol{\kappa}},\widehat{x}_2\widehat{x}_1}^{K_m}\big)_{rs}$ $ +
       \big(\boldsymbol{S}_{\underline{\boldsymbol{\kappa}},\widehat{x}_2\widehat{x}_2}^{K_m}\big)_{rs}+       \big(\boldsymbol{S}_{\underline{\boldsymbol{\kappa}},\widehat{x}_2\widehat{x}_3}^{K_m}\big)_{rs}$ $+
       \big(\boldsymbol{S}_{\underline{\boldsymbol{\kappa}},\widehat{x}_3\widehat{x}_1}^{K_m}\big)_{rs}+       \big(\boldsymbol{S}_{\underline{\boldsymbol{\kappa}},\widehat{x}_3\widehat{x}_2}^{K_m}\big)_{rs}$ $+
       \big(\boldsymbol{S}_{\underline{\boldsymbol{\kappa}},\widehat{x}_3\widehat{x}_3}^{K_m}\big)_{rs}$\;} 
      }}
    \caption{\label{alg:06}Element-by-element assembly for the Stiffness matrix.}
\end{algorithm}

\subsection{Advection Matrix}\label{sec: Adv}
 
Finally, we focus on the computation of 
\[
  \iiint_{K} \varphi_{i}  (\boldsymbol{\beta} \cdot\nabla \varphi_{j})
\]
which is reduced, via the Change of Variables and the chain rule as: 
\begin{multline*}
\lvert\det B_K\rvert \iiint_{\widehat{K}} \widehat{N}_r 
\left[\widehat{\nabla} \widehat{N}_s B_K^{-1} \widehat{\boldsymbol{\beta}}\right] \\
= \iiint_{\widehat{K}} \widehat{\beta}_1^K \widehat{N}_r \frac{\partial \widehat{N}_s}{\partial \widehat{x}} 
+ \iiint_{\widehat{K}} \widehat{\beta}_2^K \widehat{N}_r \frac{\partial \widehat{N}_s}{\partial \widehat{y}} 
+ \iiint_{\widehat{K}} \widehat{\beta}_3^K \widehat{N}_r \frac{\partial \widehat{N}_s}{\partial \widehat{z}}, 
\end{multline*}
where
\[
\widehat{\boldsymbol{\beta}}^K:=(\widehat{\beta}_1^K,\widehat{\beta}_2^K,\widehat{\beta}_3^K)^\top=B_K^{-1}(\boldsymbol{\beta}\circ F_K).
\]
The computation is then given:

\begin{algorithm}[H]
     \For{$m=1:{\tt nTtrh}$}{
      \For{$r =0:{\tt dofK}-1$}{
      Find $i$ s.t. $F_{K_m}(\widehat{\boldsymbol{x}}_{r}) =\boldsymbol{x}_i$\;
      \For{$s =0:{\tt dofK}-1$}{
      Find $j$ s.t. $F_{K_m}(\widehat{\boldsymbol{x}}_{s}) =\boldsymbol{x}_j$\;
      $\boldsymbol{A}_{\boldsymbol{\beta},ij}  = \boldsymbol{A}_{\boldsymbol{\beta},ij} $ $+\big(\boldsymbol{A}_{\boldsymbol{\beta},\widehat{x}_1}^{K_m}\big)_{rs}+\big(\boldsymbol{A}_{\boldsymbol{\beta},\widehat{x}_2}^{K_m}\big)_{rs}
       +\big(\boldsymbol{A}_{\boldsymbol{\beta},\widehat{x}_3}^{K_m}\big)_{rs}$
      }}}
    \caption{\label{alg:07}Element-by-element assembly for the Stiffness matrix.}
\end{algorithm}
Obviously, we are implicitly working with the local advection matrix $\boldsymbol{A}_{\boldsymbol\beta}^K$ that is written as 
\[
\boldsymbol{A}_{\boldsymbol\beta}^K =\boldsymbol{A}_{\boldsymbol\beta,\widehat{x}_1}^K+\boldsymbol{A}_{\boldsymbol\beta,\widehat{x}_2}^K+\boldsymbol{A}_{\boldsymbol\beta,\widehat{x}_3}^K
\]
with
\[
\boldsymbol{A}_{\boldsymbol{\beta},\widehat{x}_i}^{K}= 
   \begin{bmatrix}
 \displaystyle\iiint_{\widehat{K}} \widehat{\beta}_i^K \widehat{N}_0  \frac{\partial \widehat{N}_0}{\partial \widehat{x}_i} & \cdots &\displaystyle
  \iiint_{\widehat{K}} \widehat{\beta}_i^K \widehat{N}_0  \frac{\partial \widehat{N}_{{\tt dofK}-1}}{\partial \widehat{x}_i} \\
  \vdots & \ddots & \vdots \\ 
 \displaystyle\iiint_{\widehat{K}} \widehat{\beta}_i^K \widehat{N}_{{\tt dofK}-1} \frac{\partial \widehat{N}_0}{\partial \widehat{x}_i} & \cdots & 
  \displaystyle\iiint_{\widehat{K}} \widehat{\beta}_i^K \widehat{N}_{{\tt dofK}-1} \frac{\partial \widehat{N}_{{\tt dofK}-1}}{\partial \widehat{x}_i} \\
\end{bmatrix}.
\] 

\section{Key notes on FEM implementation in Matlab/Octave}\label{sec: data_structure}

This section deals with an implementation of the classical 3D $\mathbb{P}_m$-FEM algorithms with GMSH integration available in \cite{fem3d_gmsh}. While a detailed discussion of the implementation is beyond the scope of this work, we refer the reader to \cite{duque2023integration} for a more comprehensive explanation of the code.

\subsection{Importing and preprocessing the mesh}

We have implemented a function ({\tt importGMSH3D}) to process a simple GMSH mesh  exported in a {\tt *.msh} file under the {\tt 4.1} version (as pointed out before, the latest version available at the time of this publication) in ASCII format. We also compute some other quantities that are used or can be used in the future for the FEM implementation. We will provide a brief commentary on these aspects.

\subsubsection{Essential Boundary Conditions }
Section \ref{sec: assembly} focuses on the computation of the matrices and vectors in \eqref{eq:systemFEM2} from where \eqref{eq: systemFEM} can be easily constructed. Certainly, with
\[
\boldsymbol{C} =      
     \boldsymbol{S}_{\underline{\boldsymbol{\kappa}}} 
    + \boldsymbol{R}_{\alpha} 
    + \boldsymbol{A}_{\boldsymbol{\beta}} 
    + \boldsymbol{M}_{c},\quad {\boldsymbol d}  
    = \boldsymbol{b}_{f} + \boldsymbol{t}_{r},
\]
we obtain the system of equations:
\begin{verbatim}
C(inD, inD)u(inD) = d(inD) - C(inD, iD)u(iD),
\end{verbatim}
where {\tt iD} are the indices of Dirichlet nodes, for which the solution, i.e., ${\tt u}({\tt iD})$, is known. We need to establish a structure to classify elements and nodes based on their boundary conditions. This classification is implemented as follows: a container structure is used, named, for instance, {\tt ComplementaryInformation}, such that {\tt ComplementaryInformation}({\tt"PhysicalName"}) retrieves a vector with the numerical labels of the entities belonging to the specified physical group. 

The construction of this structure is straightforward. It involves storing the numerical and string labels of the Physical Groups mentioned in Section \ref{sec: PhysicalGroups}, looping through the entities (described in Section \ref{sec: Entities}), and subsequently looping through the Physical Groups associated with each entity.

We then create the vector {\tt domain} (resp. {\tt domBd}) of size ${\tt nTtrh} \times 1$ (resp. ${\tt nTrBd} \times 1$), such that {\tt domain(m)} (resp. {\tt domBd(m)}) contains the numerical label of the entity to which the element belongs.

If we construct a Physical Group {\tt "Dirichlet"} for the surface entities in $\Gamma_D$, the triangular elements belonging to this boundary can be retrieved using 
\begin{verbatim}
gammaD = ismember(domBd, ComplementaryInformation("Dirichlet"));
\end{verbatim}
The indices of the Dirichlet nodes, {\tt iD}, are obtained with
\begin{verbatim}
dirichletElements = trB(gammaD,:);
iD = unique(dirichletElements(:), :);
\end{verbatim}
making the non-Dirichlet nodes available with 
\begin{verbatim}
inD = setdiff(1:nNodes, iD);
\end{verbatim}

\subsubsection{Geometric mesh quantities}

As detailed in the previous section, several (geometric) quantities of the mesh are used in the assembly of the different matrices and vectors (cf.  \eqref{eq:systemFEM2}).  We give a sketch of this procedure emphasizing the vectorized approach which gives rise to a free-loop fast code. 

First the determinant of $B_{K_m}$, i.e. the volume of the tetrahedra.
 This data is stored in a ${\tt nTtrh}\times 1$ vector  {\tt detBk} so that ${\tt detBk}(m) = \text{det}(B_{K_m})$.  Matlab/Octave can compute all these quantities in a few lines by using the formula
\[
    \det B_{K_m} = \underbrace{(\boldsymbol{x}_1^{K_m}-\boldsymbol{x}_0^{K_m})}_{\boldsymbol{v}_{0,1}^{K_{m}}} \cdot \Big( \underbrace{(\boldsymbol{x}_2^{K_m}-\boldsymbol{x}_0^{K_m}}_{\boldsymbol{v}_{0,2}^{K_{m}}}) \times \underbrace{(\boldsymbol{x}_3^{K_m}-\boldsymbol{x}_0^{K_m}}_{\boldsymbol{v}_{0,3}^{K_{m}}})\Big)
\]
as shown here: 
\begin{programcode}{Determinant in Matlab/Octave}
\vspace{-0.7cm}\begin{verbatim}
v01 = coord(ttrh(:,2),:)-coord(ttrh(:,1),:);
v02 = coord(ttrh(:,3),:)-coord(ttrh(:,1),:);
v03 = coord(ttrh(:,4),:)-coord(ttrh(:,1),:);
detBk = dot(v01, cross(v02, v03, 2),2);
\end{verbatim}
\end{programcode}
Let us briefly comment these lines. The matrix {\tt coord} (see  \eqref{eq:coord}) contains the coordinates of the all the nodes of the mesh, while {\tt ttrh(:,i+1)} for ${\tt i} = 0,1, 2, 3$ returns the global indices of  the vertices ($\boldsymbol{x}_{{\tt i}}^K$) as a column vector.  

The first three lines compute ${\tt nTtrh} \times 3$ matrices, where each row $m$ contains ${\tt v}_{0,i}^{K_{m}}$, $i=1,2,3$, as described above. The functions {\tt dot} and {\tt cross} compute the dot and cross products of the row vectors of the matrices involved simultaneously, eliminating the need for explicit loops. Notice that the local indexing shown in Section \ref{sub:section:3.3} always makes the determinant positive.

A similar approach is followed for computing  the non-normalized normal vectors to the triangles in the boundary and stores them in a ${\tt ntrB}\times 3$ matrix {\tt Bdnormal}:
\begin{programcode}{Normal vector to triangle in the boundary in Matlab/Octave}
\vspace{-0.7cm}\begin{verbatim}
l1 = coord(trB(:,1),:) - coord(trB(:,2),:);
l2 = coord(trB(:,1),:) - coord(trB(:,3),:);
Bdnormal = cross(l1, l2, 2);
\end{verbatim}
\end{programcode}
Precomputing these normal vectors rather than their modules allows us to extend the code to cases where the Robin/Neumann data is not the scalar field $g_R$ but a vector field $\boldsymbol{g}_R$ such that $g_R= \boldsymbol{g}_R \cdot \widehat{\boldsymbol{n}}$ (for example, scattering problems).   

On the other hand, the module of these normal vectors is simply obtained using the command
\[
{\tt vecnorm(Bdnormal, 2, 2)}.
\]
Observe that this quantity is used in the construction of the Robin vector (cf.  \eqref{eq:RobinVector}) and the Boundary Mass matrix (cf.  \eqref{eq:BoundaryMassMatrix}).

\subsubsection{Faces and tetrahedra}
Other geometric quantities, such as: (i) the faces (including inner faces) of the tetrahedra, (ii) the specific faces that form the tetrahedra, and (iii) the tetrahedra to which these faces belong, can be obtained by processing ${\tt ttrh}$ using commands like ${\tt setdiff}$, ${\tt find}$, ${\tt sort}$ and ${\tt unique}$.  It is important to note that while this information is not utilized in our current implementation, it is necessary for calculating a posteriori FEM error estimates or for more advanced finite element methods, such as Discontinuous Galerkin schemes.

For illustrative purposes, the following code computes the faces (triangles) of the tetrahedra and store them in a  ${\tt nFaces \times 3}$ matrix {\tt faces}:

\begin{programcode}{Computation of tetrahedron interfaces}
\vspace{-0.7cm}\begin{verbatim}
indAux   = [2 3 4; 1 4 3; 1 2 4; 1 3 2]; 
faces  = [ttrh(:,indAux(1,:)); ...
          ttrh(:,indAux(2,:)); ...
          ttrh(:,indAux(3,:)); ...
          ttrh(:,indAux(4,:))];
facesAux = sort(faces(:,1:3),2);  
[~,p1, p3]  = unique(facesAux,'rows', 'first');
faces     = faces(p1,:); 
\end{verbatim}
\end{programcode}
The variable {\tt indAux} establishes an ordering for the vertices of all interfaces, ensuring that the normal vectors point outwards to the elements. We then define all faces in the mesh through their vertices. With this, it is just a matter of finding rows with a unique set of indices (up to permutation). We remove the permutation requirement by sorting out the index of the faces in a new variable {\tt facesAux}. Finally, we find the first occurrence of all the interfaces to construct {\tt faces}.

Next, the relation between tetrahedra and their faces and vice versa can be set with the following code:
 
\begin{programcode}{Relationship between tetrahedrons and interfaces}
\vspace{-0.7cm}\begin{verbatim}
ind = zeros(nTtrh,4); ind(:) = 1:numel(ttrh(:,1:4));
ttrh2faces  = p3(ind);
[~,p2]      = unique(facesAux,'rows', 'last');
p1  = mod(p1,nTtrh)    ; p1(p1==0) = nTtrh;   
p2  = mod(p2,nTtrh)    ; p2(p2==0) = nTtrh;  
faces2ttrh = [p1, p2];
ind = p1==p2;
faces2ttrh(ind,2) = nan;
\end{verbatim}
\end{programcode} 
We construct the ${\tt nTtrh \times}4$ matrix {\tt ttrh2faces} and the ${\tt nFaces \times}2$ matrix {\tt faces2ttrh}. The former matrix contains, on each row, the indices of the four faces in {\tt faces} belonging to a tetrahedron. The latter contains at each row the indices of the two tetrahedra (or one tetrahedron if the face lies in the boundary) in {\tt ttrh}.

The variable {\tt ind} enumerates the faces bounding every tetrahedron ({\tt ind(i,j)} is the {\tt (i+j)}-th row of {\tt facesAux}, i.e., the {\tt j}-th face belonging to the {\tt i}-th tetrahedron). The index {\tt p3} defined in the previous computation contains the mapping of the indices of the unique faces in {\tt faces} and all the faces. Hence we can construct {\tt ttrh2faces} using {\tt p3(ind)}. 

Notice that {\tt facesAux} is a $4\cdot {\tt nTtrh \times} 3$ matrix. Each {\tt nTtrh} block contains one of the faces of every tetrahedron in the mesh. With {\tt mod}, we restore the index of the tetrahedron associated with that face (except for the special case of the last tetrahedron, which we readjust to {\tt nTtrh}). The pair of tetrahedra containing each face are the indices of the first and last appearance of a face in {\tt facesAux}. To make boundary faces easily accessible and avoid redundant information, we set the second value to {\tt nan} by looking for the cases where the tetrahedron indices repeat in a row. We store this information in {\tt faces2ttrh}.

\subsection{Vector and matrix assembly}
Let us consider the assembly of the load vector. To achieve this, we apply a quadrature formula using appropriate weights $\{\omega_\ell\}_{\ell=1}^{\tt nQ}$ and nodes $\{\widehat{\boldsymbol{\lambda}}_\ell\}_{\ell=1}^{\tt nQ}$, provided in barycentric coordinates within the reference element $\widehat{K}$, such that
\[
\iiint_{\widehat{K}} (f\circ F_{K_{m}})\widehat{N}_{r} \approx \frac{1}{6}\det B_{K_m} \sum_{\ell=1}^{{\tt nQ}} \omega_\ell f\left(\begin{bmatrix}
    \boldsymbol{x}_0^{K_m}& \boldsymbol{x}_1^{K_m}&\boldsymbol{x}_2^{K_m}&\boldsymbol{x}_3^{K_m}
\end{bmatrix} \widehat{\boldsymbol{\lambda}}_\ell \right)
\widehat{N}_r(\widehat{\boldsymbol{\lambda}}_\ell).
\]
(Notice that, for the sake of convenience, we have written  $\widehat{N}_r$ in terms of barycentric coordinates instead of Cartesian variables, as they were originally introduced.)
The computation of these terms in Matlab can be implemented as follows:

\begin{programcode}{Quadrature formula approach for computing load vector}
\vspace{-0.7cm}\begin{verbatim}
Nr = zeros(dofK, nQ);
for i = 1:dofK
    Nr(i,:) = Nj3D{i}(nodesQuad(1,:),nodesQuad(2,:),...
                            nodesQuad(3,:),nodesQuad(4,:)); 
end
px = coord(:,1); py = coord(:,2); pz = coord(:,3);
nodesX = px(ttrh(:,1:4))*nodesQuad;
nodesY = py(ttrh(:,1:4))*nodesQuad;
nodesZ = pz(ttrh(:,1:4))*nodesQuad;

val = f(nodesX, nodesY, nodesZ).*detBk;
Nr = repmat(Nr,nTtrh,1); 
val = repelem(val,dofK,1); 
val = val.*Nr;
val = val*weights(:)/6;
\end{verbatim}
\end{programcode}

Let us describe the code briefly since it illustrates the techniques we use to assemble the right-hand-side vector. 

To compute the integral, we need to evaluate both the shape functions and the source term $f$ at the quadrature nodes, which are specified in barycentric coordinates in the ${\tt nQ}\times 4$ matrix {\tt nodesQuad}. The shape functions have been previously stored in the cell structure {\tt Nj3D} as four-variable functions of the barycentric coordinates. This evaluation is performed in the first four lines of the previous code. 

To evaluate $f$, we first compute the quadrature nodes for all elements of the mesh, i.e., the action of the affine mapping $F_K$ on the quadrature nodes in the reference element via  \eqref{eq: baryc_mapping} in Lines 6-9. The expression {\tt px(ttrh(:,1:4))} produces a matrix where each row contains the $x$-coordinate of the vertices of all tetrahedra. As a result, {\tt nodesX} is an ${\tt nTtrh} \times {\tt nQ}$ matrix that stores, row by row, the $x$-coordinate of the quadrature nodes within all elements. A similar structure applies to {\tt nodesY} and {\tt nodesZ} for the y- and z-coordinates, respectively.

Lines 11–14 evaluate $(f \widehat{N}_i) \circ F_K$ at the quadrature nodes for each element $K$. Note that a column-wise matrix-vector operation ({\tt ".*"}) is applied twice: first, to incorporate ${\tt detBK}$ into the values of $f$, and second, to multiply these values by the values of all shape functions at the quadrature nodes. 

The final line computes the action of the quadrature rule and stores the results in ${\tt nTtrh}$ blocks, each of size ${\tt dofK}$. These blocks are then concatenated into a vector of length ${\tt nTtrh} \cdot {\tt dofK}$.

Let us emphasize that this code can be easily modified to handle functions $f$ with an extra (fourth) variable to indicate the physical group to which each element belongs. This facilitates the straightforward definition of piecewise-defined functions for the load vector.

The assembly of the load vector can be implemented in a vectorised way in just three lines:
\begin{programcode}{Assembly of load vector}
\vspace{-0.7cm}
\begin{verbatim}
indi = T.ttrh'; 
indi = indi(:); 
b =  accumarray(indi, val, [nNodes, 1]);  \end{verbatim}
\end{programcode}
The validity of these lines is far from trivial. The construction of {\tt indi} is somewhat delicate, requiring careful attention to order the data stored in {\tt val}. The function {\tt accumarray} facilitates the construction of vectors by inserting {\tt val(i)} (recall Matlab uses column-major order) at the position {\tt indi(i)}, adding up values in cases of repeating indices. These lines complete the computation of the inner loop in Algorithm \ref{alg:02} in a vectorized manner.

This process is extended to compute the matrices in the system of equations, following similar steps. For instance, for the mass matrix, we proceed as follows:  
(a) compute $\widehat{N}_i \widehat{N}_j$ at the quadrature nodes and store the results in a ${\tt dofK}^2 \times {\tt nQ}$ matrix;  
(b) evaluate the function $c$ at the quadrature nodes for all tetrahedra  and store the result in a ${\tt nTtrh}\times {\tt nQ}$ vector precisely in the same manner as $f$ was evaluated in the previous code;  
(c) use {\tt repmat} and {\tt repelem}, column-wise product ``${\tt .*}$" to construct a $({\tt dofK}^2 \cdot {\tt nTtrh}) \times {\tt nQ}$ matrix, which is then multiplied by ${\tt weights}$ and ${\tt detBk}$ to compute the action of the quadrature rule (this matrix contains the required quantities conveniently stored);  
(d) assemble the matrix using {\tt sparse} (in place of {\tt accumarray}) with suitable indices vector {\tt indi},{\tt indj}.

We also need to compute $B_{K_m}^{-1}$ for the stiffness and advection matrices. We can exploit vectorization using the following representation of the inverse matrix:
\[
B^{-1}_{K_m}=\frac{1}{\det{B_{K_m}}}\begin{bmatrix}
\left[(\boldsymbol{x}_2^{K_m}-\boldsymbol{x}_0^{K_m})\times (\boldsymbol{x}_3^{K_m}-\boldsymbol{x}_0^{K_m})\right]^\top\\
\left[(\boldsymbol{x}_3^{K_m}-\boldsymbol{x}_0^{K_m})\times (\boldsymbol{x}_1^{K_m}-\boldsymbol{x}_0^{K_m})\right]^\top\\
\left[(\boldsymbol{x}_1^{K_m}-\boldsymbol{x}_0^{K_m})\times (\boldsymbol{x}_2^{K_m}-\boldsymbol{x}_0^{K_m})\right]^\top
\end{bmatrix},
\]
which results in
\begin{programcode}{Computation of coefficients of inverse matrices}
\vspace{-0.7cm}\begin{verbatim}
v01 = coord(ttrh(:,2),:)-coord(ttrh(:,1),:);
v02 = coord(ttrh(:,3),:)-coord(ttrh(:,1),:);
v03 = coord(ttrh(:,4),:)-coord(ttrh(:,1),:);
b1 = cross(v02,v03,2)./detBk;
b2 = cross(v03,v01,2)./detBk;
b3 = cross(v01,v02,2)./detBk;
\end{verbatim}
\end{programcode}
Matrices ${\tt nTtrh}\times 3$ {\tt b1}, {\tt b2} and {\tt b3} contain, row-wise, the first, second, and third column of $B_K$ for all elements $K$,

\subsection{Postprocessing of the FEM solution} \label{sub:postprocess}
After solving the linear system, we may want to evaluate the solution at arbitrary points. 
This requires, for each $\boldsymbol{x}_0$, identifying first the element $K$ containing $\boldsymbol{x}_0$ and evaluating next the polynomial  $u_h|_K \in\mathbb{P}_m$. The implementation  follows these steps: 
\begin{enumerate}
    \item Compute the barycentric coordinates $\boldsymbol{\lambda}^K(\boldsymbol{x_0})$ for all elements $K$ (or a subset of them where $\boldsymbol{x}_0$ is known to be).
    \item Find the element $j$ of nodes $\{i_1,i_2,\ldots, i_{\tt dofK}\}$ for which all barycentric coordinates $\boldsymbol{\lambda}^{K_j}(\boldsymbol{x}_0)$ are non-negative.
    \item Compute
    \[
    u_h (\boldsymbol{x}_0) = \sum_{r=0}^{\tt dofK} u_{i_r} \widehat{N}_r(\boldsymbol{\lambda}^{K_j}(\boldsymbol{x}_{0})).
    \]
\end{enumerate}
Step three can be implemented as a dot-product as follows: Set ${\boldsymbol m}$ the $1 \times {\tt nNodes}$ vector with 
\[
m_i = \begin{cases}
\widehat{N}_r(\boldsymbol{\lambda}^{K_j}(\boldsymbol{x}_0)), & \text{if $i=i_r$}, \qquad r=1, \ldots, {\tt dofK},\\
0 & \text{otherwise}.
\end{cases}
\]
Then
\[
u_h(\boldsymbol{x}_0)= {\boldsymbol m}\cdot {\boldsymbol u}.
\]
For multiple ${\tt nP}$ points, such an evaluation can be performed by stacking these vectors into an ${\tt nP} \times {\tt nNodes}$ matrix, ${\boldsymbol M}$. Note that ${\boldsymbol M}$ is sparse. Once constructed, it can be reused to evaluate multiple solutions at the same points by simply performing the matrix-vector multiplication ${\boldsymbol M}\boldsymbol{u}$.

This approach can be implemented with manageable resource requirements for a modest number of points. However, if the solution needs to be evaluated at a large number of points, the memory requirements for performing this computation in a vectorized manner increase significantly. This challenge can be addressed if some a priori knowledge about the locations of the points is available.

\section{Numerical examples}

We present numerical examples across a range of different phenomena to demonstrate the implementation of our software \cite{fem3d_gmsh}. The first two numerical simulations were performed on a MacBook Pro equipped with an Apple M3 chip, 16 GB of RAM, and running macOS Sonoma 14.5. The third experiment, however, was conducted on a Linux machine running CentOS Linux release 7.3.1611  with processors Intel(R) Xeon(R) CPU E5-2620 v4 @ 2.10GHz  and 128 GB of RAM. All linear systems were solved using direct methods: either with the ``\verb!\!'' command  ({\tt mldivide}), or, in the third experiment, through an $LDL^\top$ factorization of the matrix (combined with suitable sorting algorithm to minimize the fill-in effect in the sparse matrix).

\subsection{Heat equation}
\begin{figure}
    \centering
    \begin{subfigure}{0.33\textwidth}
        \centering
        \includegraphics[scale=0.29]{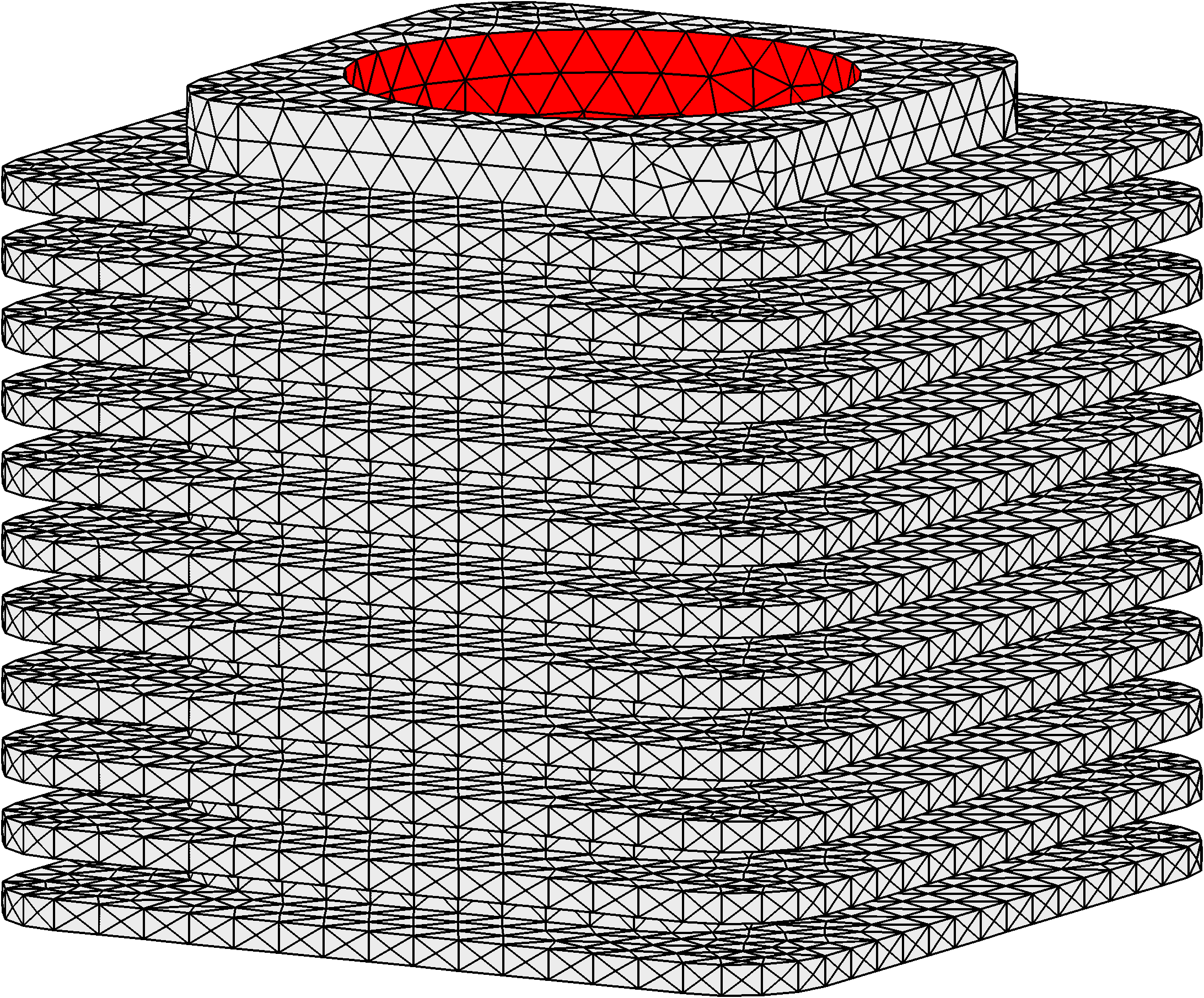}
        \caption{Finned cylinder \label{fig:geom_cylinder}}
    \end{subfigure}\hfill
    \begin{subfigure}{0.33\textwidth}
        \centering
        \includegraphics[scale=0.29]{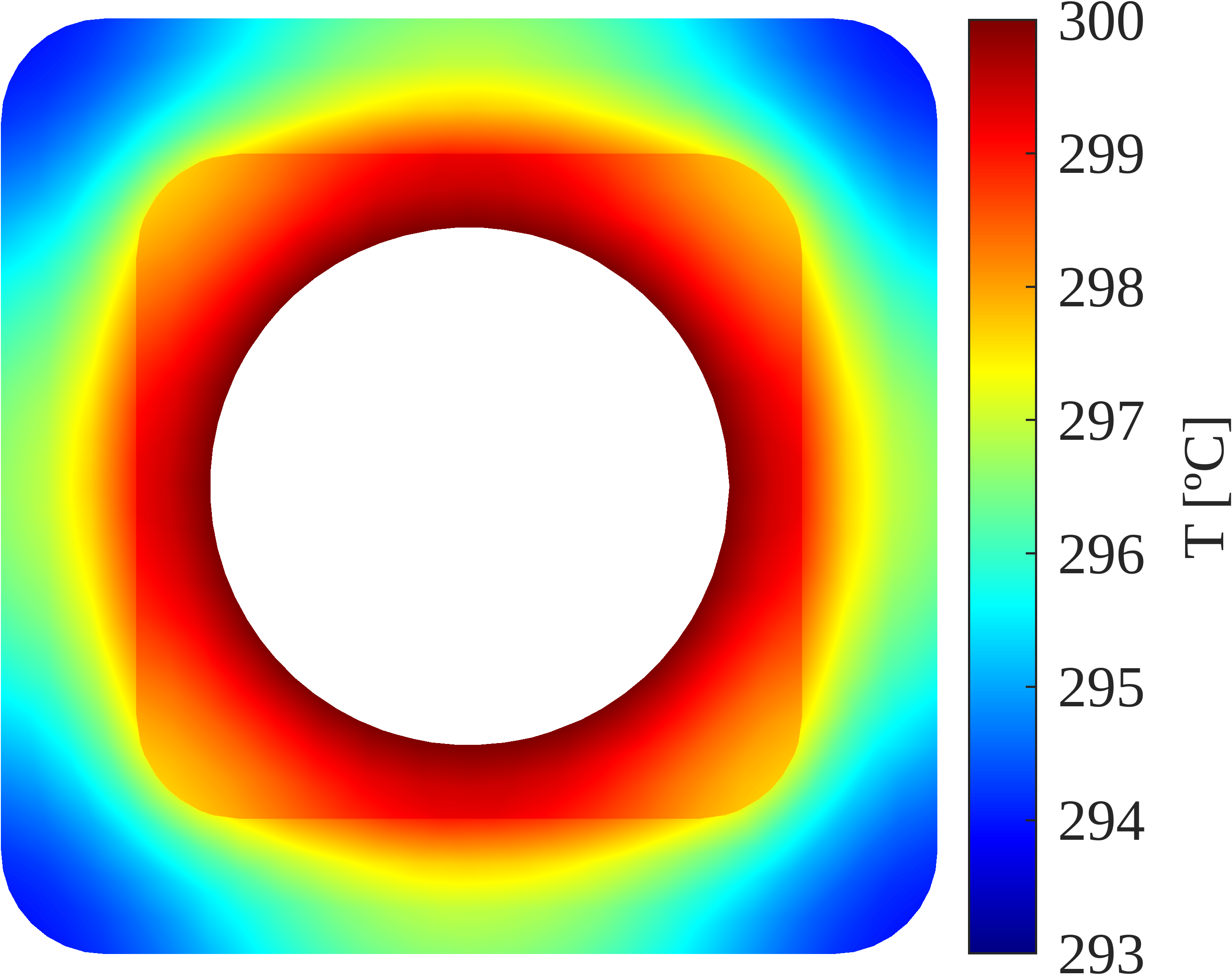}
        \caption{Top view}
    \end{subfigure}\hfill
    \begin{subfigure}{0.33\textwidth}
        \centering
        \includegraphics[scale=0.29]{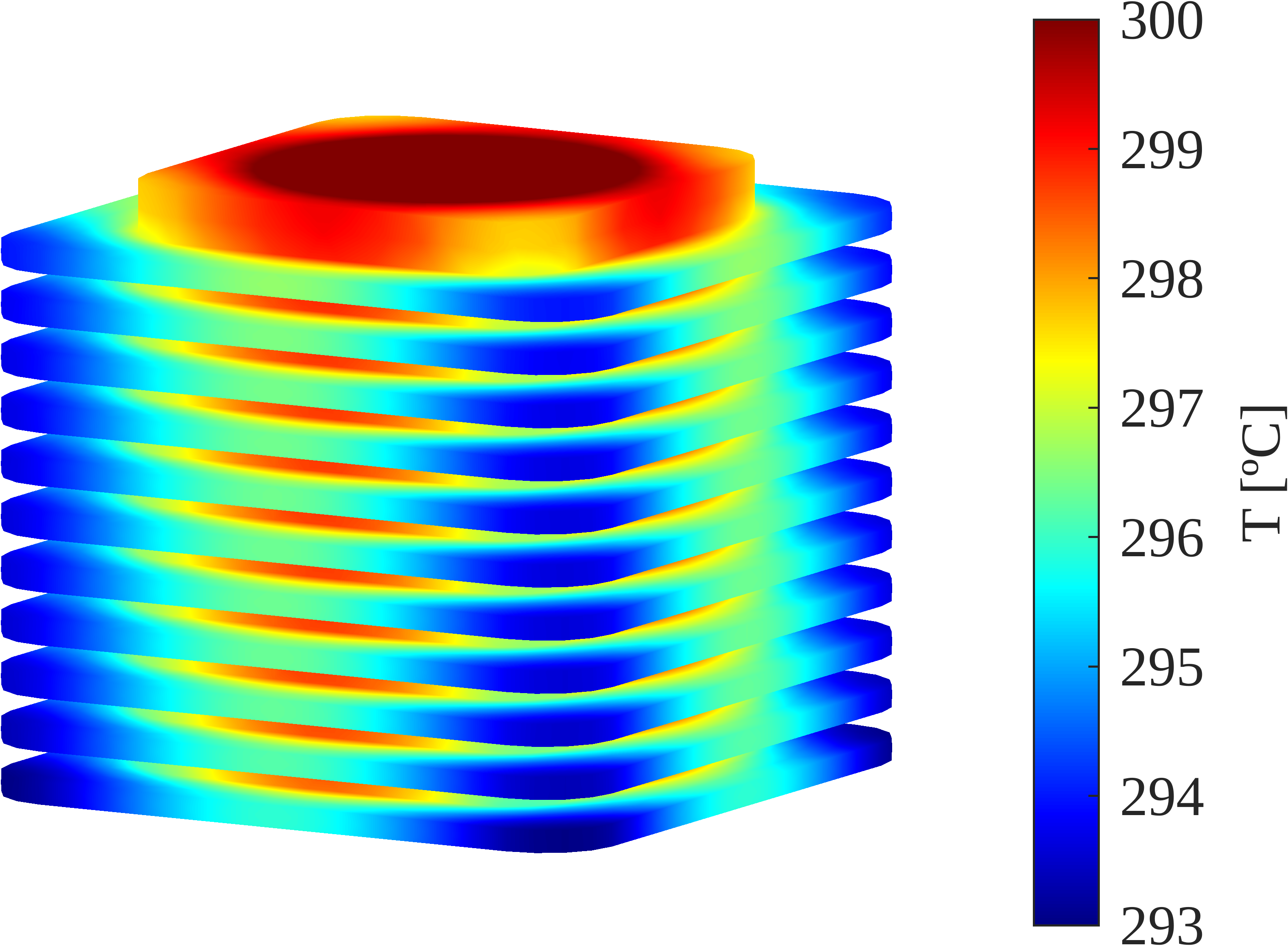}
        \caption{Perspective view}
    \end{subfigure}
    \caption{Heat equation in the steady state.}
    \label{fig:steady_state}
\end{figure}

We first solve the heat equation
\[
    \rho C_p\frac{\partial T}{\partial t}- \kappa \Delta T = 0,
\]
over an aluminum-finned engine cylinder (see Figure \ref{fig:geom_cylinder}). The properties of this material are: a density $\rho = 7200\,  \text{kg/m}^3$, heat capacity $C_p=447 \, \text{J}/^\circ\text{C}$, and thermal conductivity $\kappa = 52\, \text{W}/^\circ\text{C}$.

 Such a complex model has been constructed using free CAD software, FreeCAD \cite{freecad},  exported in STL format, and processed in GMSH.

We use $\mathbb{P}_1$ elements for this simulation. Constructing the entire data structure takes 0.82 seconds. The mesh comprises 28,979 tetrahedral elements, 17,322 triangular elements on the boundary, and 9,229 nodes. Given the entity-based structure of the {\tt *.msh} file and the nature of the pre-processing computations, the cost of this function primarily depends on the number of elements. It is, therefore, remarkable how efficient vectorized code is in Matlab.

We consider an initial temperature over the whole body of $T_0=20^\circ$C except for the Dirichlet condition (in red) with a constant temperature $T_D= 300^\circ$C. For the rest of the surfaces, we consider that there is dissipation due to convection:
\[
q_{\rm conv} = h(T-T_\infty),
\]
where $T_\infty=27^\circ$C is the temperature of the surrounding air and $h=5 \,\text{W/m}^2\,^\circ$C is the convection coefficient. This is implemented as a Robin condition by equalizing this term with the heat flux at the Robin boundary. 

The stiffness matrix is computed exactly since $\kappa$ is constant. The computation of the stiffness matrix and the Robin mass matrix and the Robin vector takes 0.03 s for their assembly. Figure \ref{fig:steady_state}b-c the solution to the steady state. In this case, we notice that neither $\rho$ and $C_p$ play a role in heat diffusion. All boils down to a balance of the diffusion term with the convection term at the boundary.  With this solution, we see the importance of fins in bodies. They allow extending the heat flow region, hence reducing the overall maximum temperature by increasing the surface of contact with the surrounding air.

\begin{figure}
    \centering
    \begin{subfigure}{0.33\textwidth}
        \centering
        \includegraphics[scale=0.29]{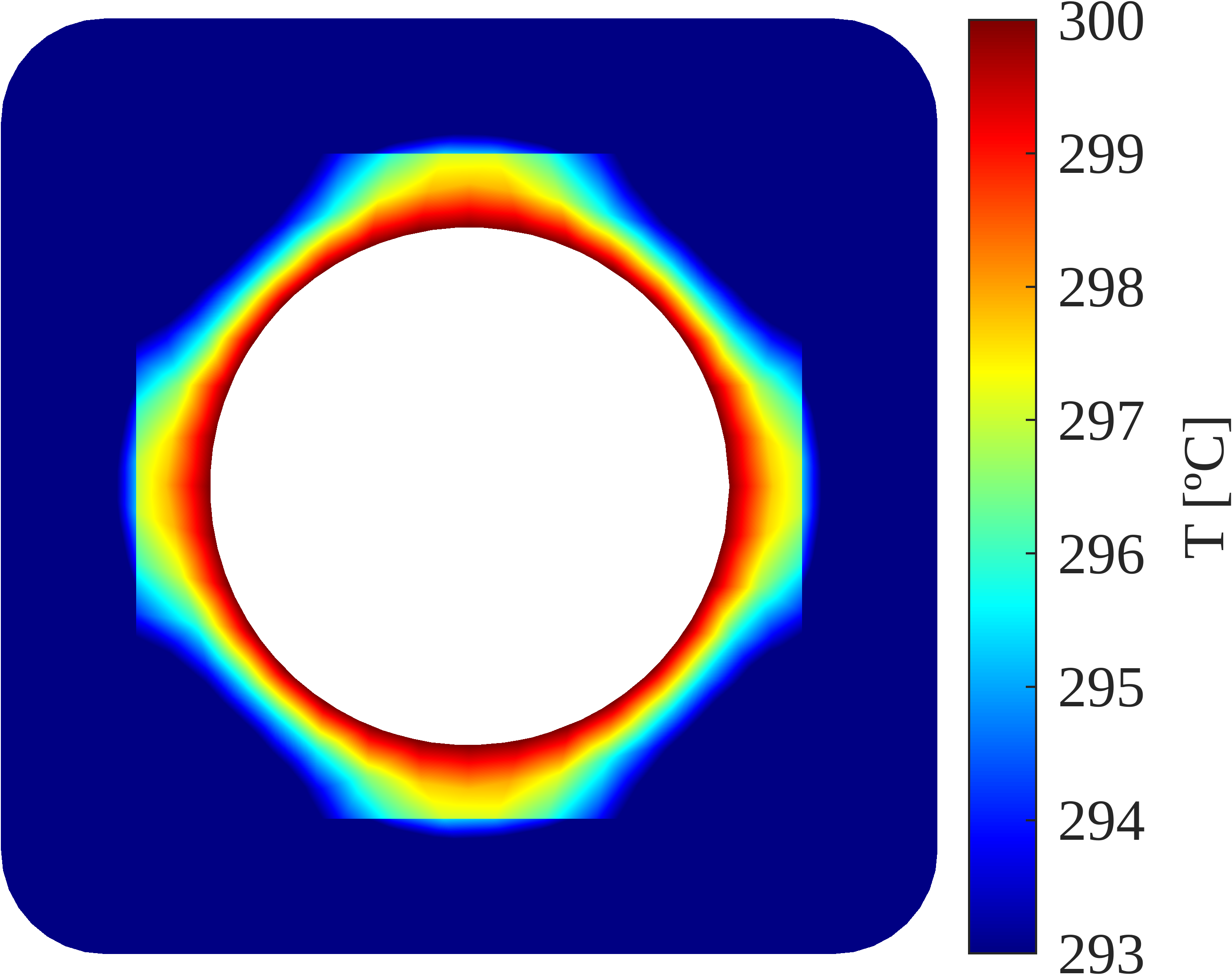}
        \caption{Top view,\\
        $t=40$ s}
    \end{subfigure}\hfill
    \begin{subfigure}{0.33\textwidth}
        \centering
        \includegraphics[scale=0.29]{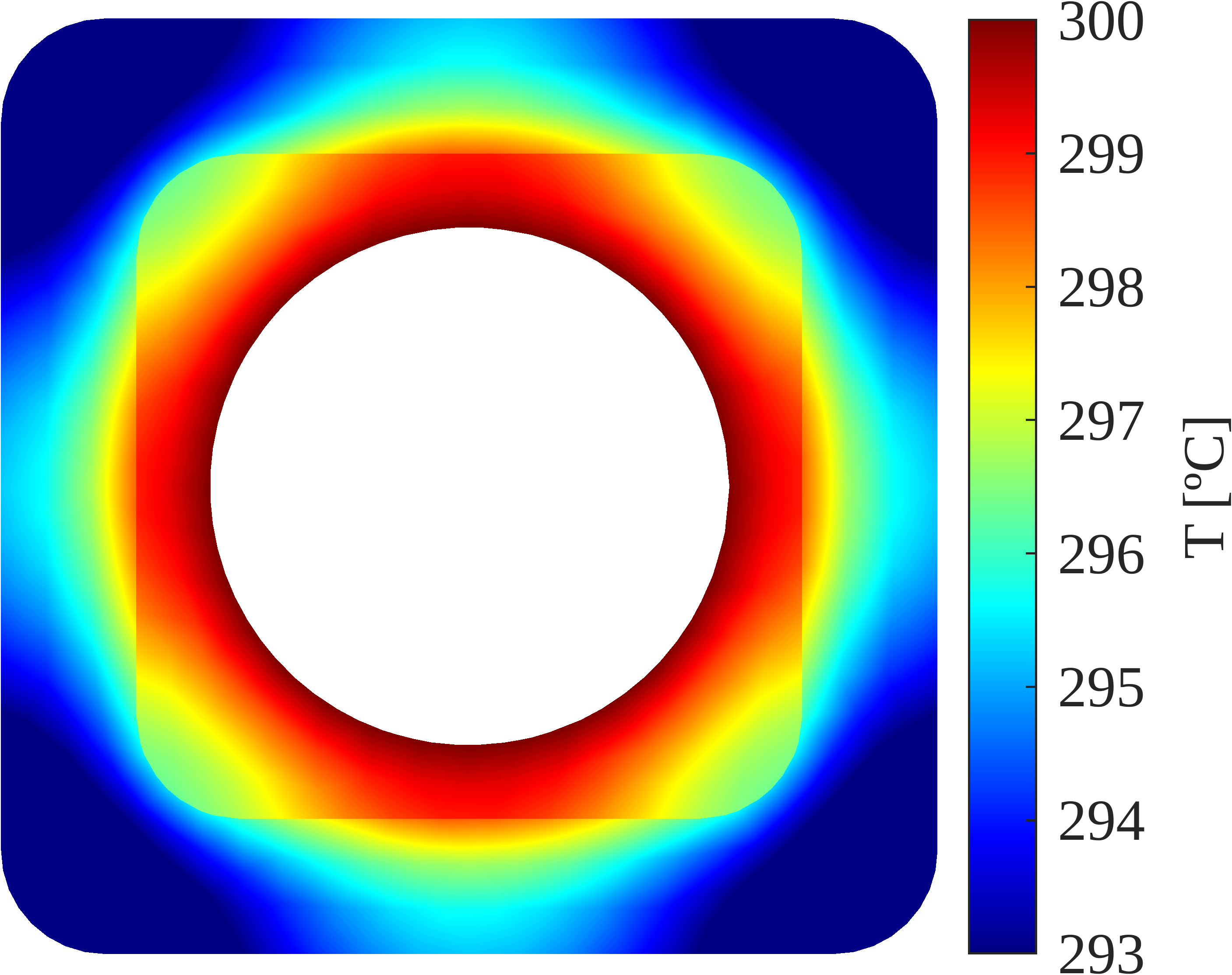}
        \caption{Top view,\\
        $t=70$ s}
    \end{subfigure}\hfill
    \begin{subfigure}{0.33\textwidth}
        \centering
        \includegraphics[scale=0.29]{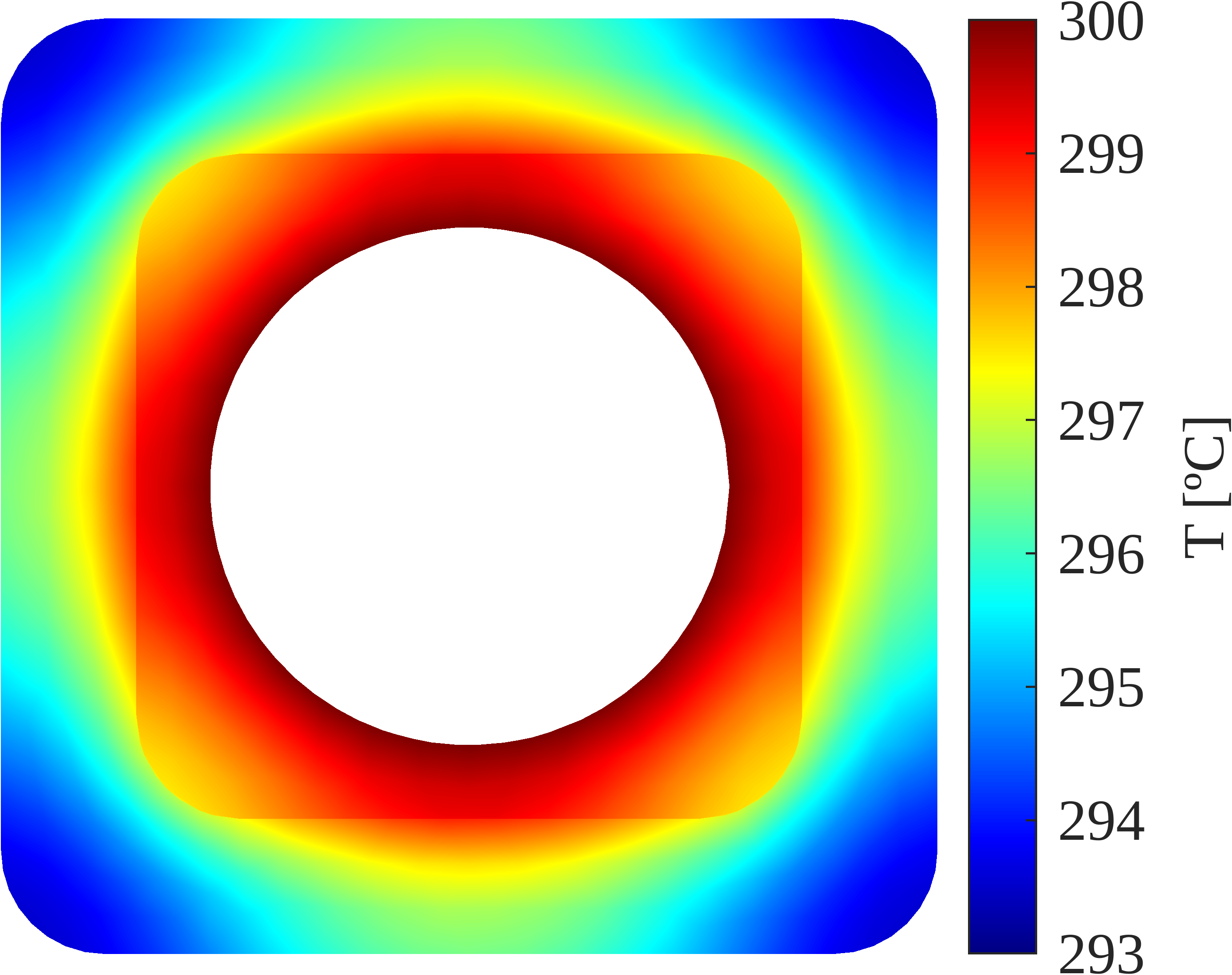}
        \caption{Top view,\\
        $t=100$ s}
    \end{subfigure}
    \\
    \begin{subfigure}{0.33\textwidth}
        \centering
        \includegraphics[scale=0.29]{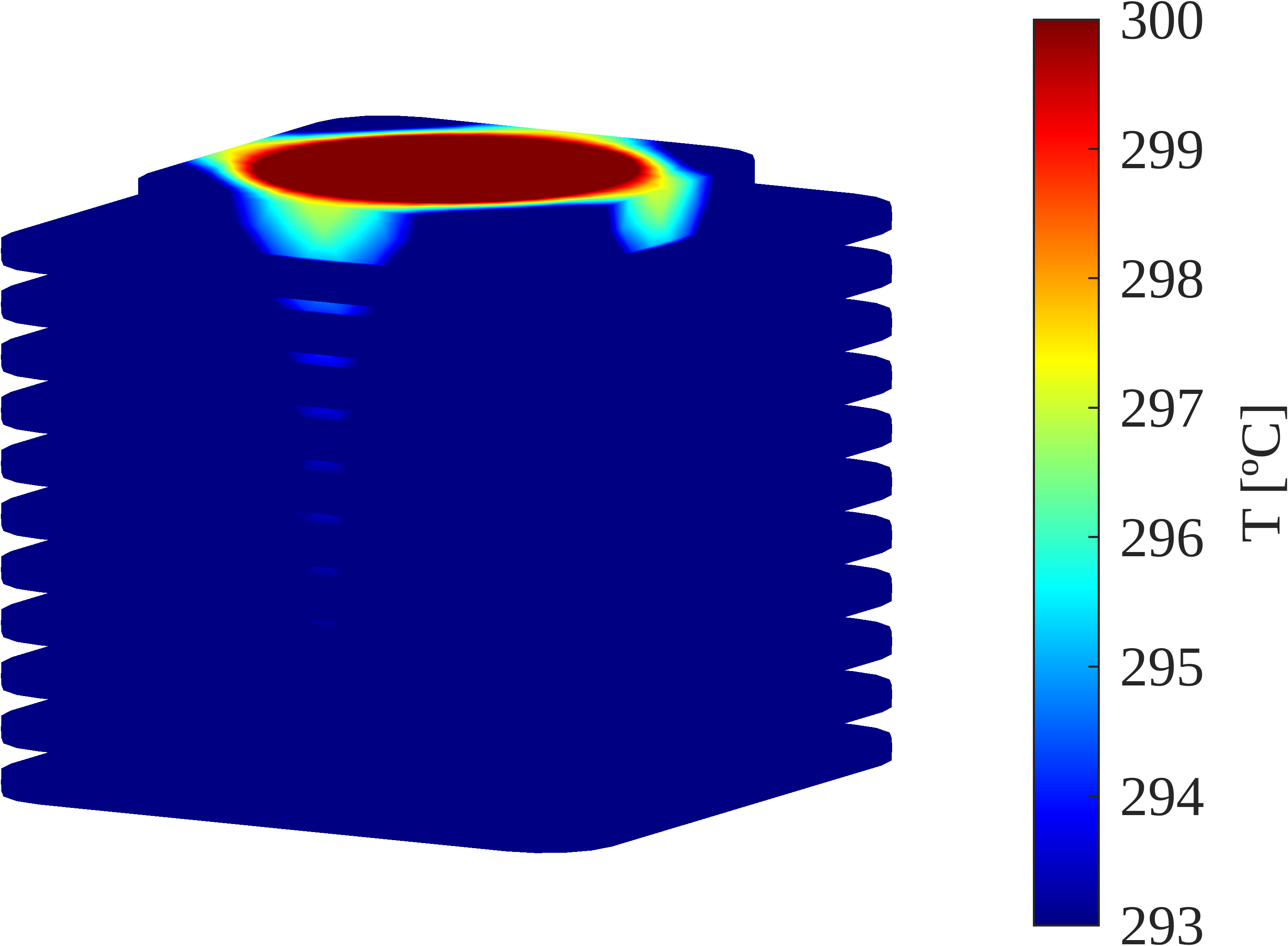}
        \caption{Perspective view,\\
        $t=40$ s}
    \end{subfigure}\hfill\begin{subfigure}{0.33\textwidth}
        \centering
        \includegraphics[scale=0.29]{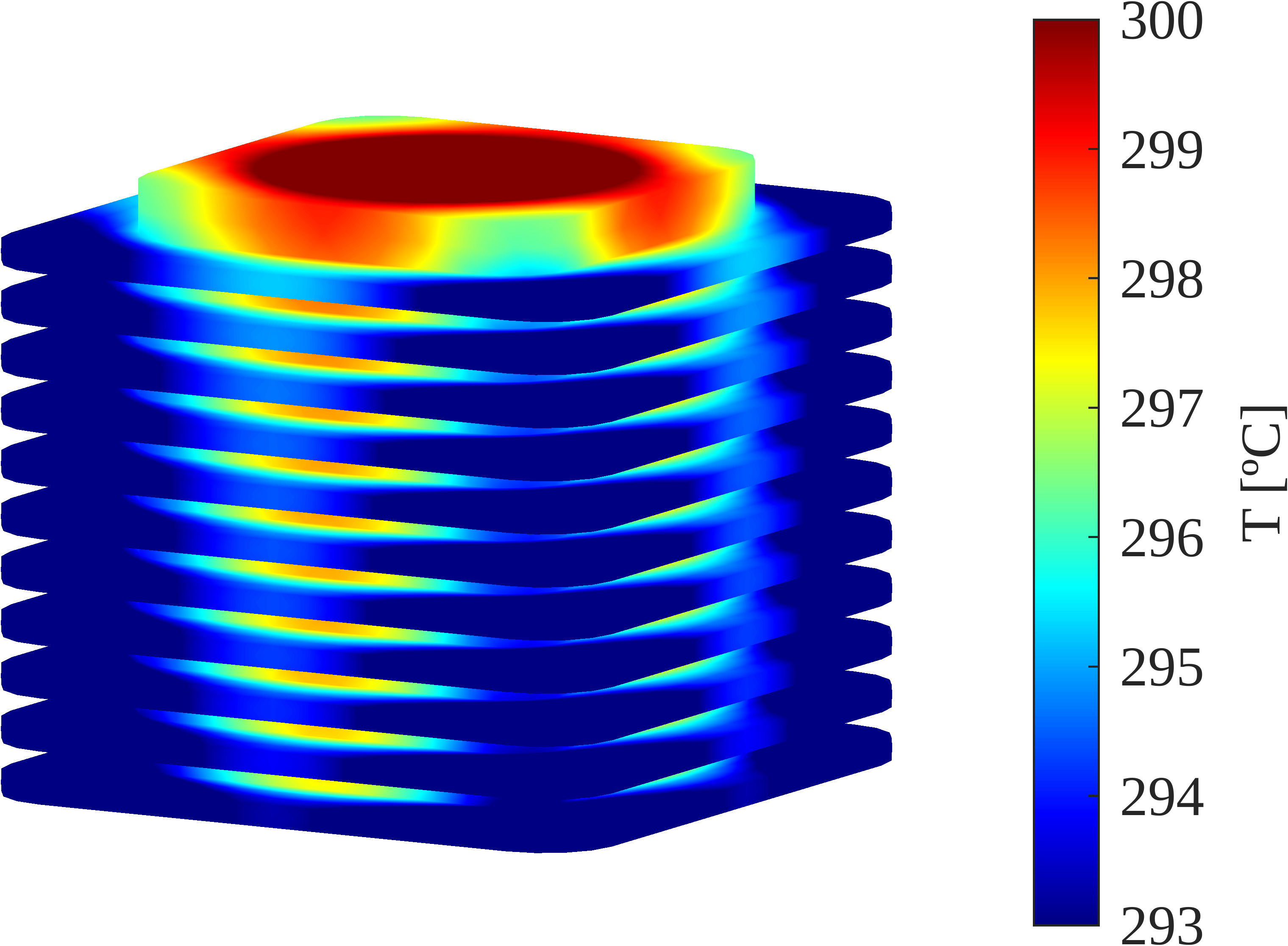}
        \caption{Perspective view,\\
        $t=70$ s}
    \end{subfigure}
    \begin{subfigure}{0.33\textwidth}
        \centering
        \includegraphics[scale=0.29]{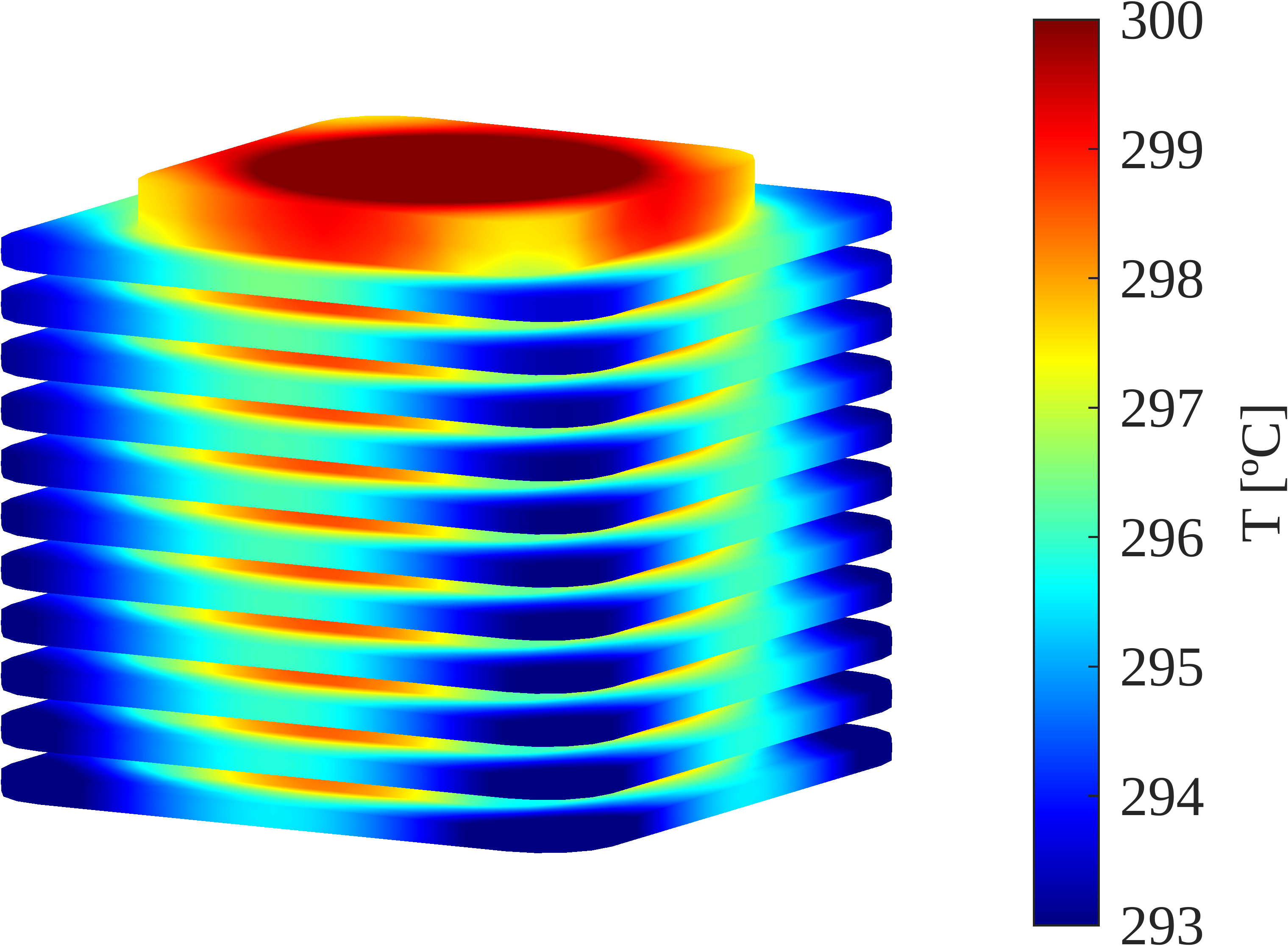}
        \caption{Perspective view,\\
        $t=100$ s}
    \end{subfigure}
    \caption{Transient solution.}
    \label{fig:transient_heatEq}
\end{figure}

We also perform a transient simulation to see the evolution of the heat flow through the cylinder. The assembly of the mass matrix takes around 0.026 s. We integrate from $t=0$ to $t=120$ s.  We use a Crank-Nicolson scheme to obtain the same second-order spatial and temporal discretization convergence. We choose a time step $\tau= 0.05$ s.  In total, it takes 4 minutes and 50 seconds to solve around 2400 time steps. The results of the numerical simulations are given in Figure \ref{fig:transient_heatEq}.

The temperatures shown in the plots are bounded by the extreme values in the steady state to represent how each point in the cylinder reaches its steady state temperature.
\subsection{Eigenvalue problem in elasticity}

 To demonstrate the capabilities and flexibility of this package, we also illustrate how it can be extended to other problems, such as the linear elasticity problem. This problem consist in finding the deformation field  $\boldsymbol{u}:\mathbb{R}^3\to \mathbb{R}$ satisfying
\[\label{eq: elasticityProblem}
   \left\{\begin{array}{cccl}
      - \nabla\cdot\underline{\boldsymbol\sigma} \left(\boldsymbol{u}\right) &=& \boldsymbol{f}, & \boldsymbol{x}\in\Omega,  \\
\underline{\boldsymbol\sigma} \left(\boldsymbol{u}\right)\widehat{\boldsymbol{n}} &=& \boldsymbol{g}_R(\boldsymbol{x}), & \boldsymbol{x}\in\Gamma_R, \\
        \boldsymbol{u}(\boldsymbol{x}) &=& \boldsymbol{u}_D(\boldsymbol{x}), & \boldsymbol{x} \in \Gamma_D, 
   \end{array} \right.
\]
where
\[\label{eq: constEq}\underline{\boldsymbol\sigma}\left(\boldsymbol{u}\right) = 2\mu  \underline{\boldsymbol\varepsilon}\left(\boldsymbol{u}\right) +\lambda \mathcal{I}_3\mathop{\rm tr}\,  \underline{\boldsymbol\varepsilon}\left(\boldsymbol{u}\right),
\]
is the stress tensor. Here, $\mathcal{I}_3$ is the $3\times 3$ identity matrix, $\mu$ and $\lambda$ are the Lamé constants and
\[
 \underline{\boldsymbol\varepsilon} = \frac{1}{2}\left( \nabla \boldsymbol{u} + \nabla \boldsymbol{u}^\top\right),
\]
is the linearised deformation tensor. 

The basis used to describe the finite element space is given by $\{\boldsymbol\varphi_j^n\}_{j=1,\ldots,{\tt nNodes}}^{n=1,2,3}$, where $\boldsymbol\varphi_j^n = \varphi_j\widehat{\boldsymbol{e}}_n$, with $\widehat{\boldsymbol{e}}_n$ the unit vector in the $n-$th direction.  Ignoring the Dirichlet condition, and integrating by parts, the weak formulation can be written as
\[
    \int_\Omega \underline{\boldsymbol\sigma}\left(\boldsymbol{u}\right) \boldsymbol{:}  \underline{\boldsymbol\varepsilon}\left(\boldsymbol\varphi_i^n\right)-\int_{\Gamma_R} \underbrace{\left(\underline{\boldsymbol\sigma} \left(\boldsymbol{u}\right)\widehat{\boldsymbol{n}} \right)\cdot \boldsymbol\varphi_i^n}_{=\boldsymbol{g}{_R }\cdot \boldsymbol\varphi_i^n}=\int_\Omega \boldsymbol{f}\cdot \boldsymbol\varphi_i^n, \qquad \forall i=1,\ldots,{\tt nNodes}, \,n=1,2,3,
\]
where $\boldsymbol{:}$ denotes the usual Frobenius product. The solution is approximated as 
\[
\boldsymbol{u} \approx \boldsymbol{u}_h = \sum_{j=1}^{\tt nNodes} \left( u_{x,j}\boldsymbol\varphi_j^{1}+ u_{y,j}\boldsymbol\varphi_j^{2}+ u_{z,j}\boldsymbol\varphi_j^{3}\right).
\]
We can easily show that for $n_1,n_2=1,2,3$:
\[
\begin{aligned}
    \underline{\boldsymbol\sigma}\left(\boldsymbol\varphi_{j}^{n_1}\right)\boldsymbol{:} \underline{\boldsymbol\varepsilon}\left(\boldsymbol\varphi_{i}^{n_2}\right)&=\mu\left[\delta_{n_1n_2}(\nabla \varphi_j\cdot \nabla\varphi_i)+\partial_{\widehat{\boldsymbol{e}}_{n_2}}\varphi_j\partial_{\widehat{\boldsymbol{e}}_{n_1}}\varphi_i\right]+\partial_{\widehat{\boldsymbol{e}}_{n_1}}\varphi_j\partial_{\widehat{\boldsymbol{e}}_{n_2}}\varphi_i, \\&= \nabla\varphi_j \underline{\boldsymbol{\kappa}}_{n_1,n_2}(\nabla\varphi_i)^\top.
\end{aligned}
\]
 The discretized system of equations become
\[
\begin{bmatrix}
        \boldsymbol{S}_{11} & \boldsymbol{S}_{12} & \boldsymbol{S}_{13}\\
        \boldsymbol{S}_{12}^\top & \boldsymbol{S}_{22} & \boldsymbol{S}_{23}\\
        \boldsymbol{S}_{13}^\top & \boldsymbol{S}_{23}^\top & \boldsymbol{S}_{33}\\
    \end{bmatrix} \begin{bmatrix}
        \boldsymbol{u}_x \\ \boldsymbol{u}_y\\ \boldsymbol{u_z} 
    \end{bmatrix}= \begin{bmatrix}
        \boldsymbol{b}_{\boldsymbol{f}\cdot \widehat{\boldsymbol{e}}_1} + \boldsymbol{t}_{\boldsymbol{g}_R\cdot \widehat{\boldsymbol{e}}_1}\\
        \boldsymbol{b}_{\boldsymbol{f}\cdot \widehat{\boldsymbol{e}}_2} + \boldsymbol{t}_{\boldsymbol{g}_R\cdot \widehat{\boldsymbol{e}}_2}\\
        \boldsymbol{b}_{\boldsymbol{f}\cdot \widehat{\boldsymbol{e}}_3} + \boldsymbol{t}_{\boldsymbol{g}_R\cdot \widehat{\boldsymbol{e}}_3}
    \end{bmatrix},
\]
where $\boldsymbol{S}_{n_1,n_2}$ is the stiffness matrix with diffusion matrix $\underline{\boldsymbol{\kappa}}_{n_1,n_2}$.

We solve the eigenvalue problem in elasticity. This corresponds to finding the natural frequencies of the geometry.  We use $\mathbb{P}_2$-elements since  $\mathbb{P}_1-$FEM requires a large number of elements to see convergence. 

\begin{figure}
    \centering
    \begin{subfigure}{0.33\textwidth}
        \centering
        \includegraphics[scale=0.5]{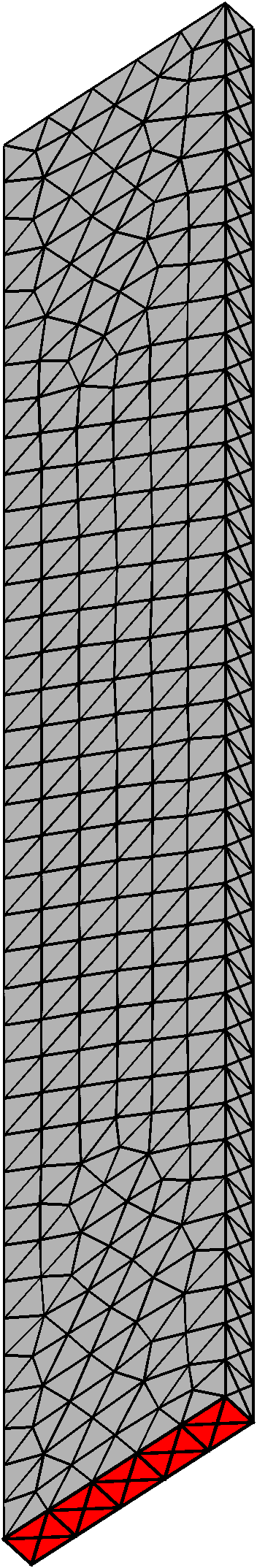}
        \caption{Beam \\geometry \label{fig:geom}}
    \end{subfigure}\hfill
    \begin{subfigure}{0.33\textwidth}
        \centering
        \includegraphics[scale=0.5]{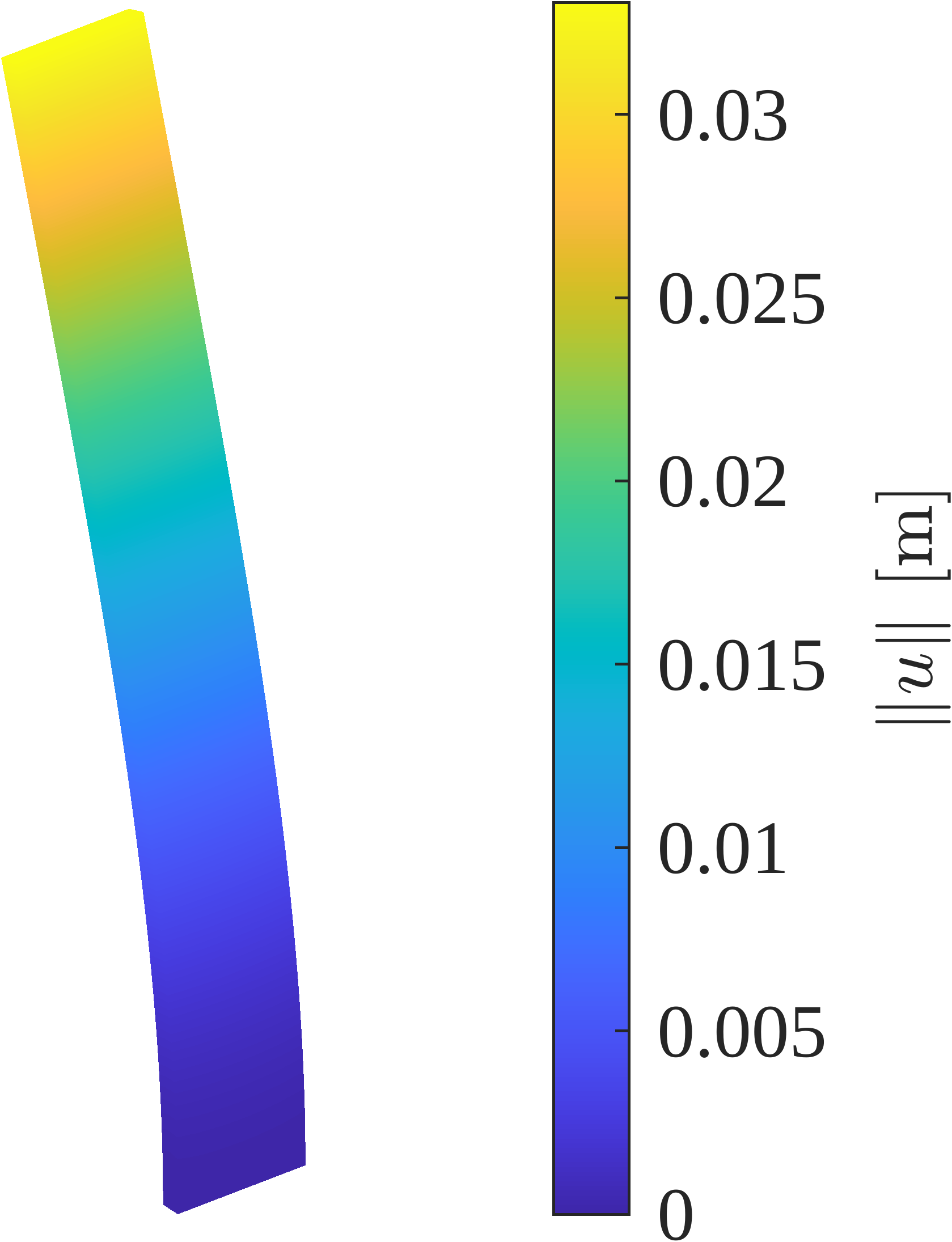}
        \caption{First mode, \\$f_1 = 112.4$ Hz}
    \end{subfigure}\hfill
    \begin{subfigure}{0.33\textwidth}
        \centering
        \includegraphics[scale=0.5]{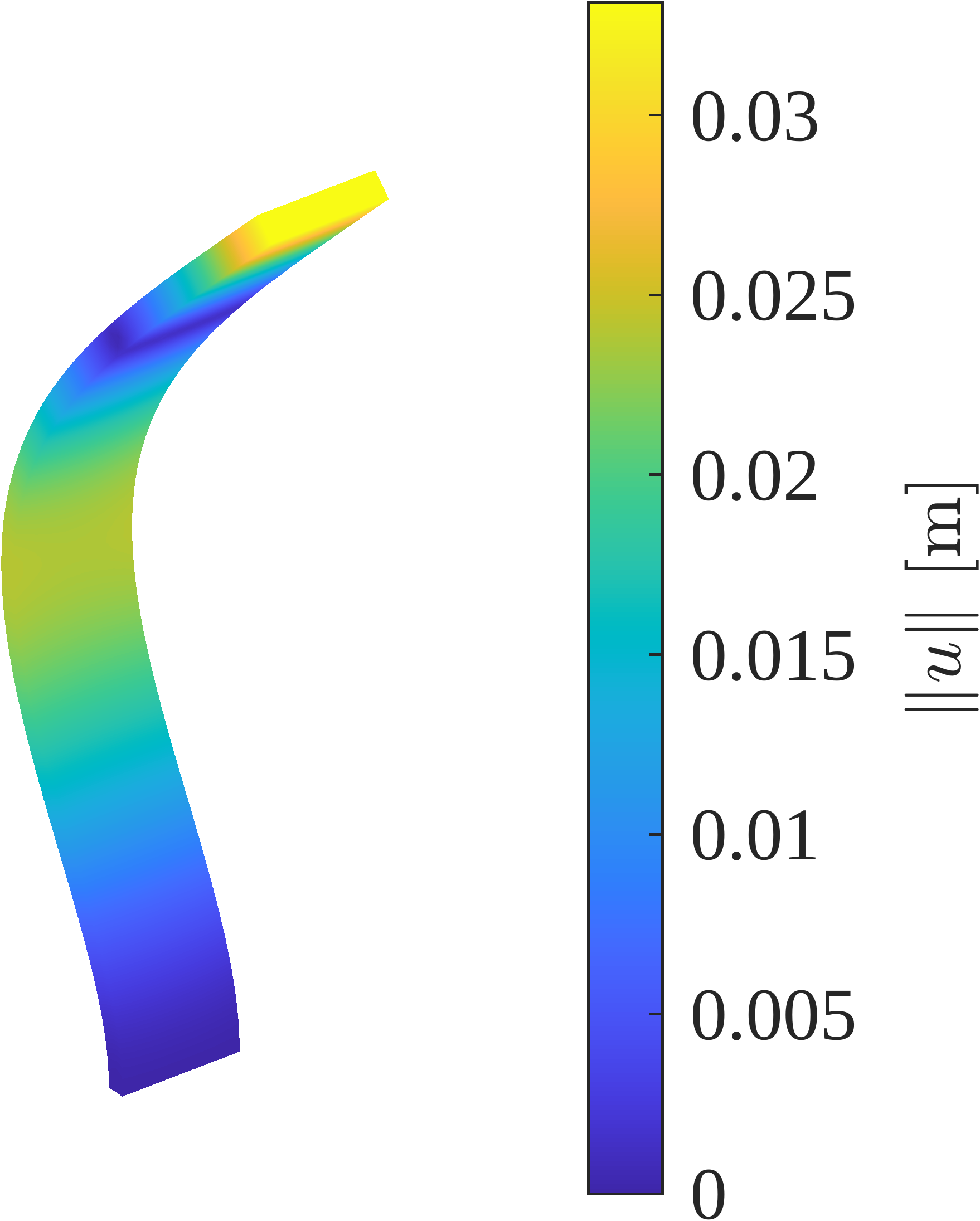}
        \caption{Second mode,\\ $f_2=703.15$ Hz}
    \end{subfigure}
    \\
    \begin{subfigure}{0.33\textwidth}
        \centering
        \includegraphics[scale=0.5]{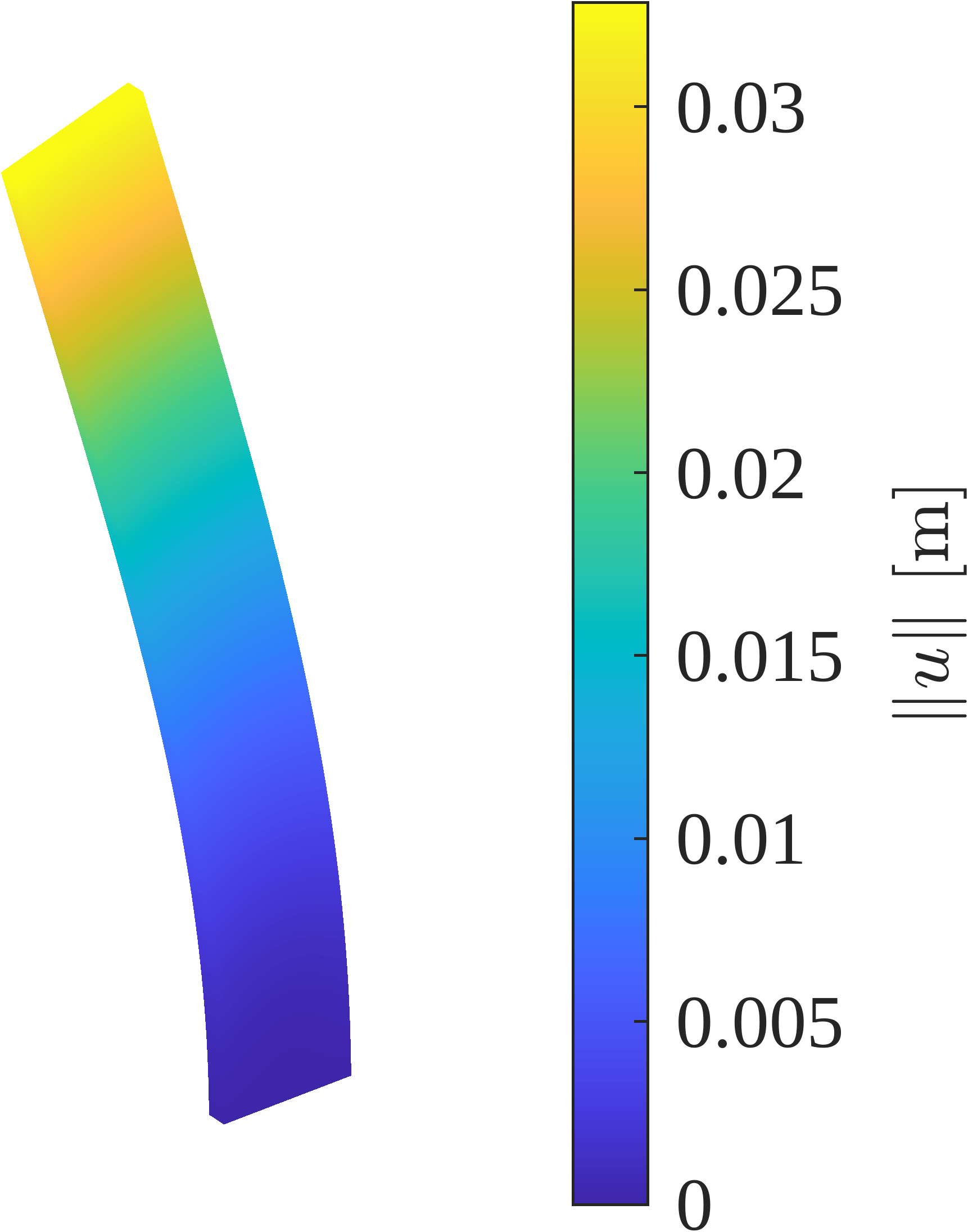}
        \caption{Third mode,\\ $f_3=734.71$ Hz}
    \end{subfigure}\hfill\begin{subfigure}{0.33\textwidth}
        \centering
        \includegraphics[scale=0.5]{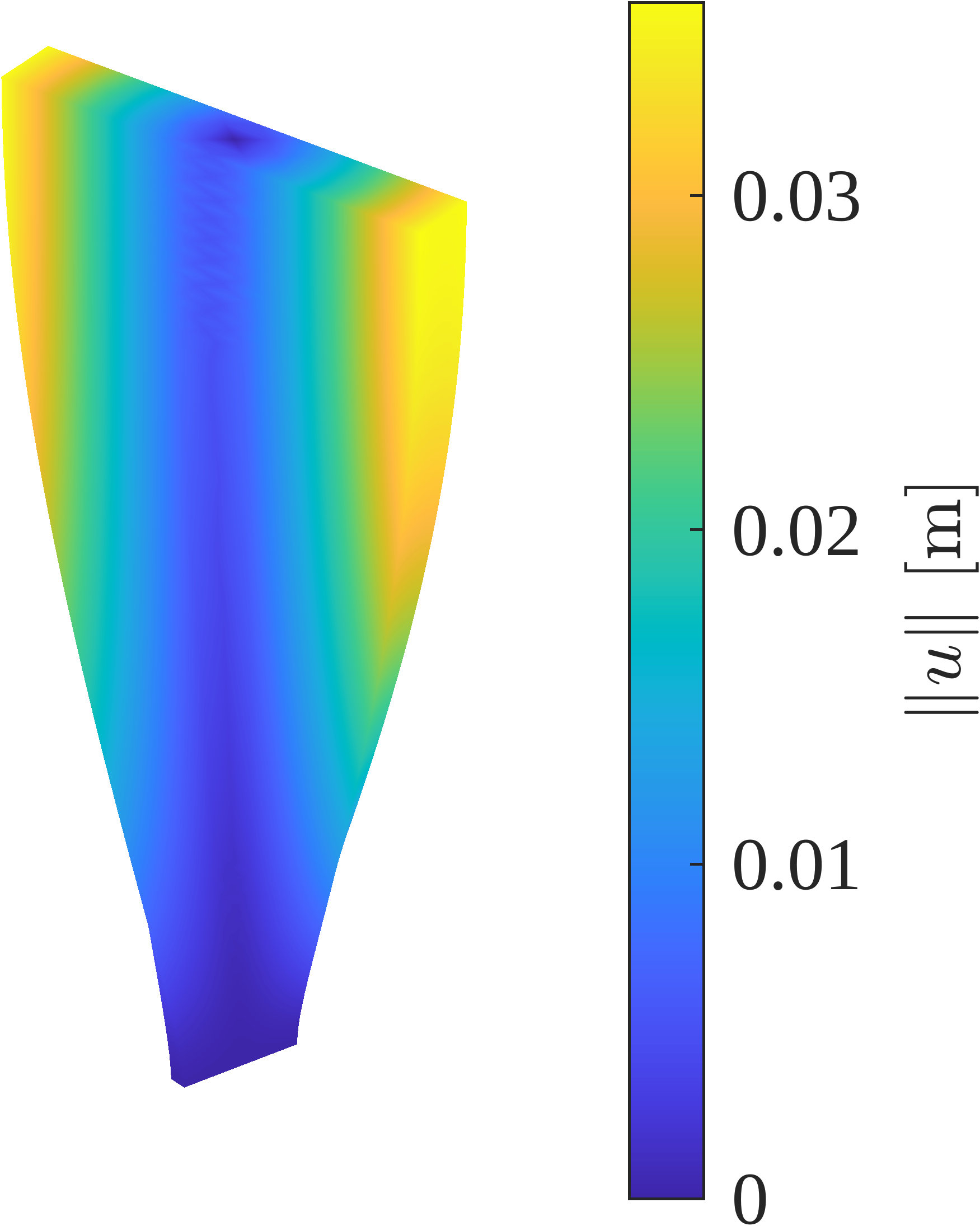}
        \caption{Fourth mode, \\$f_4 = 1576.55$Hz}
    \end{subfigure}
    \begin{subfigure}{0.33\textwidth}
        \centering
        \includegraphics[scale=0.5]{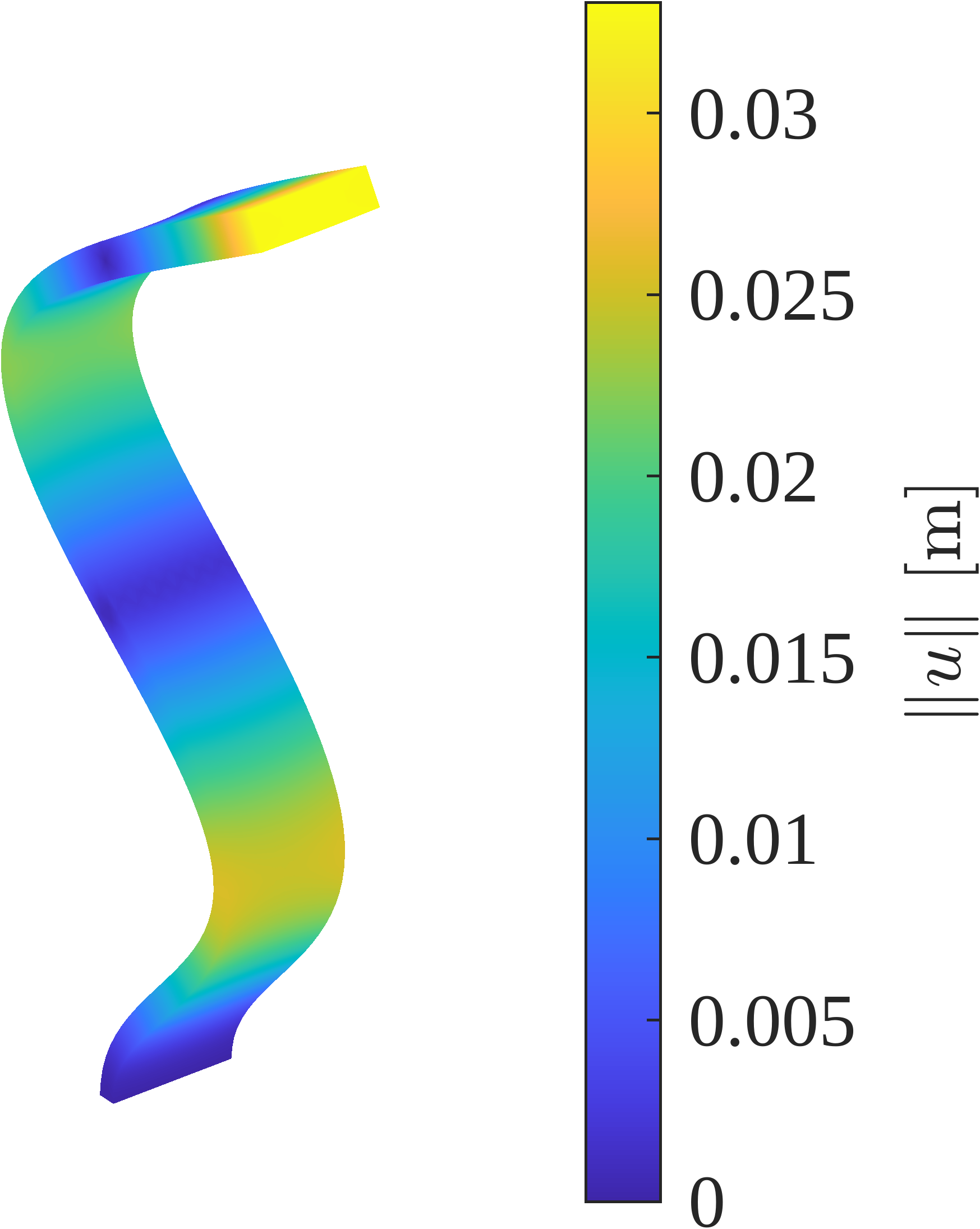}
        \caption{Fifth mode,\\$f_5=1966.22$ Hz}
    \end{subfigure}
    \caption{Eigenvalue problem in a rectangular beam.}
    \label{fig:eigenvalue_problem}
\end{figure}

We consider a rectangular beam (see Figure \ref{fig:geom}) with a length $L=0.2$ m, a width $w = 0.04$ m, and thickness $e=0.002$ m. We fix one of the ends of the beam (highlighted in red). For material properties, we consider steel with density $\rho = 7,850 \,\text{kg/m}^3$, Young modulus $E = 210$ GPa, and poisson ratio $\nu = 0.3$. The Lamé parameters are determined from the relations
\[
\lambda = \frac{\nu E}{(1+\nu)(1-2\nu)},\qquad \mu = \frac{E}{2(1+\nu)}.
\]
Therefore, $\lambda \approx 64$ GPa, which is not a very large value, makes using more sophisticated elements to address the locking phenomenon unnecessary.

The eigenvalue problem becomes finding $\omega$ such that
\[
    \nabla \cdot \underline{\boldsymbol{\sigma}}(\boldsymbol{u}) = \rho\omega^2 \boldsymbol{u} .
\]

The right-hand side term adds a block-diagonal matrix with three mass matrices on each non-zero block term on the discretized problem. As usual in engineering applications, we compute the natural frequencies as
\[
f = \frac{\omega}{2\pi}.
\]

The mesh consists of 1,656 tetrahedra and 3,557 nodes. The computation of the stiffness matrix is more complex, as it now requires nine calls to the function that assembles the stiffness matrix in the scalar case. The assembly of the enlarged stiffness matrix takes 0.257 seconds, while the assembly of the enlarged mass matrix takes 0.012 seconds.

We compute the five lowest natural frequencies using the Matlab function {\tt eigs}. The solutions are shown in Figures \ref{fig:eigenvalue_problem}b-f. We see two bending modes parallel to the axis with the minimum moment of inertia (i.e., the width of the beam). In the third mode, a bending mode occurs along the perpendicular axis to the previous one (parallel to the thickness of the beam). A torsional mode is observed for the fourth mode. The increase in beam size arises from the inherent limitation of linear elasticity, which assumes small deformations and does not accurately model large rotations. Finally, the fifth mode represents another bending mode parallel to the width of the beam.

\subsection{3D harmonic acoustic scattering }

Let us consider the following acoustic scattering problem by a bounded domain $\Omega$ in harmonic regimen: 
\begin{equation}
    \label{eq:3D}
\left| 
\begin{array}{l}
\Delta u  + (k n)^2 u  = 0\\
\partial_{\widehat{x}} u_{\rm sc}-{\rm i}k u_{\rm sc} = o(|{\boldsymbol{x}|^{-1}})
\end{array}
\right.,\,
u_{\rm sc} :=u-u_{\rm inc},\quad u_{\rm inc}= \exp({\rm i}k {\boldsymbol{x}\cdot{\boldsymbol d}})
\end{equation}
Here $k>0$ is the wave-number, and $n = n(\boldsymbol x)>0$ is the so-called function of refractive index, with $n\equiv 1$ in the exterior $\Omega$. On the other hand, $u_{\rm inc}$ is a plane incident wave, a solution of the Helmholtz equation (with $n\equiv 1$ in the full space $\mathbb{R}^3$). Therefore, $\boldsymbol{d}$ is a unit vector indicating the direction of propagation of the incident wave, and $u_{\rm sc}$ is the so-called scattered wave.

This problem is naturally posed in the full domain, which makes the imposition of the boundary condition at infinity—that is, the radiation condition—challenging. Among several options, a very attractive choice is to combine the Finite Element Method (FEM) with Boundary Element Methods (BEMs). In BEM, the solution in the exterior of a bounded domain is expressed via an integral formula involving new unknowns defined only on the boundary. The radiation condition is then satisfied naturally. This BEM-FEM combination leverages the strengths of both methods: the flexibility of FEM for handling inhomogeneities and the capability of BEM for computing the solution sufficiently far away from the object $\Omega$, as well as for determining quantities of interest such as the far field, the asymptotic behavior of the scattered wave $u_{\rm sc}$ at infinity. 

We do not intend to provide even a brief introduction to BEM-FEM coupling here. Instead, our goal is to illustrate how GMSH can handle complex domains, demonstrate the flexibility of our code in processing the returned meshes and show how our FEM solver can be integrated with available BEM codes. For details on the BEM solver used in this experiment, definition, implementation and convergence properties (which, under ideal conditions like those presented here, are superalgebraic) we refer to \cite{GaGr:2004}. The BEM-FEM coupling strategy followed in this experiment is detailed in \cite{DoGaSa:2020}.

\begin{figure}[!htb]
\[
\includegraphics[width=0.75\textwidth]{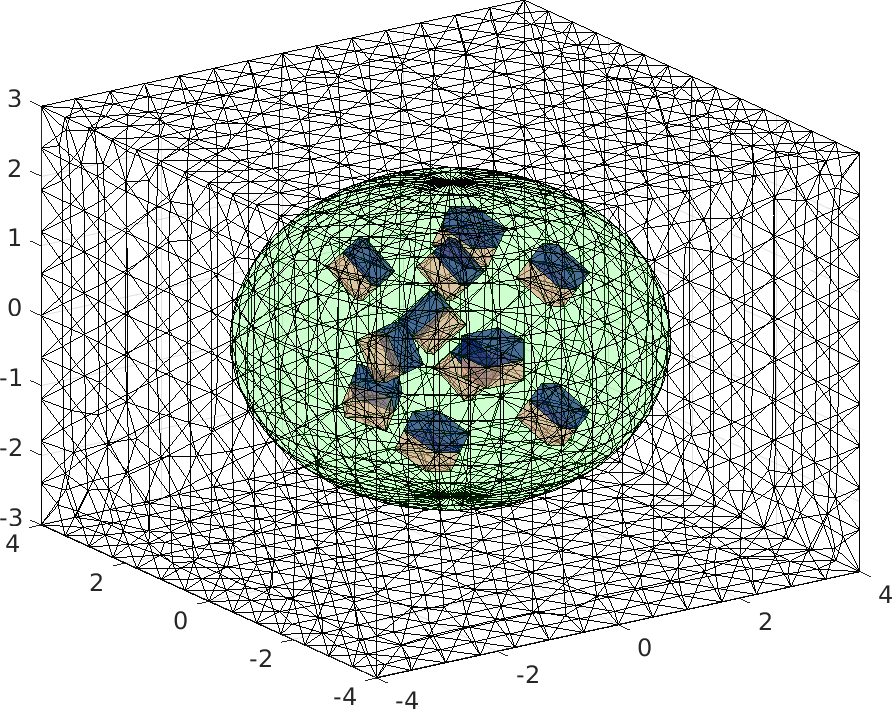}
\]
\caption{\label{fig:domain:BEM-FEM}Domain of our BEM-FEM problem. The solution is approximated by cubic finite elements (FEs) inside the prism and by the Boundary Element Method (BEM) in the exterior of the ellipsoid. There is an overlapping region where the solution is approximated simultaneously by both the FEM and the BEM.}
\end{figure}

Hence, consider the index refractive function $n$ given by
\[
 n({\bf x} = \begin{cases}
                  1.333, & \mathbf{x}\in B_{i}\\
              1.496, & \mathbf{x}\in N_{i}\\
              1,     & \text{otherwise}
             \end{cases}
\]
Here $B_i$ and $N_i$ are two halves (blue and red) of each of eleven cubes (so $i=1,2,\ldots,11$) of edge 0.3. This configuration represents a simple model of what is known in the literature as a Janus particle: particles formed by two substances of the same shape glued together. The choice of these parameters is not accidental: 1.333 and 1.496 are, respectively, the refractive indices of water and toluene at 20$^\circ$C, which are commonly used in such particles. A schematic of the domains as well as an initial mesh (only on the boundary of the FEM domain) is sketched in Figure \ref{fig:domain:BEM-FEM}. For the incident wave we took $\boldsymbol{d}=(1,0,0)$ and $k=1$. 

We used quartic elements (denoted as ${\cal P}_h^4$ in our notation) for the FEM part and compared the solutions obtained at two {\em consecutive} levels of refinement to assess the quality of our numerical solution. To clarify, by "consecutive refined meshes," we mean the following:

In 3D, unlike in 2D, it is impossible to split a tetrahedron into eight or any arbitrary number of similar elements. Instead, in Gmsh, 3D mesh generators first create a mesh for the (2D) boundaries, which serves as the foundation for generating the corresponding 3D mesh. The process starts with the previous boundary (triangular) mesh for a refined mesh, refines it uniformly (using RGB refinement), and then constructs the next-level 3D mesh based on this refined boundary.

The number of elements is approximately 10,700 and 76,000 for the coarse and fine meshes, respectively, corresponding to around 124,000 and 845,000 nodes. The problem is large enough to require a machine with substantial memory capacity since the linear system in the FEM part is solved using a direct method, specifically $L D L^\top$ factorization. The effect of the BEM scheme on convergence is negligible for two reasons: the fast convergence of the BEM scheme itself and the sufficiently large number of degrees of freedom.

\begin{figure}[!htb]
\includegraphics[width=0.51\textwidth]{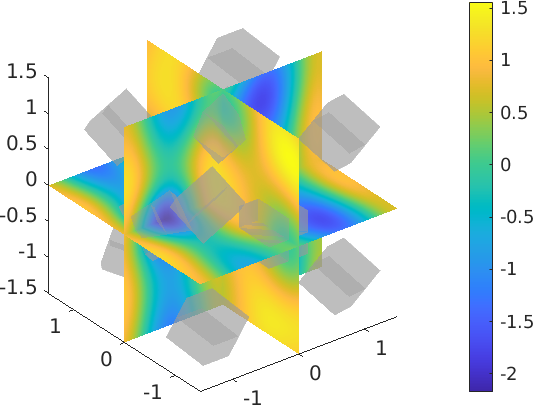}\ 
\includegraphics[width=0.51\textwidth]{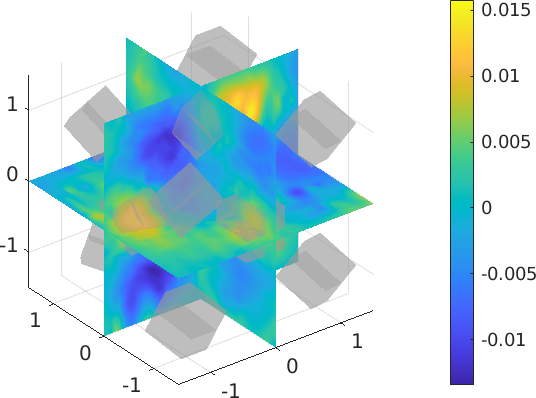}
\caption{\label{fig:BEMFEMSolution} Real part of quartic FE solution (${\cal P}_h^4$) and estimated error.}
\end{figure}

 The (real part of the) solution on the coarse mesh and the discrepancy between the FE solutions on the coarse and fine meshes is plotted on three conveniently chosen perpendicular planes, as shown in Figure \ref{fig:BEMFEMSolution}. This visualization also demonstrates the evaluation and postprocessing of the solution, as discussed in subsection \ref{sub:postprocess}.

\section*{Acknowledgments}
The first author thanks the support of projects “Adquisición de conocimiento y minería de datos, funciones especiales y métodos numéricos avanzados” from Universidad Pública de Navarra, Spain and “Técnicas innovadoras para la resolución de problemas evolutivos”, ref. PID2022-136441NB-I00 from Ministerio de Ciencia e Innovación, Gobierno de España, Spain. The second author thanks the support of the Basque Center of Applied Mathematics (BCAM) and the support of the project ``Proyecto I+D+i 2019: Inversión de datos en tiempo real usando métodos de aprendizaje profundo (DEEPINVERSE)'', ref. PID2019-108111RB-I00 from Ministerio de Ciencia e Innovación, Gobierno de España, Spain.

%
%
\bibliographystyle{plain}
\bibliography{biblioGMSH}

\begin{thebibliography}{10}

\bibitem{MR1709562}
J.~Alberty, C.~Carstensen, and S.A. Funken.
\newblock Remarks around 50 lines of {M}atlab: short finite element
  implementation.
\newblock {\em Numer. Algorithms}, 20(2-3):117--137, 1999.

\bibitem{AlnaesEtal2015}
{M}.S. Alnaes, J.~Blechta, J.~Hake, A.~Johansson, B.~Kehlet, A.~Logg, C.N.
  Richardson, J.~Ring, M.E. Rognes, and G.N. Wells.
\newblock The {FEniCS} project version 1.5.
\newblock {\em Archive of Numerical Software}, 3, 2015.

\bibitem{MR2194203}
C.~Bahriawati and C.~Carstensen.
\newblock Three {MATLAB} implementations of the lowest-order {R}aviart-{T}homas
  {MFEM} with a posteriori error control.
\newblock {\em Comput. Methods Appl. Math.}, 5(4):333--361, 2005.

\bibitem{brenner_scott}
S.C. Brenner and L.R. Scott.
\newblock {\em The Mathematical Theory of Finite Element Methods}, volume~15 of
  {\em Texts in Applied Mathematics}.
\newblock Springer, 3rd edition, 2008.
\newblock A comprehensive guide on the mathematical foundations of the finite
  element method.

\bibitem{MR1935965}
C.~Carstensen and R.~Klose.
\newblock Elastoviscoplastic finite element analysis in 100 lines of {M}atlab.
\newblock {\em J. Numer. Math.}, 10(3):157--192, 2002.

\bibitem{Ciarlet1978}
P.G. Ciarlet.
\newblock {\em The Finite Element Method for Elliptic Problems}, volume~4 of
  {\em Studies in Mathematics and its Applications}.
\newblock North-Holland, Amsterdam, 1978.
\newblock A classic foundational text in finite element theory.

\bibitem{freecad}
{FreeCAD} Community.
\newblock {FreeCAD:} an open-source parametric 3d cad modeler, 2025.

\bibitem{DoGaSa:2020}
V.~Dom{\'i}nguez, M.~Ganesh, and F.~J. Sayas.
\newblock An overlapping decomposition framework for wave propagation in
  heterogeneous and unbounded media: formulation, analysis, algorithm, and
  simulation.
\newblock {\em J. Comput. Phys.}, 403:109052, 20, 2020.

\bibitem{fem3d_gmsh}
{A}. Duque-Salazar.
\newblock {FEM3D\_GMSH: Finite Element Method with 3D meshes generated by
  GMSH}.
\newblock \url{https://github.com/aleduques/FEM3D\_GMSH}, 2023.

\bibitem{duque2023integration}
{A}. Duque-Salazar.
\newblock Integration of the meshing tool gmsh with matlab/octave for the
  resolution of 3d boundary value problems with simplicial finite elements,
  2023.

\bibitem{MR2050138}
A.~Ern and J.-L. Guermond.
\newblock {\em Theory and practice of finite elements}, volume 159 of {\em
  Applied Mathematical Sciences}.
\newblock Springer-Verlag, New York, 2004.

\bibitem{MR4242224}
A.~Ern and J.-L. Guermond.
\newblock {\em Finite elements {I}---{A}pproximation and interpolation},
  volume~72 of {\em Texts in Applied Mathematics}.
\newblock Springer, Cham, [2021] \copyright 2021.

\bibitem{MR4269305}
A.~Ern and J.-L. Guermond.
\newblock {\em Finite elements {II}---{G}alerkin approximation, elliptic and
  mixed {PDE}s}, volume~73 of {\em Texts in Applied Mathematics}.
\newblock Springer, Cham, [2021] \copyright 2021.

\bibitem{teampancho}
F.-J.~Sayas et~collaborators.
\newblock Team pancho website, 2025.
\newblock Accessed: 2025-01-17.

\bibitem{funken2011fem}
S.A. Funken, D.~Praetorius, and L.K. R{\"u}de.
\newblock Efficient implementation of adaptive {P}1-{FEM} in {M}atlab.
\newblock {\em Comput. Methods Appl. Math.}, 11(4):460--490, 2011.

\bibitem{GaGr:2004}
M.~Ganesh and I.~G. Graham.
\newblock A high-order algorithm for obstacle scattering in three dimensions.
\newblock {\em J. Comput. Phys.}, 198(1):211--242, 2004.

\bibitem{geuzaine2009GMSH}
C.~Geuzaine and J.-F. Remacle.
\newblock Gmsh: A 3-d finite element mesh generator with built-in pre- and
  post-processing facilities.
\newblock {\em International Journal for Numerical Methods in Engineering},
  79(11):1309--1331, 2009.
\newblock An open-source tool for finite element mesh generation and
  post-processing.

\bibitem{MR3043640}
F.~Hecht.
\newblock New development in {FreeFem++}.
\newblock {\em J. Numer. Math.}, 20(3-4):251--265, 2012.

\bibitem{sayas2015fem}
F.-J. Sayas.
\newblock An introduction to the finite element method, 2015.
\newblock Accessed: 2024-11-28.

\bibitem{sayas_fem_tools_3d}
F.-J. Sayas and collaborators.
\newblock {FEM} tools in three dimensions.
\newblock Unpublished document, personal communication, n.d.

\bibitem{doi:https://doi.org/10.1002/9781119176817.ecm2008}
{E}.P. Stephan.
\newblock {\em Coupling of Boundary Element Methods and Finite Element
  Methods}, pages 1--40.
\newblock John Wiley \& Sons, Ltd, 2017.

\bibitem{MR3674245}
{O}.J. Sutton.
\newblock The virtual element method in 50 lines of {MATLAB}.
\newblock {\em Numer. Algorithms}, 75(4):1141--1159, 2017.

\bibitem{zienkiewicz_fem}
O.C. Zienkiewicz and R.L. Taylor.
\newblock {\em The Finite Element Method: Its Basis and Fundamentals}.
\newblock Elsevier, 7th edition, 2013.
\newblock A comprehensive introduction to FEM with a focus on fundamentals and
  applications.

\end{thebibliography}

\end{document}